\documentclass{amsart}
\usepackage{graphicx,amsfonts}

\newcommand{\ind}{\mbox{\mbox{ } \mbox{ } \mbox{ }}}

\newcommand{\be}{\begin{enumerate}}
\newcommand{\ee}{\end{enumerate}}
\newcommand{\bit}{\begin{itemize}}
\newcommand{\eit}{\end{itemize}}
\newcommand{\bc}{\begin{center}}
\newcommand{\ec}{\end{center}}
\newcommand{\beq}{\begin{equation}}
\newcommand{\eeq}{\end{equation}}
\newcommand{\bqn}{\begin{eqnarray}}
\newcommand{\eqn}{\end{eqnarray}}
\newcommand{\bqns}{\begin{eqnarray*}}
\newcommand{\eqns}{\end{eqnarray*}}
\newcommand{\barr}{\begin{array}}
\newcommand{\earr}{\end{array}}
\newcommand{\N}{{\bf N}}

\newcommand{\R}{{\bf R}}

\newcommand{\bd}{\mbox{$\partial$}}

\newcommand{\rarrow}{\rightarrow}

\newcommand{\sC}{{\mathcal C}}

\newcommand{\sE}{{\mathcal E}}
\newcommand{\sL}{{\mathcal L}}

\newcommand{\se}{{\mathcal E}}

\newcommand{\bdoc}{\begin{document}}
\newcommand{\edoc}{\end{document}}
\newcommand{\docart}{\documentclass[12pt]{article}}
\newcommand{\vs}{\vspace{.16in}}
\newcommand{\half}{\frac{1}{2}}
\newcommand{\third}{\frac{1}{3}}
\newcommand{\fourth}{\frac{1}{4}}
\newcommand{\diffeo}{\mbox{ diffeomorphism }} 
\newcommand{\dnorm}[1]{\mbox{$\mid \mid #1 \mid \mid$}}
\newcommand{\dstyle}[1]{\mbox{$\displaystyle{ #1 }$}}
\newcommand{\zero}{ \mbox{$ {\bf 0 } $}}
\newcommand{\PL}{\mbox{$Per_L(\R)$}}
\newcommand{\htop}{\mbox{$h_{top}$}}

\newcommand{\bfE}{\mbox{${\bf E}$}}
\newcommand{\bfF}{\mbox{${\bf F}$}}

\newcommand{\hrel}{\mbox{ ${\simeq_{\bd}}$ }}

\newcommand{\0}{{\bf 0}}

\newcommand{\twoover}[2]{{\footnotesize \begin{array}[t]{c}\mbox{ {\normalsize #1}}  \\
                               #2 
        \end{array} } \hspace{1cm} }
\newcommand{\threeover}[3]{ {\footnotesize \begin{array}[t]{c}\mbox{ {\normalsize #1}}  \\
                               #2  \\
                                 #3
        \end{array} } \hspace{1cm} }
\newcommand{\fourover}[4]{ {\footnotesize \begin{array}[t]{c}\mbox{ {\normalsize #1}}  \\
                               #2  \\
                                 #3 \\
				#4
        \end{array} } \hspace{1cm}  }

\newenvironment{vover3}[3]{   
       {\footnotesize \begin{array}[t]{c}
                               {\normalsize #1} \vspace{.2cm} \\ 
                                 #2 \\
                                 #3  
                       \end{array} }
                            } { \hspace{1cm} } 
\newenvironment{vover4}[4] {
       {\footnotesize \begin{array}[t]{c} 
                               {\normalsize #1} \vspace{.2cm} \\  
                                 #2 \\
                                 #3 \\
                                 #4  
                       \end{array} }
                            } {  \hspace{1cm} }

\newenvironment{seteqover}[3]{#1 \ \ = \hspace{.6cm}
                                #2 #3} {\ind} 

\newcommand{\norm}[1]{\mbox{$\mid #1 \mid$}}
\newcommand{\Norm}[1]{\mbox{$\parallel #1 \parallel$}}
\newcommand{\imp}{\mbox{$\Longrightarrow$}}
\newcommand{\pr}{\mbox{$\prime$}}

\newcommand{\twocvec}[2]{\mbox{$\left( \begin{array}{c}
				 #1 \\
				 #2
				\end{array} \right)$} }

\newcommand{\twomat}[4]{\mbox{$ \left( \begin{array}{cc}
                              #1 & #2 \\
                           #3 & #4
                       \end{array} \right) $ } }

\newcommand{\Itwomat}[4]{\mbox{$ \left( \begin{array}{cc}
                              #1 & #2 \\
                           #3 & #4
                       \end{array} \right)^{-1} $ } }

\newcommand{\threemat}[9]{ \[ \left( \begin{array}{ccc}
				 #1 & #2 & #3 \\
				 #4 & #5 & #6 \\
				 #7 & #8 & #9 
				\end{array} \right) \]  }

\newcommand{\newthreemat}[9]{ \left( \begin{array}{ccc}
				 #1 & #2 & #3 \\
				 #4 & #5 & #6 \\
				 #7 & #8 & #9 
				\end{array} \right) }

\newcommand{\Ithreemat}[9]{\mbox{ $\left( \begin{array}{ccc}
				 #1 & #2 & #3 \\
				 #4 & #5 & #6 \\
				 #7 & #8 & #9
				\end{array} \right)^{-1} $ } }

\newcommand{\ba}{\mbox{${\bf a}$}}
\newcommand{\bfa}{\mbox{${\bf a}$}}
\newcommand{\bb}{\mbox{${\bf b}$}}
\newcommand{\bfb}{\mbox{${\bf b}$}}
\newcommand{\bfc}{\mbox{${\bf c}$}}
\newcommand{\bi}{\mbox{${\bf i}$}}
\newcommand{\bj}{\mbox{${\bf j}$}}
\newcommand{\bfi}{\mbox{${\bf i}$}}
\newcommand{\bfj}{\mbox{${\bf j}$}}
\newcommand{\bfk}{\mbox{${\bf k}$}}
\newcommand{\bfx}{\mbox{${\bf x}$}}
\newcommand{\bfy}{\mbox{${\bf y}$}}
\newcommand{\bfz}{\mbox{${\bf z}$}}
\newcommand{\bfu}{\mbox{${\bf u}$}}
\newcommand{\bfv}{\mbox{${\bf v}$}}
\newcommand{\bfw}{\mbox{${\bf w}$}}
\newcommand{\bfg}{\mbox{${\bf g}$}}
\newcommand{\onef}{\mbox{${\bf 1}$}}
\renewcommand{\bfc}{\mbox{${\bf c}$}}
\newcommand{\bfe}{\mbox{${\bf e}$}}
\newcommand{\bff}{\mbox{${\bf f}$}}
\newcommand{\bS}{\mbox{${\bf S}$}}
\newcommand{\bV}{\mbox{${\bf V}$}}
\newcommand{\bw}{\bar{w}} 
\newcommand{\fh}{\hat{f}}
\newcommand{\wh}{\hat{w}}
\newcommand{\Wh}{\hat{W}}
\newcommand{\Ih}{\hat{I}}
\newcommand{\Lh}{\hat{L}}
\newcommand{\grgh}{\hat{\gamma}}
\newcommand{\br}{\mbox{${\bf r}$}}
\newcommand{\bs}{\mbox{${\bf s}$}}

\newcommand{\oz}{\overline{z}}
\newcommand{\ogz}{\overline{gz}}
\newcommand{\oxi}{\overline{\xi}}

\newcommand{\ninf}{n \rarrow \infty}
\newcommand{\linf}{\lim_{n \rarrow \infty}}
\newcommand{\fon}{ \frac{1}{n} }
\newcommand{\skn}{ \sum_{k=0}^{n-1} }
\newcommand{\siin}{ \sum_{i=0}^{n-1} }

\newcommand{\alp}{\mbox{$\alpha$}}
\newcommand{\gra}{\mbox{$\alpha$}}

\newcommand{\bet}{\mbox{$\beta$}}
\newcommand{\grb}{\mbox{$\beta$}}

\newcommand{\gam}{\mbox{$\gamma$}}
\newcommand{\Gam}{\mbox{$\Gamma$}}
\newcommand{\Del}{\mbox{$\Delta$}}
\newcommand{\sig}{\mbox{$\sigma$}}
\newcommand{\Sig}{\mbox{$\Sigma$}}

\newcommand{\grg}{\mbox{$\gamma$}}
\newcommand{\Grg}{\mbox{$\Gamma$}}
\newcommand{\grd}{\mbox{$\delta$}}
\newcommand{\Grd}{\mbox{$\Delta$}}
\newcommand{\grs}{\mbox{$\sigma$}}
\newcommand{\Grs}{\mbox{$\Sigma$}}
\newcommand{\grf}{\mbox{$\phi$}}
\newcommand{\Grf}{\mbox{$\Phi$}}
\newcommand{\grth}{\mbox{$\theta$}}
\newcommand{\Grth}{\mbox{$\Theta$}}
\newcommand{\gro}{\mbox{$\omega$}}
\newcommand{\Gro}{\mbox{$\Omega$}}
\newcommand{\gre}{\mbox{$\epsilon$}}
\newcommand{\bgre}{\mbox{$\bar{\epsilon}$}}
\newcommand{\grve}{\mbox{$\varepsilon$}}
\newcommand{\grr}{\mbox{$\rho$}}
\newcommand{\Grr}{\mbox{$\Rho$}}
\newcommand{\grt}{\mbox{$\tau$}}
\newcommand{\gri}{\mbox{$\iota$}}
\newcommand{\grp}{\mbox{$\pi$}}
\newcommand{\grP}{\mbox{$\Pi$}}
\newcommand{\grl}{\mbox{$\lambda$}}
\newcommand{\grlh}{\mbox{$\hat{\lambda}$}}
\newcommand{\grL}{\mbox{$\Lambda$}}
\newcommand{\Grl}{\mbox{$\Lambda$}}
\newcommand{\grz}{\mbox{$\zeta$}}

\newcommand{\thet}{\mbox{$\theta$}}
\newcommand{\Thet}{\mbox{$\Theta$}}
\newcommand{\ome}{\mbox{$\omega$}}
\newcommand{\Ome}{\mbox{$\Omega$}}
\newcommand{\eps}{\mbox{$\epsilon$}}

\newcommand{\iot}{\mbox{$\iota$}}

\newcommand{\lam}{\mbox{$\lambda$}}
\newcommand{\Lam}{\mbox{$\Lambda$}}
\newcommand{\zet}{\mbox{$\zeta$}}

\newcommand{\kap}{\mbox{$\kappa$}}

\newcommand{\grx}{\mbox{$\xi$}}
\newcommand{\Grx}{\mbox{$\Xi$}}
\newcommand{\grc}{\mbox{$\chi$}}
\newcommand{\Grc}{\mbox{$\Chi$}}
\newcommand{\grn}{\mbox{$\nu$}}
\newcommand{\grm}{\mbox{$\mu$}}
\newcommand{\grk}{\mbox{$\kappa$}}

\newtheorem{theorem}{Theorem}[section]
\newtheorem{lemma}[theorem]{Lemma}

\newtheorem{definition}[theorem]{Definition}
\newtheorem{example}[theorem]{Example}
\newtheorem{xca}[theorem]{Exercise}

\newtheorem{remark}[theorem]{Remark}

\newtheorem{Cor}[theorem]{Corollary}
\newtheorem{Corollary}[theorem]{Corollary}
\newtheorem{Prop}[theorem]{Proposition}
\newtheorem{Proposition}[theorem]{Proposition}
\newtheorem{Lem}[theorem]{Lemma}
\newtheorem{Lemma}[theorem]{Lemma}
\newtheorem{Rem}[theorem]{Remark}
\newtheorem{Remark}[theorem]{Remark}
\newtheorem{Conjecture}[theorem]{Conjecture}
\newtheorem{Theorem}[theorem]{Theorem}
\newtheorem{Th}[theorem]{Theorem}
\newtheorem{Fact}[theorem]{Fact}
\newtheorem{Cond-(*)}[theorem]{(*)}
\newtheorem{Cond-(**)}[theorem]{(**)}
\newtheorem{Cond_NoLabel}[theorem]{}
\newtheorem{Cnol}[theorem]{}
\newtheorem{cnol}[theorem]{}
\newtheorem{Def}[theorem]{Definition}
\newtheorem{Alg}[theorem]{Algorithm}

\newcommand{\btheor}{\begin{theorem}}
\newcommand{\etheor}{\end{theorem}}
\newcommand{\bProp}{\begin{Prop}}
\newcommand{\eProp}{\end{Prop}}
\newcommand{\bCor}{\begin{Cor}}
\newcommand{\eCor}{\end{Cor}}
\newcommand{\bLem}{\begin{Lem}}
\newcommand{\eLem}{\end{Lem}}

\newcommand{\ta}{\tilde{a}}
\newcommand{\tgra}{\tilde{\alpha}}
\newcommand{\tGrl}{\tilde{\Grl}}
\newcommand{\tb}{\tilde{b}}
\newcommand{\tgrb}{\tilde{\beta}}
\newcommand{\tgrg}{\tilde{\gamma}}
\newcommand{\tc}{\tilde{c}}
\newcommand{\td}{\tilde{d}}
\newcommand{\te}{\tilde{e}}
\newcommand{\tf}{\tilde{f}}
\newcommand{\teta}{\tilde{\eta}}
\newcommand{\tg}{\tilde{g}}
\newcommand{\tih}{\tilde{h}}
\renewcommand{\th}{\tilde{h}}
\newcommand{\ti}{\tilde{i}}
\newcommand{\tj}{\tilde{j}}
\newcommand{\tk}{\tilde{k}}
\newcommand{\tl}{\tilde{l}}
\newcommand{\tm}{\tilde{m}}
\newcommand{\tn}{\tilde{n}}
\newcommand{\tio}{\tilde{o}}
\newcommand{\tp}{\tilde{p}}
\newcommand{\tq}{\tilde{q}}
\newcommand{\tr}{\tilde{r}}
\newcommand{\ts}{\tilde{s}}
\newcommand{\tit}{\tilde{t}}
\newcommand{\tu}{\tilde{u}}
\newcommand{\tv}{\tilde{v}}
\newcommand{\tw}{\tilde{w}}
\newcommand{\tx}{\tilde{x}}
\newcommand{\ty}{\tilde{y}}
\newcommand{\tz}{\tilde{z}}
\newcommand{\tA}{\tilde{A}}
\newcommand{\tB}{\tilde{B}}
\newcommand{\tC}{\tilde{C}}
\newcommand{\tD}{\tilde{D}}
\newcommand{\tE}{\tilde{E}}
\newcommand{\tF}{\tilde{F}}
\newcommand{\tG}{\tilde{G}}
\newcommand{\tH}{\tilde{H}}
\newcommand{\tI}{\tilde{I}}
\newcommand{\tJ}{\tilde{J}}
\newcommand{\tK}{\tilde{K}}
\newcommand{\tL}{\tilde{L}}
\newcommand{\tM}{\tilde{M}}
\newcommand{\tN}{\tilde{N}}
\newcommand{\tO}{\tilde{O}}
\newcommand{\tP}{\tilde{P}}
\newcommand{\tQ}{\tilde{Q}}
\newcommand{\tR}{\tilde{R}}
\newcommand{\tS}{\tilde{S}}
\newcommand{\tT}{\tilde{T}}
\newcommand{\tU}{\tilde{U}}
\newcommand{\tV}{\tilde{V}}
\newcommand{\tW}{\tilde{W}}
\newcommand{\tX}{\tilde{X}}
\newcommand{\tY}{\tilde{Y}}
\newcommand{\tZ}{\tilde{Z}}

\newcommand{\tTh}{\tilde{\Theta}}

\newcommand{\fl}{f_{\grl}}
\newcommand{\el}{E_{\grl}}
\newcommand{\sel}{\se_{\grl}}
\newcommand{\szl}{\sz_{\grl}}
\newcommand{\zl}{Z_{\grl}}
\newcommand{\famfl}{\{ \fl \} }
\newcommand{\sT}{{\mathcal T}}
\newcommand{\hsT}{{\hat{\mathcal T}}}
\newcommand{\hsH}{{\hat{\mathcal H}}}
\newcommand{\hsCZ}{{\hat{\mathcal CZ}}}
\newcommand{\sCG}{{\mathcal CG}}
\newcommand{\sFG}{{\mathcal FG}}
\newcommand{\sTemp}{{\mathcal TM}}
\newcommand{\ab}{{\bf a}}
\newcommand{\xb}{{\bf x}}
\newcommand{\yb}{{\bf y}}
\newcommand{\Pf}{{\bf Proof}}
\newcommand{\uG}{\bar{{\mathcal G}}}
\newcommand{\sG}{{\mathcal G}}
\newcommand{\sF}{{\mathcal F}}
\newcommand{\stG}{\tilde{{\mathcal G}}}
\newcommand{\tmu}{\tilde{\mu}}
\newcommand{\tpi}{\tilde{\pi}}
\newcommand{\tgrl}{\tilde{\grl}}
\newcommand{\sD}{{\mathcal D}}
\newcommand{\sDh}{\hat{\sD}}
\newcommand{\sP}{{\mathcal P}}
\newcommand{\tsP}{{\tilde{\mathcal P}}}
\newcommand{\sM}{{\mathcal M}}
\newcommand{\bsM}{\bar{{\mathcal M}}}
\newcommand{\sN}{{\mathcal N}}
\newcommand{\sS}{{\mathcal S}}
\newcommand{\tsS}{\tilde{{\mathcal S}}}
\newcommand{\tgrt}{\tilde{{\tau}}}
\newcommand{\sI}{{\mathcal I}}
\newcommand{\sW}{{\mathcal W}}
\newcommand{\sCZ}{{\mathcal C}{\mathcal Z}}
\newcommand{\tsT}{\tilde{\sT}}
\newcommand{\tsC}{\tilde{\sC}}
\newcommand{\tsD}{\tilde{\sD}}
\newcommand{\tgre}{\tilde{\gre}}
\newcommand{\sH}{{\mathcal H}}
\newcommand{\tsH}{\tilde{\sH}}
\newcommand{\tsG}{\tilde{\sG}}
\newcommand{\tsW}{\tilde{\sW}}
\newcommand{\tDel}{\tilde{\Grd}}
\newcommand{\sA}{{\mathcal A}}
\newcommand{\sB}{{\mathcal B}}
\newcommand{\sCF}{{\mathcal CF}}
\newcommand{\bsH}{\bar{\sH}}
\newcommand{\ntwob}{\bar{\nu_2}}
\newcommand{\nb}{\bar{\nu}}
\newcommand{\bE}{\bar{E}}
\newcommand{\bF}{\bar{F}}
\newcommand{\bG}{\bar{G}}
\newcommand{\bI}{\bar{I}}
\newcommand{\bZ}{\bar{Z}}
\newcommand{\bh}{\bar{h}}
\newcommand{\bm}{\bar{m}}
\newcommand{\ml}{m_{\ell}}
\newcommand{\bT}{\bar{T}}
\newcommand{\bK}{\bar{K}}
\newcommand{\bx}{\bar{x}}
\newcommand{\bB}{\bar{B}}
\newcommand{\bxi}{\bar{\xi}}
\newcommand{\bu}{\bar{u}}
\newcommand{\bv}{\bar{v}}
\newcommand{\bg}{\bar{g}}
\newcommand{\bN}{\bar{\N}}
\newcommand{\bX}{\bar{X}}
\newcommand{\bY}{\bar{Y}}
\newcommand{\bGrs}{\bar{\Grs}}
\newcommand{\bgrs}{\bar{\grs}}
\newcommand{\bmu}{\bar{\mu}}
\newcommand{\bnu}{\bar{\nu}}
\newcommand{\bTh}{\bar{\Theta}}
\newcommand{\breta}{\bar{\eta}}
\newcommand{\db}{\bar{d}}
\newcommand{\dx}{\dot{x}}
\newcommand{\dy}{\dot{y}}
\newcommand{\dz}{\dot{z}}
\newcommand{\ddx}{\partial_x} 
\newcommand{\ddy}{\partial_y} 
\newcommand{\ddz}{\partial_z} 
\newcommand{\dxn}{\dot{x}_n}
\newcommand{\n}{\norm}
\newcommand{\No}{\dnorm}
\newcommand{\no}[1]{\mbox{$\norm{#1}_1$}}
\newcommand{\grfh}{\hat{\grf}}

\newcommand{\ol}{\overline}
\newcommand{\ob}{\bar}

\newcommand{\eqdef}{\mbox{ $\stackrel{{\rm def}}{=}$ }}
\newcommand{\defeq}{\mbox{ $\stackrel{{\rm def}}{=}$ }}

\newcommand{\starline}
 { \vspace{0.12in} 
   \noindent
 /*********************************************************************/
   \vspace{0.12in}
 }

\newcommand{\Lpfk}{\mbox{$\cap \hspace{-.44cm}\mid$}}

\dedicatory{Dedicated to the memory of Dmitri Anosov and Philip Hartman} 

\title{On a Differentiable  Linearization Theorem of Philip Hartman}

\author{Sheldon E. Newhouse}

\date{\today}
         
\begin{document}

\bibliographystyle{plain}

\maketitle

\begin{abstract}
Given a linear automorphism $L$ of a Banach space $E$, let $\rho(L)$ 
denote its spectral radius and let $c(L) = \rho(L)\rho(L^{-1})$ denote
its spectral condition number. Given the direct sum decompostion $L =
A \oplus B$ let $c_h = \max(c(A), c(B)), \rho_h =
\max(\rho(A^{-1}),\rho(B))$. For  $0 < \alpha < 1$, the map $L$ is
called $\alpha$-hyperbolic if $L$ can be written as $L = A \oplus B$
so that $c_h \rho_h^{\alpha}$ is less than one.
A fixed point of a $C^1$ diffeomorphism $T$ is called
$\alpha-hyperbolic$  if its derivative $DT(p)$ is $\alpha-$hyperbolic.
A well-known theorem of Philip Hartman states that if $E$ is finite
dimensional with an $\alpha-$hyperbolic fixed point and, in addition,  
 the derivative $DT$ is uniformly Lipschitz near
$p$, then the map $T$ is locally $C^{1,\beta}$ linearizable near $p$ for some $\beta > 0$. We obtain the
same result under the weaker assumption that $DT$ is uniformly
$\alpha-$H{\"{o}}lder near $p$. We also  extend
the result to Banach spaces with $C^{1,\alpha}$ bump functions and obtain
continuous  dependence of the linearization on
parameters. The results apply to give simpler proofs under weaker
regularity assumptions of classical theorems of 
L. P. Shilnikov and others on  the existence of horseshoe type dynamics
and bifurcations near homoclinic curves. 
\end{abstract} 

\section{Introduction}

Linearization theorems are of fundamental
 importance in Dynamical
Systems.  On the one hand they provide simple descriptions of the
behavior near critical points and periodic motions, and on the other,
they can often be applied with non-local techniques to give
substantial information about global structures, e.g. as in
\cite{Guckenheimer-Holmes},\cite{Ilyashenko-Li}, 
\cite{Shil-Shil-Tur-Chua-I}, \cite{Shil-Shil-Tur-Chua-II}, \cite{Palis-Takens}, \cite{Chow-Hale-Bifurcation-Theory}.

The so-called {\it Grobman-Hartman Theorem} states that a $C^r$ (with 
$r$ a positive integer) diffeomorphism or flow can be (locally) $C^0$ linearized near a hyperbolic
fixed point\footnote{This theorem, first proved in Euclidean space
  independently by Hartman and Grobman, was extended to Banach spaces,
  independently by Palis \cite{Palis-on-Hartman} and Pugh \cite{Pugh-on-Hartman}.}.  While this
theorem gives topological information about  
orbits which remain near the fixed
point, it is  inadequate for the study of orbits which recur near
the fixed point after traveling a relatively large distance away from
it.  The analysis of such orbits often requires estimates of
derivatives, and hence, makes good use of linearizations with various
amounts of smoothness.  

Standard smooth (local) linearization theorems near a hyperbolic fixed
point in a finite dimensional manifold have the following form.  First,
by a suitable choice of coordinates,
we may assume that the fixed point is the origin in Euclidean space.
Next, given a desired smoothness $k$ for the
linearization, there is an integer $N(k)$ such that if the $N(k)-$jet
of the system at $0$ is linear (the non-linear
terms up to order $N(k)$ vanish), then there is a linearization of
order $k$ near the fixed point.  This vanishing of the nonlinear terms
is most often guaranteed via a preliminary change of coordinates under
so-called {\it diophantine inequalities} or {\it non-resonance conditions}.  These are polynomial inequalities involving the
eigenvalues of the derivative of the map or flow at the fixed
point.  See Hartman \cite{Hartman-ODE}, Sternberg \cite{Sternberg-Loc-Con-Real-Line},
\cite{Sternberg-Loc-Con-Theorem-Poincare}, Bronstein and Kopanskii \cite{Bronstein-Kopanskii}, and the
references contained therein for more details.  Several papers of 
Chaperon \cite{Chaperon-Steklov}, \cite{Chaperon-Stable-Man-Perron-Irwin}
describe a beautiful development of Invariant Manifold Theory and its applications to linearization
theorems.  

In this paper we deal  with local linearizations under rather weak smoothness
assumptions and, in the bi-circular case defined below,  {\it without
diophantine inequalities}. 

 Since the results and arguments hold in Banach spaces with
sufficiently smooth bump functions, let
us set up the notation in that case. 

Let $E$ and $F$ be real Banach spaces and let $L(E,F)$ denote the
Banach space of bounded linear maps from $E$ to $F$ with the usual
supremum norm, for $L \in L(E,F)$, 

\[ \n{L} = \sup_{\mid x \mid=1} \n{Lx}. \]

Let $Aut(E)$ denote the linear automorphisms of $E$; i.e., the
invertible elements of $L(E,E)$.

Let $U$ be  a non-empty open, connected  subset of $E$ and let $\alpha > 0$ be a
positive real number\footnote{In this paper, unless otherwise stated,
  all open sets will be assumed to be connected.}. 

The map $f: U \rarrow F$ is called $\alpha-$H{\"{o}}lder continuous at
a point $x \in U$ (or simply $\alpha-$H{\"{o}}lder at $x$) if 

\beq \label{pt-Holder}
Hol(f,x,U) \eqdef \sup_{y \in U, y \neq x}\frac{\n{f(y) - f(x)}}{\n{y-x}^{\alpha}} <
\infty. 
\eeq

The map $f$ is called {\it $\alpha-$H{\"{o}}lder} in $U$ if 

\beq \label{unif-Holder}
Hol(f,U) \eqdef \sup_{x,y \in U, y \neq x}\frac{\n{f(y) - f(x)}}{\n{y-x}^{\alpha}} <
\infty. 
\eeq

The latter condition is sometimes called {\it uniformly
$\alpha-$H{\"{o}}lder} in $U$.

If $\alpha = 1$ in the above inequalities, we say that $f$ is,
respectively, {\it Lipschitz at $x$} or {\it Lipschitz in $U$}. 

We refer to the constant $Hol(f,U)$ in (\ref{pt-Holder}) (or (\ref{unif-Holder}))
as the {\it
  H{\"{o}}lder constant} of $f$ at $x$ (in $U$).  When $\alpha = 1$, we use the term
{\it Lipschitz constant} of $f$, and we write this as $Lip(f,U)$.

 The map $f$ is {\it differentiable} at a point $x \in U$ if there is a
bounded linear map $Df(x) \in L(E,F)$ such that  

\[ \lim_{\mid h \mid  \rarrow 0} \frac{\mid f(x+h) - f(x) - Df(x)h
  \mid}{\mid h \mid } = 0. \]

Sometimes this concept is referred by saying that $f$ is {\it
Fr{\'{e}}chet differentiable} at $x$. For notational convenience, we
sometimes write $Df(x)$ as $Df_x$.

If $f$ is differentiable at each $x \in U$, and, in addition, the map $x \rarrow Df(x)$ is continuous from $U$ into
$L(E,F)$, then we say that $f$ is {\it continuously differentiable} or
$C^1$ in $U$.  

If the $C^1$ map $f:U \rarrow F$ is such that $y \rarrow Df(y)$ is
$\alpha-$H{\"{o}}lder at $x$, (resp., in $U$), then we say that $f$ is
$C^{1,\alpha}$ at $x$ (resp.,  in $U$).  Unless otherwise stated, we
will typically use this for $0 < \alpha <  1$.  When $\alpha = 1$, we will
say that $Df$ is {\it Lipschitz} at $x$ (resp.,  in $U$). 

For a $C^{1,\alpha}$ map $f:U \rarrow E$, it is convenient to define its
{\it $D-$H{\"{o}}lder constant} $Hol(Df,U)$ by

\[  Hol(Df,U) = \sup_{x \neq y \in U}\frac{\n{Df(x) -  Df(y)}}{\n{x-y}^{\alpha}}. \]

Recall that two norms $\n{\cdot}$ and $\dnorm{\cdot}$ on a linear space
$E$ are {\it equivalent} if there is a positive constant $C>0$ such
that 

\[ C^{-1} \leq \frac{\n{x}}{\dnorm{x}} \leq C \]

for every $x \neq 0$  in $E$.

Note that differentiability of functions $f:E \rarrow F$ is
independent of the choice of equivalent norms on $E$ or $F$.

Throughout this paper, we denote the ball of radius $\delta$ about $0$
in $E$ by $B_{\delta} = B_{\delta}(0)$. 

Let $0 < \alpha < 1$.  We will call a Banach space $(E, \n{\cdot})$ a 
$C^{1,\alpha}$ {\it Banach space} if there is a $C^{1,\alpha}$ bump function on $E$. This
is a $C^{1,\alpha}$ function $\grf:E \rarrow \mathbb{R}$ such that 
$\grf(E) = [0,1]$, and there are positive numbers $c_1 <
c_2$ such that $\grf(x) = 1$ for $x \in B_{c_1}$ and
$\grf(x) = 0$ for $x \notin B_{c_2}$. Replacing $\grf(x)$ by $x
\rarrow \grf(c_2x)$, we may assume that $c_2=1$. 
This is a somewhat restrictive condition on Banach spaces. 
For instance, it holds for the $L_p$ spaces with $p > 1$ but fails for $L_1$. If
there is an equivalent norm $\dnorm{ \cdot }$ on $E$ which is $C^{1,
\alpha}$ on open sets which do not contain $0$, then it is easy to
see that such bump functions exist. In particular, this is true in
Hilbert Spaces since their norms are $C^{\infty}$ away from $0$.  The study of the existence or non-existence of various
smooth norms or bump functions in Banach spaces is a central part of the geometric
theory of Banach spaces. We refer the reader to \cite{Fry-McManus} 
and \cite{Johnson-Lindenstrauss-Handbook} 
for more information.  

Given a linear automorphism  $L:E \rightarrow E$, the 
{\it spectral radius } of $L$, denoted $\rho(L)$, is defined to be 

\beq \label{spectral-radius-formula}
\rho(L) = \lim_{n \rarrow \infty} \n{L^n}^{\frac{1}{n}} 
\eeq

To see that this limit exists, observe that the sequence $a_n =
\log(\n{L^n})$ is subadditive (i.e. $a_{n+m} \leq a_n + a_m$) and the
numbers $\frac{a_n}{n}$ are bounded below (in fact, by the number 
$ -log(\n{L^{-1}})$),
 so Lemma 1.18 in \cite{Bowen-75-1} implies that  the quantity

\[ h_0 \eqdef \lim_{n \rarrow \infty} \frac{a_n}{n} = \inf_{n \geq 1} \frac{a_n}{n} \]

exists.  Then, $\rho(L) = exp(h_0)$. 

If $E$ is $n-$dimensional Euclidean space and $L$ is a linear
automorphism, then $\rho(L)$ is, of course, the maximum of the absolute values
of the eigenvalues of $L$. 

A linear automorphism $L$ is called {\it contracting} if $\rho(L) < 1$.
 It is {\it expanding} if $\rho(L^{-1}) < 1$. 

A linear automorphism $L$ is called {\it hyperbolic} if it satisfies
exactly one of the three conditions
\be
 \item $L$  is contracting,
 \item $L$  is expanding, or
 \item $L$  can be written as a direct sum $L = A \oplus B$ in which $A$ is
expanding and $B$ is contracting.
\ee

    In the last case, we say that $L$ is {\it hyperbolic of saddle
type} or, simply {\it of  saddle type}, and we call $A$ the
expanding part of $L$ and $B$ the contracting part of $L$.  They are
unique. 

The {\it spectral condition number} of the linear automorphism $L$ is
the quantity $c(L) = \rho(L)\rho(L^{-1})$. It is an easy consequence of the spectral radius formula
(\ref{spectral-radius-formula}) that $c(L) \geq 1$. Indeed, for each
positive integer $m$, we have $I = L^{-m}L^m$, so $1 \leq
\n{L^{-m}}\n{L^m}$. Now, take $m-$th roots and the limit as $m$
approaches infinity. 

A linear automorphism $L$ is called {\it circular} if
$c(L) = 1$. In the finite dimensional case, this means
that all of its eigenvalues lie in a single circle in the complex
plane. If $L$ is circular and contracting, then we call it
$c-$contracting. Analogously, we say that $L$ is $c-$expanding if $L^{-1}$ is
$c-$contracting.  We say that $L$ is {\it bi-circular} if it is
hyperbolic of saddle type,  and its expanding and contracting parts are both circular. 

Trivially, hyperbolic linear automorphisms of saddle type on the plane
are bi-circular as are hyperbolic linear automorphisms of saddle type
in $\mathbb{R}^3$ which have a pair of non-real complex conjugate
eigenvalues. 

An obvious, but important remark, is that {\it every
hyperbolic linear automorphism of saddle type on a finite dimensional
space is bi-circular on the direct sum of its leading eigenspaces}. These are the subspaces on
which the eigenvalues are closest to 1 in absolute value (see Section
\ref{Motivation} for details on the analogous conditions for vector
fields and applications). 

Given the real number $\alpha$ with $0 < \alpha <  1$, we say that $L$ is
{\it $\alpha-$contracting} if $c(L)\rho(L)^{\alpha} < 1$. 
Since $c(L)$ is always at least 1, this implies that
$\rho(L) < 1$.  It is well-known (see Lemma \ref{re_norm} below) that this
implies the existence of an equivalent norm on $E$ whose induced norm
on $L$ is less than 1 (whence the name contracting). 
In the finite dimensional case, this
implies that the absolute values of the eigenvalues lie in a complex
annulus whose bounding circles have radii in the open-closed interval
$(\rho(L)^{1+\alpha}, \rho(L)]$. 

We say that $L$ is $\alpha$-expanding if $L^{-1}$ is
$\alpha$-contracting. 

Let $L$ be a hyperbolic linear automorphism of saddle type, and let $L = A \oplus B$ with $A$ expanding
and $B$ contracting.

Let 

\beq \label{c_h-def}
c_h = c_h(L) = c_h(A,B) = \max(c(A),c(B)).
\eeq

We call this the {\it hyperbolic condition number of $L$}.

Similarly, we define the {\it hyperbolic spectral radius of $L$} to be

\beq \label{rho_h-def}
\rho_h = \rho_h(L) = \rho_h(A,B) = \max(\rho(A^{-1}), \rho(B)).
\eeq

We say that $L$ is {\it $\alpha-hyperbolic$} if it can be written as a
direct sum $L = A \oplus B$ so that 

\beq \label{gra-hyperbolic-def}
c_h \rho_h^{\alpha} = c_h(L) \rho_h(L)^{\alpha} < 1.
\eeq

\begin{Rem}
{\rm
Since definition of $\rho_h$ in (\ref{rho_h-def}) requires that the inverse
operator appearing is the {\it expanding part} in a hyperbolic
decomposition of $L$, it follows that $\rho_h(L^{-1}) = \rho_h(L)$. 
}
\end{Rem}

\begin{Rem}
{\rm Since the condition (\ref{gra-hyperbolic-def}) implies that
 \[ c(A)\rho(A^{-1})^{\alpha} < 1 \mbox{ and } c(B)\rho(B)^{\alpha} < 1,
\]
 we see   that the contracting and expanding parts of an $\alpha-$hyperbolic
  automorphism must, in fact, be $\alpha-$contracting and
  $\alpha-$expanding, respectively. 
}
\end{Rem}

Note that if $L$ is bi-circular, then it is $\alpha-$hyperbolic for
every $\alpha \in (0,1)$. 

Let $Hyp_{\alpha}(E)$ be the set of $\alpha-$hyperbolic linear
automorphisms of $E$.  Using the formula
(\ref{spectral-radius-formula}) on the contracting and expanding parts
of an element of $Hyp_{\alpha}(E)$ and the techniques in Section 4 of
Hirsch and Pugh \cite{Hirsch-Pugh}, it can be shown that, for each $\alpha \in
(0,1)$, $Hyp_{\alpha}(E)$ is an open subset of $L(E,E)$. 

Let $r \geq 1, k \geq 1$ be real numbers (not-necessarily integers). 
Let $[r]$ denote the greatest integer less than or equal to $r$.  When
$r$ is not an integer, we say that a map $T$ is $C^r$ if it is
$C^{[r]}$ and its $[r]-$th derivative is $r-[r]$-H{\"{o}}lder.
Setting $r-[r] = \beta$, we
will also write this as $C^{[r],\beta}$.

Given the Banach space $E$ and a point $p \in E$, let $T$ be a
$C^r$ diffeomorphism from an open neighborhood $U$ of $p$ onto
its image such that $T(p) = p$.  We say that $p$ is an {\it
$\alpha-$contracting} fixed point of (or for) $T$ if the derivative $DT(p)$ is
$\alpha-$contracting. In a similar way we define {\it $\alpha-$expanding},
{\it $\alpha-$hyperbolic}, and {\it bi-circular} fixed points. 

A {\it $C^k$ linearization} of $T$ at $p$ (or near $p$) is a
$C^k$ diffeomorphism $R$ from a neighborhood of $p$ onto a
neighborhood of $0$ such that $R T R^{-1} = L$ on some neighborhood of
$0$.  An equivalent condition is that $LR = RT$ on some neighborhood
of $p$.  When such a diffeomorphism $R$ exists, we say that $T$ is
$C^k$ {\it linearizable} at (or near) $p$.  Sometimes the term
{\it locally $C^k$ linearizable} at (or near) $p$ is used and the map $R$ is
called a {\it local linearization} at (or near) $p$.  At times it will
be convenient to use the statement {\it $p$ is $C^k$ linearizable} to
mean that $p$ is a fixed point of a $C^r$ diffeomorphism $T$ (for some
$r \geq k$) which has a $C^k$ linearization at $p$.  Since the
neighborhoods in the domains of definition of the various maps
considered often change, the concept of {\it germ} of a map at a point
is often used (as in \cite{Chaperon-Steklov}, \cite{Ilyashenko-Li}). 

As far as we know, S. Sternberg was the first to obtain linearization
results for finitely smooth systems without explicit diophantine
inequality assumptions.  In the paper
\cite{Sternberg-Loc-Con-Real-Line}, published in 1957,  he proved
that, for $k \geq 2$,  $C^k$ maps of an interval with
contracting or expanding fixed points are locally $C^{k-1}$ linearizable
near the fixed points.  

In \cite{Hartman-saddle}, Philip Hartman extended the $C^2$ case of Sternberg's one dimensional
results in the following significant ways. 

\begin{Theorem}(Hartman).\label{Hartman-Theorem} Let $E$ denote the Euclidean space $\R^n$,
and let $U$ be an open neighborhood of $0$ in $E$.  Let $T$ be a
$C^{1,1}$ diffeomorphism from $U$ to its image such that $T(0) = 0$ and $DT(0) =
L$ where either 

\be
 \item[(a)] $L$ is contracting\footnote{Since $(RTR^{-1})^{-1} =
RT^{-1}R^{-1}$, the theorem also applies when $L$ is expanding.}, or 
 \item[(b)] $L$ is $\alpha-$hyperbolic for some $\alpha \in (0,1)$. 
\ee

Then, $T$ is $C^1$ linearizable at $0$.
\end{Theorem}

As is well-known, the more general case in which $T(p) = p$ and $p$ is not necessarily $0$ can
be reduced to the case of the theorem by replacing $T$  by $S^{-1}T S$ where
$S(x) = x+p$ is the translation by $p$.  

\begin{Rem} 
{\rm In the contracting case, Hartman proved that the
linearization $R$ was $C^{1,\beta}$ for some $\beta > 0$. He stated
that this was true in the $\alpha-$hyperbolic case, but did not present
the proof. 
}
\end{Rem}

It is natural to ask if one can reduce the regularity 
assumption on $T$ and still  get a $C^1$ linearization. 

On page 101 in \cite{Sternberg-Loc-Con-Real-Line} Sternberg shows
that, for any $a \in (0,1)$,  the $C^1$ map $T$ on the real line 
defined by

\[ T(x) = 
   \left\{ \barr{ll} 
          ax(1 - \frac{1}{log (\mid x \mid)}) & \mbox{ if } x \neq 0 \\
          0                             & \mbox{ if } x = 0 
           \earr
   \right.
\]

has no $C^1$ linearization at $0$. 

Hence, the most natural
assumption is that we consider the case the $DT$ is uniformly H{\"{o}}lder
instead of uniformly Lipschitz.  That is, we assume that $T$ is
$C^{1,\alpha}$ (on some neighborhood of $0$) for some positive $\alpha$ with $0 < \alpha < 1$. 

It is known (even in arbitrary Banach spaces) that a $C^{1,\alpha}$
diffeomorphism $T$ with an $\alpha-$contracting fixed point at $0$ has a $C^{1,\alpha}$ linearization
near $0$.  This is a consequence of Corollary
1.3.3 in Chaperon \cite{Chaperon-Steklov}.  Since this
result is fundamental for our work here, we will include an elementary
proof (see Theorem \ref{Lin_Contraction_2}) below.
In the finite dimensional case, this result was stated (under the
stronger assumption that $T$ was $C^{1,1}$) with different
notation in the last sentence of the first paragraph in section 8 on page
235 in \cite{Hartman-saddle}. There is also a statement on page 222 of
\cite{Hartman-saddle} that the regularity assumption on $T$ could be
weakened to some $C^{1,\grg}$ with $\grg > \alpha_0$ where $\alpha_0$
depends on the eigenvalues of $L$.  We suspect that this $\alpha_0$
equals our $\alpha$, but, since it is not given explicitly, we cannot be
sure of this.  In the infinite dimensional case, the weaker result
that there is a $C^1$ linearization whose derivative  is 
H{\"{o}}lder at $0$ was proved by Mora and Sol{\`{a}}-Morales in Theorem 3.1 in \cite{Mora-Sola-Morales}. 

Since any circular linear contraction is $\alpha-$contracting for any
$\alpha > 0$, it follows that {\it any $C^{1,\alpha}$ map $T$ with a c-contracting fixed
point $p$ has a $C^{1,\alpha}$ linearization at $p$}. 

The next case to consider is contractions which are not circular.
Here, in  \cite{Zhang-Zhang-Lin}, Zhang and Zhang study such $C^{1,\alpha}$
contractions in the plane. They obtain $C^{1,\grb}$ linearizations for various $\grb$
depending on $\alpha$.  They also show that any linear non-circular contraction can be perturbed to
give a $C^{1,\alpha}$ diffeomorphism fixing the origin which has no 
$C^{1,\grb}$ linearizations for any $\grb > 0$.   It is not known in these examples
if there is a $C^1$ linearization. 

In \cite{Rodrigues-Sola-Morales-Survey}, Rodrigues and
Sol{\`{a}}-Morales give examples of $C^{1,1}$ diffeomorphisms in
infinite dimensional Banach spaces with contracting fixed points which
are not $C^1$ linearizable.

An alternative approach to linearizations of $C^2$ flows via the Lie
Derivative is in Chicone-Swanson \cite{Chicone-Swanson-Lie}. 

The main result in the present paper can be paraphrased as follows. 

{\it  In a $C^{1,\alpha}$ Banach
space with $0 < \alpha < 1$, every $\alpha-$hyperbolic fixed point of a
$C^{1,\alpha}$ diffeomorphism is $C^{1,\beta}$ linearizable for some $\beta > 0$. 
}

After translating the fixed point to the origin as above, we have the
more precise statement. 

\begin{Theorem} \label{Main_Theorem_1}
Let $0 < \alpha < 1$ and let $E$ be a $C^{1,\alpha}$ real Banach space. 
Let $U$ be a neighborhood of $0$ in $E$
and let $T:U \rarrow E$ be a $C^{1,\alpha}$ map such that $T(0) = 0$
and $L=DT(0)$ is an $\alpha-$hyperbolic automorphism of $E$. 

Then, there are a subneighborhood $V \subset U$ of $0$, a real number
$\beta \in (0, \alpha)$, 
and a $C^{1,\beta}$ diffeomorphism $R$ from $V$ onto its image such
that $R(0) = 0$ and  

\beq \label{Lin_L_f}
 L(R(x))= R(T(x))  \mbox{ for all } x \in V \bigcap T^{-1} V 
 \eeq

\end{Theorem}

{\bf Remarks.} 
\be
\item 
In the contracting $C^{1,1}$ case (i.e., part (a) of Theorem \ref{Hartman-Theorem}), Hartman proved that the linearization
could be made $C^{1,\grb}$ for some $\grb > 0$.  In the
$\alpha-$hyperbolic case, he stated that his proof could be modified to
give a linearization for some $\grb > 0$, but he did not
present the proof.  Thus, our proof, even in the finite dimensional
case, may  be the first available proof that such a  $C^{1,\grb}$
linearization exists. 
\item 
Given real numbers $a > b > 1 > c > 0$, and $\gre \neq 0$,
Hartman shows in \cite{Hartman-saddle} that, if $b = a c$, then
the quadratic polynomial map  
\beq \label{counter-example}
T_{a,b,c}(x,y,z) = (a x,\  b(y + \gre x z),\  c z)
\eeq 

has no $C^1$ linearization at $(0,0,0)$.  

It is interesting to note that Hartman's example is sharp in the
following sense.  If $b \neq ac$, and $R(x,y,z) = (x, y +
\frac{b\epsilon}{b-ac}xz, z)$, then $LR = RT$ so the map $R$ gives even a
quadratic polynomial linearization of $T$ at $(0,0,0)$.

\item In  \cite{Rodrigues-Sola-Morales-Linearization-Saddles-Banach-Spaces}, Rodrigues and
Sol{\`{a}}-Morales consider $C^{1,1}$ diffeomorphisms having a
hyperbolic saddle fixed point at $0$ in a $C^{1,1}$ Banach space. They
prove that, under a spectral condition which is equivalent to our
notion of $\alpha$-hyperbolicity, the diffeomorphisms are $C^1$ linearizable {\it provided $\alpha$
is sufficiently close to 1}.  See  Theorem 2 in
\cite{Rodrigues-Sola-Morales-Linearization-Saddles-Banach-Spaces}. 
\item Given real numbers $r \geq 1$ and $r \geq k \geq 0$, let us say that a fixed
point $p$ of a $C^r$ diffeomorphism $T$ is {\it robustly $(r,k)$-
linearizable} if, for any $S$ sufficiently $C^r$ close to $T$,
there is a pair $(p_S,R_S)$ with $p=p_T$, depending continuously on $S$, such that
$S(p_S) = p_S$ and $R_S$ is a $C^k$ linearization of $S$ at $p_S$.

More precisely, let $U$ be an open set in
the Banach space $E$.  Consider the space $D^r(U, E)$ of $C^r$
diffeomorphisms $T$ from $U$ to $T(U)$ with the uniform $C^r$
topology.  For $k=0$, let $D^k(U,E)$ denote the space of
homeomorphisms $T$ from $U$ onto $T(U)$ with the uniform $C^0$
topology. The fixed point $p$ of the
diffeomorphism $T \in D^r(U,E)$  is said to be 
{\it robustly $(r,k)$-linearizable} if there are neighborhoods $W$ of $T$ in $D^r(U,E)$
and $V$ of $p$ in $E$ such that, for any $S \in W$, there is a pair
$(p_S, R_S)$ such that

\be
 \item $p = p_T$,
 \item $S(p_S) = p_S$, 
 \item $p_S \in V$, 
 \item $R_S \in D^k(V,E)$,
 \item $R_S(p_S) = 0$, 
 \item $DS_{p_S}(x) = (R_S \circ S \circ R_S^{-1})(x)$ for $x \in
R_S(U \bigcap S^{-1}V)$, and 
 \item the map $S \rarrow (p_S, R_S)$ from $W$ into $E \times
D^k(V,E)$ is continuous. 
\ee
 
The typical linearization results that we are aware of 
(e.g. those given by the Grobman-Hartman Theorem, Sternberg's
Theorems, Hartman's Theorem \ref{Hartman-Theorem} above, etc) yield
fixed points which are robustly $(r,k)-$linearizable for some
$(r,k)$. 

    As a consequence of  Theorem \ref{par-Main-Theorem} below we
also obtain that, for $r = 1+\alpha$, the $\alpha-$hyperbolic and
bi-circular fixed points considered in this paper, are, in fact,
robustly $(r,1+\beta)$ linearizable for some $\beta > 0$. 
 \item 
This paper owes much to Hartman's paper \cite{Hartman-saddle}.  The
scheme of the proof of Theorem \ref{Main_Theorem_1} is similar to that indicated in
pages 235-238 in \cite{Hartman-saddle} for $C^{1,1}$ maps $T$.  In actuality, Hartman
only gave the proof that the linearization was $C^0$.  He stated that
similar methods could be used to prove that it was $C^1$, and, as we
have already mentioned,  he also stated that the
proof could be modified to give a linearization whose derivative was
uniformly H{\"{o}}lder continuous. 
\item 
 In the fifty-five  years since the paper \cite{Hartman-saddle} appeared,  
two significant developments have occurred that make our arguments
possible. 

\be
 \item We now know that no smoothness is lost for stable and unstable
manifolds at a hyperbolic fixed point.  In particular, if $T$ is
$C^{1,\alpha}$ with a hyperbolic fixed point at $p$, then the stable and
unstable manifolds at $p$ are locally the graphs of $C^{1,\alpha}$
maps. 
 \item A $C^{1,\alpha}$ map $T$ with fixed point $p$ such
that $DT(p)$ is $\alpha-$contracting has a $C^{1,\alpha}$ linearization. See Theorem
\ref{Lin_Contraction_2} below. 
\ee
\item There are several additional papers in the literature which are
relevant to the work presented here. In particular, the papers by Samovol \cite{Samovol-1990}, \cite{Samovol-1991},
\cite{Samovol-1995},  Belitskii \cite{Belitskii-1973}, \cite{Belitskii-Counter-Example-1976}, 
\cite{Belitskii-1978}, and Stowe \cite{Stowe-1986} all contain
interesting results on smooth linearizations with finite  class of differentiability.  Many of these
are discussed and generalized in the book of Bronstein and Kopanskii 
\cite{Bronstein-Kopanskii}. 
\item The contents of this paper are as follows.  The proof of Theorem \ref{Main_Theorem_1} will be carried out in
Sections (\ref{Associated-Ops})-(\ref{Main_Proof}).  Section
\ref{par-dep} considers the dependence of the linearization $R$ in
Theorem \ref{Main_Theorem_1} on external parameters, and Section
\ref{Motivation} presents some motivation and applications of the results. 
\ee

{\bf Acknowledgements.} 
\be
 \item Being relatively unaware of the recent
literature on the subject, we gave a lecture at Penn State in October,
2014, in which we discussed the generalization to the $C^{1, \alpha}$
case of part (b) in  Hartman's Theorem \ref{Hartman-Theorem} above for
bi-circular fixed points.The result was still for finite dimensional systems, including the finite 
dimensional version of Theorem \ref{Lin_Contraction_2} below, and only
considered the existence of $C^1$ linearizations.  We were not aware that the general result in Theorem \ref{Lin_Contraction_2} was known.
We wish to thank Misha Guysinsky for subsequently
providing many references to the recent literature, including, in
particular, \cite{Guysinsky-Hasselblatt-Rayskin} and
\cite{Zhang-Zhang-Lin}. His remarks and references provided substantial impetus for us to
study the recent literature, and eventually extend our results to the infinite dimensional case as presented here. 
 \item We have already indicated the importance of  Hartman's paper
\cite{Hartman-saddle} in connection with this work. This was a special
time for him in that, on May 16, 2015, he reached his 100th
birthday. We recently learned, sadly, that he passed away on August
28, 2015. 
 We dedicate this paper to both Anosov and Hartman--two icons in
Dynamical Systems whose mathematical achievements will have a long
lasting legacy. 
\ee

\section{The Associated Algebraic Problem} \label{Associated-Ops}

The standard way to approach the functional equation (\ref{Lin_L_f}) is to
convert it into a related algebraic equation on a suitable Banach space $\sE$ as follows.

Let $I$ denote the identity map on $\sE$. 

Writing $T = L+f$, we try to find $R$ of
the form $R = I + \grf$ with $\grf \in \sE$ so that 
the equation 

\[ LR = RT \]

becomes

\[ L(I + \grf) = (I+ \grf)T = (I + \grf) (L + f) \]

which leads to 

\bqns
 L + L \grf & = & L+f + \grf  (L+f) \\
 L \grf & = & f + \grf (L+f) \\
   \grf & = & L^{-1} f + L^{-1} \grf (L+f),
\eqns
 
or 

\beq \label{Alg-eq-1}
\grf - L^{-1} \grf (L+f)  =  L^{-1} f.
\eeq

Defining the operator $H$ on functions $\grf$ by 

\beq \label{H-def}
H(\grf) = L^{-1} \grf (L+f) = L^{-1} \grf \circ T,
\eeq

we see that $H$ is linear, and equation (\ref{Alg-eq-1}) 
becomes 

\beq \label{alg-prob}
(I - H)(\grf) = L^{-1}f. 
\eeq

The problem, then,  is reduced to finding a Banach space $\sE$ so that the
operator $H$ is a  well-defined bounded linear map such that $I - H$
has a bounded linear inverse, in which case we obtain 

\beq \label{alg-soln} 
\grf =  (I - H)^{-1}L^{-1}f.
\eeq

Theorem \ref{Main_Theorem_1} will be proved by using the map $H$
defined  by  (\ref{H-def}) and the ensuing equations 
 (\ref{alg-prob}) and (\ref{alg-soln}) in several different function spaces and then combining the results. 

A rough outline (assuming the hypotheses of Theorem \ref{Main_Theorem_1}) is as follows.

{\bf Step 1:} We show that if $L$ is contracting; i.e., $\rho(L) < 1$, then
$T = L + f$ can be $C^{1,\alpha}$ linearized near $0$.   That is,
there is a local $C^{1,\alpha}$ diffeomorphism $R$ near $0$ satisfying $R T R^{-1} = L$
near $0$.  It follows, of course, that $R T^{-1} R^{-1} = L^{-1}$ near
$0$. 

{\bf Step 2:} Assuming $L$ is hyperbolic of saddle type, let $E = E^u
\oplus E^s$ be the associated unstable and stable subspaces of $E$.
We choose $C^{1,\alpha}$ coordinates near $0$ so that the local stable
and unstable manifolds are flattened; i.e., are contained in 
$E^s$ and $E^u$, respectively. Applying the result in Step 1, 
we can locally $C^{1,\alpha}$ conjugate $T$ to a map $S_1$ which 
locally preserves $E^s$ and $E^u$ and is linear when restricted to
those subspaces. Then, using 
a bump function, we extend the map $S_1$ (actually
its restriction to a small neighborhood of $0$) to a
$C^{1,\alpha}$ diffeomorphism $S$ from the whole Banach space $E$ onto
 itself such that $S(0) = 0$, $DS(0) = L$, $Lip(S-L)$ is small, and 
 $S = L$ on the union of $E^u$, $E^s$, and the complement of a small
 neighborhood $V$ of $0$. 

{\bf Step 3:} For $V$ and $Lip(S-L)$ small enough, we define an
  appropriate Banach space $\sE$ of $C^1$ functions $\phi$ defined on all of $E$ and solve the
  equation (\ref{alg-soln}) to obtain a global $C^1$ conjugacy $I +
  \phi$ with $\phi \in \sE$ from $S$   to $L$ defined on all of
  $E$.

{\bf Step 4:} Choosing $V$ small enough, and making use of
(\ref{alg-soln}) in two separate Banach spaces of functions, we show
that, for appropriately chosen small $\beta > 0$, the map  $\phi$
is $C^{1,\beta}$ on bounded subsets of the orbit of
$V$. 

This gives that $S$ is locally $C^{1,\beta}$ linearizable at $0$, and, hence, so is $T$.

The operator $H$ occurs often enough in this kind of problem that it
deserves a name. We call it the {\it associated linearization operator for $T$} or the
$T-$linearization operator.

\section{The $\alpha-$Contracting Case} \label{Contracting_Case}

Let $E$ be a Banach space, and let $U$ be an open neighborhood of $0$
in $E$.  Let
$0 < \alpha < 1$ and consider the set $C^{1,\alpha}_U = C^{1,\alpha}_{0,U}$ of $C^1$
functions $g: U  \rarrow E$ such that 

\beq \label{g_Dg_0}
g(0) = 0, \ Dg(0) = 0, 
\eeq 

and

\beq \label{unif_gra_holder}
Hol(Dg,U) = \sup_{x, y \in U, x \neq y}
\frac{\n{Dg(x) - Dg(y)}}{\n{x-y}^{\alpha}} < \infty. 
\eeq

Using the techniques for statement (8.6.3) in \cite{Dieudonne-Analysis} 
one can verify that the D-H{\"{o}lder constant of $g$, 

\[ \n{g}_U = \n{g}_{U,\alpha} = Hol(Dg,U), \] 

defines a norm on $C^{1,\alpha}_U$ making it into a Banach
space. 

Observe that the usual way of defining a norm on  a space of $C^r$
functions into a Banach space involves taking the supremum of the
$C^0$ size of the function as well as that of its
derivatives. Otherwise one only gets a semi-norm in that the vanishing
of this semi-norm does not imply the vanishing of the functions.
However, our functions already vanish at $0$ (and $U$ is connected), so this issue does not
occur.

Note that if $V \subset U$ and $\grf \in C^{1,\alpha}_U$, then the restriction map from
$\grf \rarrow \grf | V$ 
induces a continuous  map  from $(C^{1,\alpha}_U, \n{\cdot}_{U,\alpha})$ into $(C^{1,\alpha}_V,
\n{\cdot}_{V,\alpha})$. We will use the same notation for the maps $\grf$ and
$\grf \mid V$ letting the context determine which space we are using
at a given moment. 

When $U$ is the open ball $B_{\delta}(0)$ we will denote $C^{1,\alpha}_U$
by $C^{1,\alpha}_{\delta}$. 

Let $I$ denote the identity map on $E$. 

The goal of this section is to prove the  following theorem

\begin{Theorem} \label{Lin_Contraction_2} Let $U$ be an open
neighborhood of $0$, let $\alpha$ be a positive
real number with $0 < \alpha < 1$ and let $T$ be a map of the form
$T(x) = Lx + f(x)$ where $L:E \rarrow E$ is a linear
$\alpha-$contracting automorphism and $f \in C^{1,\alpha}_U$. 

Then, for sufficiently small $\delta$, there is a unique map $\grf \in C^{1,\alpha}_{0,\delta}$ such that, setting $R = I + \grf$, we have 

\beq \label{R-eq-1}
 L(R(x)) = R(T(x)) \mbox{ for all } x \in
B_{\delta}(0) \bigcap T^{-1} B_{\delta}(0).
\eeq

\end{Theorem} 

\begin{Rem} \label{Rem-c-expansion}
{\rm For $\delta$ small, the map $R$ will be Lipschitz close to the
identity. Therefore, the Inverse Function Theorem implies that $R$ is a $C^{1,\alpha}$ diffeomorphism onto its
image with inverse which is also  $C^{1,\alpha}$.  Also, $L = R T
R^{-1}$ implies that  $L^{-1} = R T^{-1} R^{-1}$, so  Theorem 
\ref{Lin_Contraction_2} also holds  when $L$ is a linear
$\alpha-$expansion. 
}
\end{Rem}

\begin{Rem}
{\rm Consider a $C^{1,\alpha}$ diffeomorphism $T:\mathbb{R}^N \rarrow \mathbb{R}^N$ on the
Euclidean space $\mathbb{R}^N$ with a hyperbolic fixed point at $0$ and
derivative $L = DT(0)$. Then, we have the direct sum 
decomposition 

\[ \mathbb{R}^N = E_1 \oplus E_2 \oplus \ldots \oplus E_s \]

into $L-$invariant subspaces such that $L \mid E_i$ is circular for
each $1 \leq i \leq s$. 

From invariant manifold theory,  \cite{Hirsch-Pugh-Shub}, \cite{Llave-Wayne}, for some $0 < \alpha < 1$, there
are $T$ locally invariant $C^{1,\alpha}$ submanifolds $W_i$ tangent at
$0$ to $E_i$ for each $i$.  These submanifolds are generally not 
unique, of course, and, even if $T$ is linear, they are not generally 
$C^{1,1}$. Nevertheless, using Theorem
\ref{Lin_Contraction_2}, Remark \ref{Rem-c-expansion}, and the techniques in Section \ref{Flattening-Subsection} below,
one can find an $\alpha \in (0,1)$ and a  local $C^{1,\alpha}$ coordinate
chart $(U, \zeta)$ near $0$ such that $\zeta(0)=0$, $\zeta(W_i)
\subset E_i$ for each $i$, and $\zeta T \zeta^{-1}$ restricted  to
each $E_i$ is linear near $0$. 
}
\end{Rem}

As we mentioned in the introduction, Theorem \ref{Lin_Contraction_2}
is not new. It is a consequence of a more general result proved by
Chaperon (see Corollary 1.3.3 in \cite{Chaperon-Steklov}). 
Since the result is needed in an essential way in the present paper and
it is relatively simple to prove, we give an alternate proof here for completeness. 

We begin with a standard result which shows that, given a Banach space
$(E, \n{\cdot})$, and a bounded linear operator $L:E \rarrow E$, we can approximate the spectral
radius $\rho(L)$ by the supremum norm $\n{L}_1$ induced by a
norm $\n{\cdot}_1$ on $E$ which is equivalent to $\n{\cdot}$. 

\begin{Lem} \label{re_norm}
 Let $(E, \n{\cdot})$ be a Banach space and let $L: E \rarrow E$ be a bounded linear map with spectral radius
$\rho(L)$.   Then, given $\gre > 0$ there is a norm $\n{\cdot}_1$ on $E$
which is equivalent to $\n{\cdot}$ such that  the induced norm
$\n{L}_1$ defined by

\[ \n{L}_1 = \sup_{\n{v}_1 = 1} \n{Lv}_1 \]

satisfies

\[ \n{L}_1 \leq \rho(L) + \gre \]
\end{Lem}

\begin{proof}

Setting $\rho_1 = \rho(L) + \gre$ and using 
formula (\ref{spectral-radius-formula}), we can choose a positive
integer $q$ such that, for $n \geq q$, we have 

\[ \n{L^n} < \rho_1^n. \] 

This gives, for any $v$ and $n \geq q$, 

\[ \n{L^n v} \leq \n{L^n} \n{v} \leq \rho_1^n \n{v}, \]

which implies that 

\[ \rho_1^{-n} \n{L^nv} \leq \n{v}. \]

As usual, we define $L^0 = I$. 

Then, setting 

\[ K = \max_{0 \leq n < q}  \rho_1^{-n}\n{L^n}, \]

we get, for $n < q$, 

\[ \rho_1^{-n}\n{L^nv} \leq \rho_1^{-n}\n{L^n} \n{v} \leq K \n{v}. \]

For $n = 0$, and $\n{v} \neq 0$, this gives $K \geq 1$.

Putting these inequalities together, gives 

\[ \sup_{n \geq 0} \rho_1^{-n}\n{L^nv} \leq K \n{v} \]

for any vector $v \in E$, 

Defining a new norm by 

\[ \n{v}_1 = \sup_{n \geq 0} \rho_1^{-n}\n{L^nv}, \] 

we then have

\[ \n{v} \leq \n{v}_1 \leq K \n{v}\]

so $\n{\cdot}_1$ is equivalent to $\n{\cdot}$. 

Moreover, for every $v$, we have 

\bqns
 \n{Lv}_1 & = & \sup_{n \geq 0} \rho_1^{-n}\n{L^nLv} \\
    & = &  \sup_{n \geq 0} \rho_1^{-n}\n{L^{n+1}v} \\
    & = &  \rho_1 \sup_{n \geq 0} \rho_1^{-n-1}\n{L^{n+1}v} \\
    & = &  \rho_1 \sup_{n \geq 1} \rho_1^{-n}\n{L^{n}v} \\
    & \leq & \rho_1 \n{v}_1.
\eqns

Taking any $v$ with $\n{v}_1 = 1$, this gives $\n{Lv}_1 \leq \rho_1$,

whence 

\[ \n{L}_1 = \sup_{\n{v}_1 = 1} \n{Lv}_1 \leq \rho_1 \]

as required to  complete the proof of Lemma \ref{re_norm}. 
\end{proof}

Let us say that a bounded linear map $H$ of a Banach
space $E$ is {\it contracting} if its spectral radius 
is less than 1.  

Theorem \ref{Lin_Contraction_2} will be proved using the method
described in section \ref{Associated-Ops}. In this case, the
corresponding operator $H$ will, in fact, be contracting. 
 Using Lemma \ref{re_norm}, we find an equivalent norm in
which the induced norm on $H$ satisfies $\n{H} < 1$. Then, as is well known, 
$I-H$ is invertible with inverse given by 

\[ (I-H)^{-1} = \sum_{n=0}^{\infty} H^n. \]

For $\delta > 0$ sufficiently small, we take the space $\sE$ alluded to
in section \ref{Associated-Ops} to be $C^{1,\alpha}_{\delta}$, and we use
the norm $\n{\grf}_{\delta} = \n{\grf}_{\delta,\alpha}$ for elements $\grf \in \sE$. 

Let $U, L, f, T$ be as in the hypotheses of Theorem
\ref{Lin_Contraction_2}.

All numbers $\delta$ below will be chosen so that $B_{\delta}(0) \subset U$, and, hence, $f$ is defined and $C^{1,\alpha}$ in $B_{\delta}(0)$.

Recalling the procedure in section \ref{Associated-Ops}, we consider the functional equation 

\[ L(I + \grf) = (I + \grf)(L+f) \]

and the resulting algebraic equation

\[ (I - H)(\grf) = L^{-1}f. \]

where $H(\grf) = L^{-1} \grf (L+f)$.

We will show that, for $\delta$ small enough,  the operator $H$ has the
properties that 

\beq \label{H-def-C}
H \mbox{ {\it is defined as a function on} }
C^{1,\alpha}_{\delta}, 
\eeq

\beq \label{H-into-C}
H \mbox{  {\it maps} } C^{1,\alpha}_{\delta} \mbox{
{\it into itself} }
\eeq

and, for some positive integer $m$, 

\beq \label{norm-Hm}
\n{H^m} < 1. 
\eeq

In then follows that $\rho(H) < 1$, and, as we have already noted,
the required solution is obtained as $\grf = (I-H)^{-1} L^{-1}f$. 

Before going to the proof of (\ref{H-def-C})-(\ref{norm-Hm}), we need some easy
estimates obtained from the Mean Value Theorem and our other
assumptions. 

As we have already noted, the $\alpha-$contracting condition on $L$
implies that $\rho(L) < 1$, so we may use Lemma \ref{re_norm} to
renorm $E$ so that 

\beq \label{norm-L-less-1}
\n{L} < 1,
\eeq

which implies that 

\beq \label{norm-Ln-less-1}
\n{L^n} < \n{L}^n < 1
\eeq

for every positive integer $n$. 

From the fact that $L$ is $\alpha$-contracting, formula (\ref{spectral-radius-formula}) gives
us a positive integer $m$ such that

\[ \left( \n{L^{-m}} \n{L^m}^{1+\alpha} \right)^{\frac{1}{m}} < 1, \]

which, in turn, implies

\[ \n{L^{-m}} \n{L^m}^{1+\alpha} < 1. \]

Given a linear automorphism $L$, we recall that its  {\it minimum
norm}, denoted by $m(L)$ (or $m_L$)  is defined by

\beq \label{min-norm-L}
m(L) = m_L = \inf_{\mid x \mid =1} \n{Lx}.
\eeq

It is easily checked that $m(L) = \n{L^{-1}}^{-1}$. 

Let $\gre_1 > 0$ be such that

\beq \label{L-gre-less-1}
\n{L} + \gre_1 < 1,
\eeq

\beq \label{Linv-gre-1-0}
m(L) - \gre_1 > 0,
\eeq

\beq \label{Lm-gre-less-1}
\n{L^m} + \gre_1 < 1,
\eeq

and 

\beq \label{Linvm-gre-1-0}
m(L^{-m}) - \gre_1 > 0,
\eeq

and 

\beq \label{Lm-gra-est}
\n{L^{-m}} (\n{L^m} + \gre_1)^{1+\alpha} < 1. 
\eeq

Using $Df_0 = 0$, choose $\delta > 0$ small enough so that, for $x \in
B_{\delta}(0)$, 

\beq \label{Df-est-1a}
\n{Df_x} < \gre_1.
\eeq

In particular, 

\beq \label{Lip-f}
Lip(f, B_{\delta}(0)) \leq \gre_1.
\eeq

Now, for $\n{x} < \delta$, 

\beq \label{L_f_est}
\n{T(x)} = \n{Lx + f(x)} \leq (\n{L} + \gre_1)\n{x} \leq \n{x},
\eeq 

and, for $x \neq  y \in B_{\delta}$, we have

\bqns
\n{T(x) -T(y)} & = & \n{(L+f)x - (L+f)y} \\
  & \geq & \n{L(x-y)} - \gre_1\n{x-y} \\
  & > & (m_L -\gre_1)\n{x-y} > 0. 
\eqns

In particular, since $T(0) = 0$, it follows that  $\n{T(x)} > 0$ for $\n{x} > 0$.

Further, 

\bqns
\n{T(x) - T(y)} & = & \n{Lx + f(x) - Ly - f(y)} \\
 & \leq & (\n{L} + Lip(f))\n{x-y} \\
 & \leq & (\n{L} + \gre_1)\n{x-y}, 
\eqns

so,

\beq \label{Lip-T-grd}
 Lip(T,B_{\delta}(0)) \leq \n{L} + \gre_1 < 1.
\eeq

This implies that both $T$ and $T^m$ map $B_{\delta}$ into itself. 

An easy induction shows that $DT^m(0) = L^m$. So, since the
composition of $C^{1,\alpha}$ maps is $C^{1,\alpha}$, if we set $f_m = T^m -
L^m$, then we have that  $f_m \in C^{1,\alpha}_{\delta}$.  

Now, shrinking $\delta$, if necessary, we may assume that 

\beq \label{Lip-fm-1}
Lip(f_m,B_{\delta}) < \gre_1, 
\eeq

and 

\beq \label{Lip-Tm-1}
Lip(T^m,B_{\delta}) < \n{L^m} + \gre_1.
\eeq

For $n=1$ or $n=m$, let

\[ K_{n,1} = \n{L^{-n}} (\n{L^n} + \gre_1)^{\alpha} \delta^{\alpha}
\n{f_n}_{\delta}, \]

and 

\[ K_{n,2} =  \n{L^{-n}}(\n{L^n} + \gre_1)^{1+\alpha},  \]

and 

\[ \grt_n = K_{n,1} + K_{n,2}. \]

From (\ref{Lm-gra-est}), we have $K_{m,2} < 1$, so we can choose $\delta > 0$ small enough so that
$\grt_m < 1$.

Let us now, proceed to the proof of (\ref{H-def-C})-(\ref{norm-Hm}).

Let $\grf \in C^{1,\alpha}_{\delta}$. 

From (\ref{Lip-T-grd}), it follows that $T(B_{\delta}(0)) \subset
B_{\delta}(0)$, so  $H(\grf)$ is defined on
$B_{\delta}(0)$, and this is (\ref{H-def-C}). 

Again using that the composition of $C^{1,\alpha}$ maps is again
$C^{1,\alpha}$, it is clear that
$H(\grf)$ is  $C^{1,\alpha}$. It is also evident that it 
vanishes at $0$ together with its derivative. 

For $x, y \in B_{\delta}$ with $\n{x - y} > 0$ and $\n{x} > 0$, and
$n=1$ or $n=m$, we have

\bqns
\n{D(H^n(\grf))(x) -D(H^n(\grf))(y)} & = & \n{L^{-n}D\grf_{T^n(x)}DT^n_x -L^{-n}D\grf_{T^n(y)}DT^n_y} \nonumber \\
 & \leq &  \n{L^{-n}} \n{D\grf_{T^n(x)}DT^n_x -D\grf_{T^n(x)}DT^n_y} \\
 & & \ \ \ + \n{L^{-n}}  \n{D\grf_{T^n(x)}DT^n_y -D\grf_{T^n(y)}DT^n_y} \\
 & = & (A) + (B)
\eqns

where 

\bqns
(A) & = & \n{L^{-n}}\n{D\grf_{T^n(x)}DT^n_x -D\grf_{T^n(x)}DT^n_y} \\
   & \leq & \n{L^{-n}}\n{D\grf_{T^n(x)}} \mid L^n + Df_{n,x} -(L^n +
Df_{n,y}) \mid \\
  & = & \n{L^{-n}}\n{D\grf_{T^n(x)}} \mid Df_{n,x} -Df_{n,y} \mid
\eqns

and 

\bqns
 (B) & = & \n{L^{-n}}  \n{D\grf_{T^n(x)}DT^n_y -D\grf_{T^n(y)}DT^n_y} \\
    & \leq & \n{L^{-n}} (\n{L^n}+\gre_1) \n{D\grf_{T^n(x)} - D\grf_{T^n(y)}}.
\eqns

For $x \neq 0$, we have 

\bqns
(A)  & \leq & \n{L^{-n}} \frac{\n{D\grf_{T^n(x)}}}{\n{T^n(x)}^{\alpha} } \n{T^n(x)}^{\alpha} 
\n{f_n}_{\delta} \n{x-y}^{\alpha} \\
 & \leq & \n{L^{-n}} \,  \n{\grf}_{\delta} \, (\n{L^n}+\gre_1)^{\alpha}
\n{x}^{\alpha} \n{f_n}_{\delta} \n{x-y}^{\alpha} \\
 & \leq & \n{L^{-n}} \, \n{\grf}_{\delta} \, (\n{L^n}+\gre_1)^{\alpha}
\delta^{\alpha} \n{f_n}_{\delta} \n{x-y}^{\alpha} \\
& =  & K_{n,1} \, \n{\grf}_{\delta} \, \n{x-y}^{\alpha}. 
\eqns

and, for $x \neq y$, 

\bqns
(B) & \leq & \n{L^{-n}} (\n{L^n}+\gre_1) \n{D\grf_{T^n(x)} - D\grf_{T^n(y)}} \\
   & \leq & \n{L^{-n}} (\n{L^n}+\gre_1) \frac{\n{D\grf_{T^n(x)} -
D\grf_{T^n(y)}}}{\n{T^n(x) - T^n(y)}^{\alpha}} \n{T^n(x) - T^n(y)}^{\alpha} \\
   & \leq & \n{L^{-n}} (\n{L^n}+\gre_1) \, \n{\grf}_{\delta} \, \n{T^n(x) -
T^n(y)}^{\alpha} \\
   & \leq & \n{L^{-n}} (\n{L^n}+\gre_1) \, \n{\grf}_{\delta} \,
Lip(T^n)^{\alpha} \n{x-y}^{\alpha} \\
   & \leq & \n{L^{-n}} (\n{L^n}+\gre_1)^{1+\alpha}
\,  \n{\grf}_{\delta} \, \n{x - y}^{\alpha} \\
 & =  & K_{n,2}  \, \n{\grf}_{\delta} \, \n{x - y}^{\alpha}.
\eqns

Hence, 
we have

\bqn 
\n{D(H^n(\grf))(x) -D(H^n(\grf))(y)} & \leq & (A) + (B) \nonumber \\
 & \leq & \left(K_{n,1}  \n{\grf}_{\delta}  + K_{n,2} \n{\grf}_{\delta}
\right) \,  \n{x - y}^{\alpha} \nonumber \\
 & = & \grt_n \n{\grf}_{\delta} \,  \n{x - y}^{\alpha}. \label{norm-DH}
\eqn

Dividing both sides by $\n{x - y}^{\alpha}$ and taking the supremum
over $x \neq y$, gives

\[ \n{H^n(\grf)}_{\delta} \leq \grt_n \n{\grf}_{\delta}. \]

For $n=1$, this gives (\ref{H-into-C}), and for $n=m$, this gives
(\ref{norm-Hm}), completing the proof of Theorem \ref{Lin_Contraction_2}.

\section{Preliminary Constructions} \label{Preliminary_Constructions}

Let us begin with an outline of Hartman's method in the proof of part (b)
of Theorem \ref{Hartman-Theorem}. He first  changes coordinates so that the 
stable and unstable manifolds are flattened (i.e., are contained in
the stable and unstable subspaces, respectively) near the origin. Next, using his
earlier result for $C^{1,1}$ contracting $T$, he chooses
$C^{1,\grb}$ coordinates (where $\grb$ is related to the contracting
and expanding eigenvalues) so that, near the origin,  $T$ preserves
and becomes linear when restricted to the
stable and unstable subspaces. Finally, he globalizes the mapping $T$,
using a bump function to replace $T$ by a map $T_1$
which equals $T$ on a ball $B_{\gre_0}(0)$ about zero,
equals $L$ off a slightly larger ball $B_{\gre_1}(0)$, and is globally 
Lipschitz close to $L$. This globalization preserves the earlier
properties that the stable and unstable manifolds are flattened and
$T_1$ is linear when restricted to the entire stable and unstable
subspaces. 
This allows him to get sufficiently good estimates on the partial
derivatives of the non-linear part of $T_1$ to use the method of successive
approximations to obtain the desired local linearization. In fact, he only
proves that the linearization is continuous. He states that similar
methods will give that it is $C^1$ and that the proof can be modified
to show that it has H{\"{o}}lder continuous derivatives.  

We extend and modify his techniques to obtain our result.  In addition, we write
our solution in a way that it makes use  of certain
linear  contracting maps in suitable Banach spaces. This
provides the additional benefit that we can get continuous dependence
of the linearization on parameters as in Theorem \ref{par-Main-Theorem} below.

\subsection{Estimates of the Lipschitz and D-H{\"{o}}lder contants
of inverses and compositions}\label{Hol-nonlinear}

The proof of Theorem \ref{Main_Theorem_1} requires estimates of the non-linear parts
of the diffeomorphisms $T, T^m, T^{-1}$ and $T^{-m}$ for a certain
positive integer $m$ in a sufficiently small ball $B_{\delta}(0)$
about $0$ in the Banach space $E$.  The results in this section can
be used to give some information about these estimates in terms of the
non-linear part of the original map $T$.  Strictly speaking, they are
not needed if one only wants the linearization $R$ to exist on 
{\it some} small neighborhood of $0$ and one does not need to 
estimate the size of that neighborhood. 

Let us begin with a simple lemma relating the Lipschitz and
D-H{\"{o}}lder constants of a composition $S \circ T$ of  $C^{1,\alpha}$ maps in
terms of those of the maps $S$ and $T$.

\begin{Lem} \label{Comp-Lemma}
Let $U$ and $V$ be open subsets of the Banach space $E$, and consider
maps $S \in C^{1,\alpha}(U,E)$ and $T \in C^{1,\alpha}(V,E)$. 

Then,

\beq \label{Lip-Comp}
Lip(S \circ T, U \bigcap T(V)) \leq Lip(S,U) Lip(T,V),
\eeq

and 

\beq \label{DHol-Comp}
Hol(D(S \circ T),U \bigcap T(V)) \leq Lip(S,U) Hol(DT,V) + Lip(T,V)^{1+\alpha}
Hol(DS,U). 
\eeq

\end{Lem}

\begin{proof}

Letting  $x, y \in V$, and leaving out the obvious domains of the maps
involved, statement (\ref{Lip-Comp}) is immediate from

\bqns
 \n{S(T(x)) -S(T(y))} & \leq &  Lip(S)\n{T(x) - T(y)} \\
 & \leq & Lip(S) Lip(T) \n{x-y}.
\eqns

For statement (\ref{DHol-Comp}), we have 

\bqns
\n{ D(S \circ T)(x) - D(S \circ T)(y)} & = & \n{DS_{Tx}DT_x -
DS_{Ty}DT_y} \\
 & \leq & \n{DS_{Tx}DT_x - DS_{Ty}DT_x} + \n{DS_{Ty}DT_x -
DS_{Ty}DT_y} \\
 & \leq & \n{DT_x} Hol(DS)\n{Tx - Ty}^{\alpha} +
\n{DS_{Ty}}Hol(DT)\n{x - y}^{\alpha} \\
 & \leq & Lip(T)^{1+\alpha}Hol(DS)\n{x-y}^{\alpha} +
   Lip(S) Hol(DT)\n{x - y}^{\alpha}.
\eqns

Now, divide both sides by $\n{x-y}^{\alpha}$ and take the supremum over
$x,y$ to complete the proof of Lemma \ref{Comp-Lemma}. 

\end{proof}

\begin{Rem} \label{Hol-DT-Df}
Observe that, if $T$ has the form $T = L+f$ where $L$ is bounded,
$f(0) = 0$, and $Df(0) = 0$, then, since the derivative of a bounded
linear map $L$ is just $L$ at every point, it cancels in the calculation of
$DT$. Hence, 

\beq \label{Hol-DT}
Hol(DT) = Hol(Df).
\eeq

That is, the D-H{\"{o}}lder constant is determined by the nonlinear
part of $T$.  In particular, the addition of another bounded linear
map to $T$ does not change $Hol(DT)$. 
\end{Rem}

\begin{Cor} \label{Lip-Hol-Df3}
Let $S, T, U, V$ be as in the hypotheses of Lemma \ref{Comp-Lemma},
and assume that $S = L_1 + f_1$ and $T = L_2 + f_2$ where $L_1, L_2$ are
bounded linear maps on $E$ and $f_1, f_2$ vanish at $0$ together with
their derivatives. Let $f_3 = S \circ T - L_1L_2$. 

Then, 

\beq \label{Lip-f3}
Lip(f_3) \leq \n{L_1}Lip(f_2) + Lip(f_1)Lip(T), 
\eeq

and 

\beq \label{Hol-Df3}
Hol(Df_3)  \leq  Lip(S) Hol(DT) + Lip(T)^{1+\alpha} Hol(Df_1). 
\eeq

\end{Cor}

\begin{proof}

We have 

\bqns
 S T & = & (L_1 + f_1) (L_2 + f_2) \\
     & = & L_1(L_2 + f_2) + f_1(T) \\
     & = & L_1L_2 + L_1f_2 + f_1(T),
\eqns

giving  $ f_3 = L_1f_2 + f_1(T)$, and (\ref{Lip-f3}) immediately
follows using the fact that the function $Lip(\cdot)$ is subadditive and submultiplicative.

Statement (\ref{Hol-Df3}) follows immediately from (\ref{DHol-Comp})
since 

\[ Hol(Df_3) = Hol(D(S \circ T)). \]

\end{proof}

We will need the following  estimate which we first saw in Hirsch and Pugh
\cite{Hirsch-Pugh}.

\begin{Lem} \label{HP-lin}
For any linear automorphisms $h_1, h_2$, we have

\beq \label{Lip-lin-inv}
\n{h_1^{-1} - h_2^{-1}} \leq \n{h_2^{-1}}\n{h_1 - h_2}\n{h_1^{-1}}. 
\eeq

\end{Lem}

\begin{proof}

We have 

\bqns \n{h_1^{-1} - h_2^{-1}} & = & \n{h_2^{-1}h_2 h_1^{-1} - h_2^{-1}h_1 
  h_1^{-1}} \\
  & = & \n{h_2^{-1}(h_2 - h_1)h_1^{-1}} \\
 & \leq & \n{h_2^{-1}}\n{h_2   - h_1}\n{h_1^{-1}} \\
 & = & \n{h_2^{-1}}\n{h_1   - h_2}\n{h_1^{-1}} 
\eqns

\end{proof}

Let $diam(U)$ denote the diameter of a subset $U$ of a Banach space
$E$. 

In Lemma \ref{Lip-g-Hol-g}, we again use the notation

\[ m(L) = m_L = \inf_{\mid x \mid = 1} \n{Lx} = \n{L^{-1}}^{-1} \]

for a linear automorphism $L$.

\begin{Lem} \label{Lip-g-Hol-g}
Let $(E, \n{\cdot})$ be a real Banach space, and let $L:E \rarrow E$
be a linear automorphism. 
Let $0 < \alpha < 1$, and let $U$ and $V$ be neighborhoods of $0$ in $E$ such
that 

\beq \label{diam}
\max(diam(U), diam(V)) < 1.
\eeq

Let $T: U \rarrow V$ be a $C^{1+\alpha}$
diffeomorphism from $U$ onto $V$ such that $T(0) = 0$ and $DT(0) =
L$, and let $f = T - L$ and $g =  T^{-1} - L^{-1}$ be the nonlinear
parts of $T$ and $T^{-1}$, respectively. 

Assume that 

\beq \label{Df-1a-est}
Lip(f,U) = \sup_{z \in U} \n{Df(z)} < m_L. 
\eeq

Then, we have 

\beq \label{Lip-Tinv}
Lip(T^{-1},V) \leq (m_L - Lip(f,U))^{-1}, 
\eeq

\beq \label{Lip-f-1}
Lip(g,V) \leq \n{L^{-1}} (m_L - Lip(f,U))^{-1} Lip(f,U),
\eeq

\beq \label{Hol-DTinv}
Hol(DT^{-1},V) \leq (m_L-Lip(f,U))^{-(2+\alpha)} Hol(Df,U)
\eeq

and 

\beq \label{Hol-Dg-2}
Hol(Dg,V) \leq (m_L-Lip(f,U))^{-(2+\alpha)} Hol(Df,U).
\eeq

\end{Lem} 

\begin{proof}

For notational convenience, we will use the following notation.

\[ l_f = Lip(f,U), \ h_f = Hol(Df,U). \]

From (\ref{Df-1a-est}), we have that, for any $z \in U$ and non-zero vector $v$,

\[ \n{DT_zv}  =   \n{Lv + Df(z)v}  \geq  m_L \n{v} - l_f \n{v} = (m_L
- l_f) \n{v}. \]

Now, taking any non-zero $u$ and setting $v = DT^{-1}_{T(z)}u$ and $w
= T(z)$, we get 

\[ \n{u} = \n{DT_z DT^{-1}_w u} \geq (m_L -l_f) \n{DT^{-1}_wu}, \]

or 

\[ \n{DT^{-1}_w u} \leq (m_L - l_f)^{-1} \n{u}. \]

Thus, for every $w \in V$, we have 

\[ \n{DT^{-1}_w} \leq (m_L - l_f)^{-1}, \] 

and (\ref{Lip-Tinv}) follows. 

Again, writing 

\[ w = T(z) = Lz + f(z), \]

and solving for $z = z(w)$, we get 

\beq \label{z-inv}
z(w) = T^{-1}(w) = L^{-1}w - L^{-1}f(z(w)).
\eeq

Thus, 

\[ g(w) = -L^{-1}f(z(w)), \]

\[ Dg(w) = -L^{-1} Df_{z(w)}DT^{-1}_w, \]

and 

\bqns
 \n{Dg(w)} & \leq & \n{L^{-1}} \n{Df_{z(w)}} \n{DT^{-1}_w}   \\
 & = &  \n{L^{-1}} \n{Df_{z(w)}}(m_L - l_f)^{-1} \\
 & \leq &  \n{L^{-1}} l_f (m_L - l_f)^{-1} \\
\eqns 

This implies (\ref{Lip-f-1}).

Proceding to the proof of (\ref{Hol-DTinv}), from (\ref{z-inv}), we
have 

\[ DT^{-1}_w = L^{-1} - L^{-1} Df_{z(w)} DT^{-1}_w \]

\[ (I + L^{-1}Df_{z(w)}) DT^{-1}_w = L^{-1}, \]

\[ DT^{-1}_w = (I + L^{-1}Df_{z(w)})^{-1} L^{-1}. \]

From (\ref{Df-1a-est}), we have, for each $w \in V$,

\bqns
\n{(I + L^{-1}Df_{z(w)})^{-1}} & \leq & (1 - \n{L^{-1}}l_f)^{-1} \\
    &  =   &  (\n{L^{-1}}m_L - \n{L^{-1}}l_f)^{-1} \\
    &  =   &  \n{L^{-1}}^{-1}(m_l - l_f)^{-1} \\
    &  =   &  m_L (m_l - l_f)^{-1} \\
\eqns

This, together with (\ref{Lip-lin-inv}) gives 

\bqns
\n{DT^{-1}_{w_1}- DT^{-1}_{w_2}} & \leq & m_L^2(m_L - l_f )^{-2} 
 \n{L^{-1}}^2 \n{Df_{z(w_2)} - Df_{z(w_1)}} \\
& = & (m_L - l_f )^{-2} \n{Df_{z(w_2)} - Df_{z(w_1)}} \\
& \leq &  (m_L - l_f)^{-2} h_f \n{z(w_2) - z(w_1)}^{\alpha} \\
& \leq &  (m_L - l_f)^{-2} h_f l(T^{-1})^{\alpha}\n{w_2 - w_1}^{\alpha} \\
& \leq &  (m_L - l_f)^{-2} h_f (m_L - l_f)^{-\alpha}\n{w_2 - w_1}^{\alpha} \\
& \leq &  (m_L - l_f)^{-2-\alpha} h_f \n{w_2 - w_1}^{\alpha} 
\eqns

which is (\ref{Hol-DTinv}).

Also, (\ref{Hol-Dg-2}) holds since $Hol(DT^{-1},W) = Hol(Dg,W)$.
This completes the proof of Lemma \ref{Lip-g-Hol-g}. 

\end{proof}

The next lemma gives analogous estimates of the non-linear parts of $T^m$ and
$T^{-m}$ for integers $m > 1$. 

\begin{Lem} \label{Lip-Hol-Power}
Let $T, U, V, f, g$ be as in Lemma \ref{Lip-g-Hol-g}, 
let $m$ be a positive integer, and assume that $Lip(T) \geq 1$. 

Consider the maps $T^m, T^{-m}$ and $f_m, g_m$, given by 

\[ f_m = T^m - L^m, \ g_m = T^{-m} - L^{-m}. \]

with corresponding domains $U_m, V_m$ given by 

\[ U_m = \bigcap_{j=0}^{m}T^{-j}U \mbox{ and } V_m = \bigcap_{j=0}^{m}T^jV, \]
respectively. 

Then, $T^m$ and $T^{-m}$ are defined and $C^{1,\alpha}$ on $U_m$ and $V_m$,
respectively, and 

\beq \label{Lip-Tm}
Lip(T^m) \leq Lip(T)^m, 
\eeq

\beq \label{Lip-fm}
 Lip(f_m) \leq m Lip(f) Lip(T)^{m-1}, 
\eeq

\beq \label{Hol-fm}
 Hol(Df_m) \leq  m Hol(Df) Lip(T)^{(1+\alpha)(m-1)}, 
\eeq

\beq \label{Lip-gm}
 Lip(g_m) \leq mLip(g) Lip(T^{-1})^{m-1}, 
\eeq

and 

\beq \label{Hol-gm}
 Hol(Dg_m) \leq m Hol(Dg) Lip(T^{-1})^{(1+\alpha)(m-1)}.
\eeq

\end{Lem}

\begin{proof}

Since the composition of $C^{1,\alpha}$ maps is again $C^{1,\alpha}$, it
is clear that $T^m$ and $T^{-m}$ are defined and $C^{1,\alpha}$ on $U_m$
and $V_m$, respectively. 

Also, the Lipschitz composition formula (\ref{Lip-Comp}) gives (\ref{Lip-Tm})
immediately by induction. 

Using that $DT^n(0) = L^n$, we get 

\bqns
T^{n+1}  & = & (L+f)(T^n) \\
    & = & LT^n + f(T^n) \\
    & = & L(L^n + f_n) + f(T^n)\\
    & = & L^{n+1} + Lf_n + f(T^n).
\eqns

This implies that, for each $n \geq 1$, 

\beq \label{recursive-fn}
f_{n+1} = Lf_n + f(T^n).
\eeq

We will show that, for different choices of positive numbers $A_n,
B, C$, the numbers  $Lip(f_m), Lip(g_m), Hol(Df_m)$ and $Hol(Dg_m)$ all
satisfy the recursion relation

\beq \label{rec-abc}
A_{n+1} \leq B A_n + B^n C, \ A_1 = C
\eeq

for each $1 \leq n \leq m-1$. 

This recursion relation has the solution 

\beq \label{rec-soln}
A_{n+1} \leq (n+1) B^nC 
\eeq

as can be seen from 

\bqn
 A_{n+1} & \leq & B A_n + B^n C \nonumber \\
 & \leq & B(BA_{n-1} + B^{n-1}C) + B^n C \nonumber \\
 & = & B^2A_{n-1} + B^nC + B^n C \nonumber \\
 & = & B^2A_{n-1} + 2 B^n C \nonumber \\
 & \vdots & \nonumber \\ 
 & \leq & B^nA_1 + n B^n C \nonumber \\
 & \vdots &  \nonumber \\
 & = & B^n A_1 + n B^n C \nonumber \\
 & = & (n+1) B^nC \label{abc-rec}.
\eqn

Now, consider $Lip(f_m)$.

Since $\n{L} \leq Lip(T)$, we have

\bqns
 Lip(f_{n+1}) & \leq  & \n{L} Lip(f_n) + Lip(f) Lip(T^n) \\
              & \leq & Lip(T) Lip(f_n) + Lip(f) Lip(T^n) \\
\eqns

and we can take $A_n = Lip(f_n), B = Lip(T)$ and $C = Lip(f)$. 

Next, observe that $Hol(Df_m) = Hol(DT^m)$ since 

\[  \n{DT^mx - DT^my}  =  \n{L^m + Df_m(x) - (L^m + Df_m(y))} =
\n{Df^mx - Df^my}.  \] 

The assumption that $Lip(T) \geq 1$ implies that 

\[ Lip(T) \leq Lip(T)^{1+\alpha}. \]

Using this and (\ref{Hol-Df3}),  we have 

\bqns
 Hol(DT^{n+1}) & = & Hol(D(T \circ T^n)) \\
               & \leq & Lip(T) Hol(DT^n) + Lip(T^n)^{1+\alpha} Hol(DT)
\\
               & \leq & Lip(T)^{1+\alpha} Hol(DT^n) + Lip(T)^{n(1+\alpha)} Hol(DT)
\eqns

and we can take $A_n = Hol(DT^n), B = Lip(T)^{1+\alpha}$ and $C =
Hol(DT)$. 

Similar recursions hold for $Lip(g_m)$ and $Hol(Dg_m)$. 

Thus, formula (\ref{rec-soln}) applied, in turn,  to each of $Lip(f_m),
Hol(Df_m), Lip(g_m)$ and $Hol(Dg_m)$ gives (\ref{Lip-fm})-(\ref{Hol-gm}) to complete
the proof of Lemma \ref{Lip-Hol-Power}.

\end{proof}

\subsection{Global Extension Via Bump Functions} \label{Bump-Functions}

Let $(E,\n{\cdot})$ be a $C^{1,\alpha}$ Banach space. Thus, there is a
$C^{1,\alpha}$ real valued function $\grl$ defined on $E$ and a 
real number $c \in (0,1)$ such that 

\[ \lambda(E) = [0,1], \]

\[ \grl(x) = 1 \mbox{ for } \n{x} \leq c, \]

\[ \grl(x) = 0 \mbox{ for } \n{x} \geq 1, \]

\[
 Lip(\lambda,E) = \sup_{x \in E} \n{D\lambda(x)}  < \infty, 
\]

and 

\[
Hol(D\lambda,E) = \sup_{x \neq y \in E} \frac{ \n{D\lambda(x) -
D\lambda(y)}}{\n{x-y}^{\alpha}}  < \infty.
\]

We call $\grl$ a {\it unit bump function } on $E$. 

Given such a bump function $\grl$ and a positive 
real number $\delta$, we define the associated $\delta-scaled$ version of $\grl$ to be 

\beq \label{delta-scaled}
 \grl_{\delta}(x) = \grl(\frac{x}{\delta}). 
\eeq

Obviously, the function $\grl_{\delta}$ vanishes off $B_{\delta}$ and
is 1 on $B_{c\delta}$. 

The Lipschitz and H{\"{o}}lder constants of $\grl_{\delta}$ are easily
computed as follows.  

We have

\[ D\grl_{\delta}(x) = \frac{1}{\delta} D\lambda(\frac{x}{\delta}), \]

and

\bqns
 \n{D\grl_{\delta}(x) - D\grl_{\delta}(y)} & = & \n{\frac{1}{\delta}( D\lambda(\frac{x}{\delta}) - D\lambda(\frac{y}{\delta}))} \\
 & \leq & \frac{1}{\delta}Hol(D\lambda,E) \n{\frac{x}{\delta} -   \frac{y}{\delta}}^{\alpha} \\
 & = & \frac{1}{\delta^{1+\alpha}} Hol(D\lambda,E) \n{x - y}^{\alpha}
\eqns

for any $x \neq y \in E$.

Thus, 

\beq \label{Lip-Hol-est1}
Lip(\grl_{\delta}) = \frac{Lip(\lambda)}{\delta} \mbox{ and } Hol(D \grl_{\delta}) =
\frac{Hol(D\lambda)}{\delta^{1+\alpha}}
\eeq

Recall that, for a $C^{1,\alpha}$ function $f:U \rarrow E$, we have defined

\[ Lip(f,U) = \sup_{x \in U} \n{Df(x)} \]

and

\[ Hol(Df,U) = \sup_{x \neq y \in U} \frac{\n{Df(x) - Df(y)}}{\n{x -
    y}^{\alpha}}. \]

For functions $g$ defined on all of $E$, we leave out the domain $E$ and simply write

\[ Lip(g) = Lip(g,E) \mbox{ and } Hol(Dg) = Hol(Dg,E). \]

\begin{Lem} \label{Hol-Dg}
Let $U$ be an open neighborhood of $0$ in the $C^{1,\alpha}$ Banach space $E$
and let $f:U \rarrow E$ be a $C^{1,\alpha}$ map such that $f(0) = 0$ and
both $Lip(f,U)$ and $Hol(Df,U)$ are finite. 

Let $\grl:E \rarrow \R$ be the unit bump function defined above, and
define the functions  $C_1(\grl)$ and $C_2(\grl)$ by 

\[ C_1(\grl) = 1 + Lip(\grl) \]

and 

\[ C_2(\grl) = 1 + 2 Lip(\grl) + 3 Hol(D\grl). \]

For any $\delta > 0$ be such that $B_{\delta}(0) \in U$, consider the
function $g$ defined by the product 

\[ g(x) = \grl_{\delta}(x) f(x) \]
 
where $\grl_{\delta}$ is the $\delta-$scaled version of $\grl$ defined
in (\ref{delta-scaled}). 

Then, for $c$ as in the definition of $\grl(\cdot)$, the function $g$ is defined and $C^{1,\alpha}$ on all of $E$ and
satisfies the following.

\beq \label{g=f}
g(x) = f(x) \mbox{ for } \n{x} \leq c \ \delta, 
\eeq

\beq \label{g=0}
g(x) = 0 \mbox{ for } \n{x} \geq \delta, 
\eeq

\beq \label{Lipg}
Lip(g) \leq C_1(\grl) Lip(f,U) 
\eeq

and

\beq \label{Hol-Dg-1}
Hol(Dg) \leq C_2(\grl) Hol(Df,U) 
\eeq

 \end{Lem}

{\bf Proof.}
Statements (\ref{g=f}) and (\ref{g=0}) are  obvious.  

For statement (\ref{Lipg}),  using $\n{x} \leq \delta$ and $f(0)=0$, we have

\bqns
\n{D(\grl_{\delta}(x) f(x))} & = & \n{ D \grl_{\delta}(x) f(x) +
  \grl_{\delta}(x) Df(x)} \\
 &  \leq & \frac{Lip(\grl)}{\delta} \n{f(x)} + \n{\grl_{\delta}(x)}
\n{Df(x)} \\
 & \leq & \frac{Lip(\grl)}{\delta} Lip(f,U) \n{x} + Lip(f,U) \\
 & \leq & Lip(\grl)Lip(f,U) + Lip(f,U) \\
 & = & Lip(f,U) ( Lip(\grl) + 1) 
\eqns

as needed.

We now proceed to verify (\ref{Hol-Dg-1}). 

For ease of notation, let us denote $\grl_{\delta}(x)$ as $\grg(x)$. 

It suffices to prove that, for $x \neq y \in U$, 

\beq \label{Dgrg-f}
\n{D(\grg f)(x) - D(\grg f)(y)} \leq C_2(\grl) Hol(Df,U) \n{x - y}^{\alpha}. 
\eeq

We have 

\bqns
 \n{D(\grg f)(x) - D(\grg f)y} & = & \n{D\grg(x) f(x) + \grg(x) Df(x)
   - (D\grg(y) f(y) + \grg(y) Df(y)) } \\
 & \leq & \n{D\grg(x) f(x) -D\grg(y) f(y)} + \n{\grg(x) Df(x) -\grg(y) Df(y)} 
\eqns
 
Let us estimate the two expressions

\[ R_1 = \n{D\grg(x) f(x) -D\grg(y) f(y)} \]

and

\[ R_2 = \n{\grg(x) Df(x) -\grg(y) Df(y)} \]

separately.

Let $H_0(Df)$ represent the H{\"{o}}lder constant of $Df$ at $0$. 

We have

\bqns
R_1 & = & \n{D\grg(x) f(x) -D\grg(y) f(y)} \\
    & \leq & \n{D\grg(x) f(x) -D\grg(y) f(x)} +  \n{D\grg(y) f(x)
  -D\grg(y) f(y)} \\
    & \leq & Hol(D\grg) \n{x - y}^{\alpha} H_0(Df) \n{x}^{1 + \alpha}
\\
 & & \ \ \ + H_0(D\grg)\n{y}^{\alpha} \sup_{0 \leq \grt \leq 1} \n{Df((1-\grt)x + \grt
  y)} \n{x-y} \\
& \leq & \frac{Hol(D\grl)}{\delta^{1+\alpha}}
H_0(Df)\n{x-y}^{\alpha}\delta^{1+\alpha} \\
 & & \ \  + \frac{Hol(D\grl)}{\delta^{1+\alpha}}\n{\delta}^{\alpha}H_0(Df) \max(\n{x}, \n{y})^{\alpha} \n{x-y} \\
 & \leq & Hol(D\grl) H_0(Df) \n{x-y}^{\alpha} + \frac{Hol(D\grl)}{\delta}H_0(Df) 
\delta^{\alpha} \n{x-y}^{1-\alpha} \n{x-y}^{\alpha} \\
& \leq & Hol(D\grl) H_0(Df) \n{x-y}^{\alpha} + \frac{Hol(D\grl)}{\delta}H_0(Df) 
\delta^{\alpha} (2 \delta)^{1-\alpha} \n{x-y}^{\alpha} \\ 
& \leq & Hol(D\grl) H_0(Df) \n{x-y}^{\alpha} + Hol(D\grl)H_0(Df) 
 2  \n{x-y}^{\alpha} \\
& \leq & 3 Hol(D\grl) H_0(Df) \n{x-y}^{\alpha} \\
& \leq & 3 Hol(D\grl) Hol(Df,U) \n{x-y}^{\alpha} 
\eqns

and 

\bqns
R_2 & = &  \n{\grg(x) Df(x) -\grg(y) Df(y)} \\
    & \leq & \n{\grg(x) Df(x) -\grg(y) Df(x)} + \n{\grg(y) Df(x)
  -\grg(y) Df(y)} \\
& \leq & Lip(\grg) \n{x-y} \n{Df(x)} + \n{Df(x) - Df(y)} \\
& \leq & \frac{Lip(\grl)}{\delta}\n{x-y}^{\alpha}\n{x-y}^{1-\alpha}Hol(Df,U)
\n{x}^{\alpha} + Hol(Df,U) \n{x-y}^{\alpha} \\
& \leq & \frac{Lip(\grl)}{\delta}\n{x-y}^{\alpha}2^{1-\alpha}\max(\n{x},\n{y})^{1-\alpha}Hol(Df,U)
\delta^{\alpha} + Hol(Df,U) \n{x-y}^{\alpha} \\
& \leq & \frac{Lip(\grl)}{\delta}\n{x-y}^{\alpha}2\delta^{1-\alpha}Hol(Df,U)
\delta^{\alpha} + Hol(Df,U) \n{x-y}^{\alpha} \\
& \leq & (2 Lip(\grl) + 1) Hol(Df,U) \n{x-y}^{\alpha} 
\eqns 

Hence,

\bqns
 \n{D(\grg f)(x) - D(\grg f)(y)} & = & R_1 + R_2 \\
 & \leq & (3 Hol(D\grl) + 2 Lip(\grl) +1) Hol(Df,U) \n{x-y}^{\alpha} \\
 & \leq & C_2(\grl) Hol(Df,U) \n{x-y}^{\alpha} \\
\eqns

This completes the proof of Lemma \ref{Hol-Dg}.

\subsection{Convention when using direct sum decompositions.}

When the Banach space $E$ is written as a direct sum decomposition
$E^u \oplus E^s$, it is often convenient to identify $T:E \rarrow E$
with the  map
$\tilde{T}$ from $E^u \times E^s$ to $E^u \times E^s$ defined by
taking the unique representations of $z = x+y, T(z) = x_1+y_1$, with
$x, x_1 \in E^u, \ y,y_1 \in E^s$, 
and defining 

\[ \tilde{T}(x,y) = (x_1, y_1). \]

The map $\tilde{T}$ is  the conjugate $RTR^{-1}$ where $R(z) =
(x,y)$. Thus, $\tilde{T}$ is simply the map obtained using $R$ as a
linear change of coordinates. 

  Some statements have a more elegant formulation using $T$, but their proofs are best given in
terms of the map $\tilde{T}$. 

We call the map $\tilde{T}$ the {\it product representation } of $T$.
Letting $\pi^u:E^u \oplus E^s \rarrow E^u$ and $\pi^s: E^u \oplus E^s
\rarrow E^s$ denote the natural projections 

\[ \pi^u(z) = x, \ \pi^s(z) = y, \]

and writing 

\[ (\pi^u \circ T)(x,y) = f_1(x,y), \]

and 

\[ (\pi^s \circ T)(x,y) = f_2(x,y), \]

we will identify $T$ with the map $\tilde{T}$ and simply say {\it we
may write $T$ as } 

\[ T(x,y) = (f_1(x,y), f_2(x,y)). \]

Similarly, the origin will be written as $0$ or $(0,0)$, and balls centered at $0$ in $E$ will be written as $B_{\delta} =
B_{\delta}(0)$ or $B_{\delta}(0,0)$. 

\subsection{Flattening and Linearizing on Invariant Manifolds} \label{Flattening-Subsection}
Let $E$ be a real Banach space. 
For a positive real number $\delta$ we use the notation $B_{\delta} = B_{\delta}(0)$
for the open ball of radius $\delta$ centered at $0$; i.e., the set of
points $x  \in E$ such that $\n{x} < \delta$.  We will say that a property holds {\it near } $0$ if it
holds in $B_{\delta}(0)$ for some small $\delta > 0$. 

In this section, all neighborhoods of $0$ in $E$ will be assumed to be
open and convex.  If $U$ is such a neighborhood and $f:U \rarrow E$,
then the Mean Value Theorem can be applied to show that 

\[ Lip(f,U) = \sup_{x \in U} \n{Df_x}, \]

and  we will often use this fact. 

Given a neighborhood $U$ of $0$ in $E$ and an injective  map $T:U \rarrow E$, 
define 

\beq \label{Ws_U}
W^s_U = W^s_U(T) = \bigcap_{n \geq 0} T^{-n}(U)
\eeq

and 

\beq \label{Wu_U}
W^u_U  = W^u_U(T) =  \bigcap_{n \geq 0} T^n(U). 
\eeq

These sets are called the $U$-stable and $U$-unstable
sets of $T$, respectively. 

From the definitions, it is immediate that 

\[  W^s_U(T) = W^u_U(T^{-1}), \ W^u_U(T) = W^s_U(T^{-1}) \]

\[ T(W^s_U(T)) \subset W^s_U(T), \mbox{ and }
T^{-1}(W^u_U(T)) \subset W^u_U(T). 
\] 

When $U$ is a ball $B_{\delta}(0)$, we write

\[ W^s_U(T) = W^s_{\delta}(0,T), \ W^u_U(T) = W^u_{\delta}(0,T). \] 

Let $r \geq 1$ be a real number. A  {\it local $C^r$ diffeomorphism} at $0$ is a pair $(T,U)$ such that
$U$ is an open neighborhood of $0$ in $E$ and $T$ is a $C^r$
diffeomorphism from $U$ onto its image such that $T(0) =0$. Note that
$T(U)$ is an open set by the Inverse Function Theorem. As usual,
we often ignore the domain of the local diffeomorphism, and simply say
that $T$ is a local diffeomorphism at $0$.  When considering
compositions and inverses of  local $C^r$ diffeomorphisms, we simply shrink domains
as needed to make the definitions correct. When $U = E$, we sometimes
say that $T$ is a {\it global} diffeomorphism to emphasize that we are
considering a bijection from $E$ onto $E$.  Since we only consider
 diffeomorphisms with a fixed point at  $0$ in this and the next two
sections, we sometimes drop the {\it
at $0$} and simply use the terms {\it local diffeomorphism} or {\it
   global diffeomorphism}. 

Two local $C^r$ diffeomorphisms $T$ and $S$ are $C^r$
conjugate if there is a  local $C^r$ diffeomorphism $R$ such that
$RTR^{-1} = S$ near $0$.  If  $S$ happens to be linear, then, setting $L =
DT_0, M = DR_0$, and $P = M^{-1}R$, we have

\[ PTP^{-1} = M^{-1}S M, \ DP(0) = I, \mbox{ and }  D(PTP^{-1}) = L. \]

Hence, 

\begin{quote} \label{DR0=I}
{\it $T$ is $C^r$ conjugate to {\bf some}  linear map at $0$ if and only if
it is $C^r$ conjugate to its derivative at $0$, and it may be assumed
that the conjugacy $R$ has the properties that $R(0) = 0$ and $DR_0 =
I$. 
}
\end{quote} 

To emphasize this concept, we will say that $T$ is {\it strongly
$C^r$ conjugate} to $S$ if it is $C^r$ conjugate to $S$ and $DS(0) =
DT(0)$. 

We call the  local $C^r$ diffeomorphism $T$ {\it
hyperbolic} if $0$ is a hyperbolic fixed point of $T$. 

Let $(T,U)$ be a $C^r$ local hyperbolic diffeomorphism, and  let $L=DT(0)$ be its derivative at $0$.
We assume that $L$ is of saddle type with associated  $L-$invariant splitting   $E = E^u \oplus E^s$ so that
$L \mid E^u$ is expanding and $L \mid E^s$ is contracting.

We will say that $T$ is {\it flat } in $U$ if 

\beq \label{flat-def}
T(E^u \bigcap U) \subset E^u \mbox{ and } T(E^s \bigcap U) \subset
E^s. 
\eeq

We will say that $T$ is {\it locally flat} at (or near) $0$ if there is
some open neighborhood $U$ of $0$ such that $T$ is flat in $U$. 

It is well-known (see \cite{Llave-Wayne}, \cite{Hirsch-Pugh}) that,
for $\delta > 0$ small, the sets $W^s_{\delta}(0,T)$ and $W^u_{\delta}(0,T)$ are $C^r$
embedded submanifolds of $E$ which are tangent at $0$ to $E^s$ and
$E^u$, respectively. 

The following lemma is well-known. 

\begin{Lem} \label{flat-lemma}
Every  hyperbolic local $C^r$ diffeomorphism $(T,U) $ is strongly $C^r$ conjugate to 
a locally flat one. 

More  precisely, one can find a ball $B_{\delta}(0) \subset U \bigcap
T(U)$ and a $C^r$ diffeomorphism $R$ from $B_{\delta}(0)$ onto its image such that $R(0) = 0,
DR(0) = I$, and the map $T_1 = R T R^{-1}$ is flat on
$R(T^{-1}B_{\delta}(0))$. 
\end{Lem}

We call an $R$ as in Lemma \ref{flat-lemma}, a {\it flattening }
map (or diffeomorphism) for $T$. 

{\bf Proof.}

We begin by identifying  $E$ with
$E^u \times E^s$ and take the product
representation 

\beq \label{T-def}
T(x,y) = (Ax + X(x,y), By + Y(x,y))
\eeq 

where $A \in Aut(E^u)$ is expanding, $B \in Aut(E^s)$ is contracting, and $X:E^u \times E^s \rarrow
E^u$ and $Y:E^u \times E^s \rarrow E^s$ are $C^r$ maps which vanish at
$(0,0)$ and have their partial derivatives also vanishing at $(0,0)$. 

Letting $\pi^u:E^u \times E^s \rarrow E^u$ and $\pi^s:E^u \times E^s
\rarrow E^s$ be the natural projections, and writing $E^u_{\delta} = \pi^u(B_{\delta}(0,0))$ and $E^s_{\delta} =
\pi^s(B_{\delta}(0,0))$, it is proved in \cite{Llave-Wayne} and \cite{Hirsch-Pugh}) that,
for small $\delta$, there are $C^r$ functions $g_u:E^u_{\delta} \rarrow E^s$ and
$g_s: E^s_{\delta} \rarrow E^u$ such that 

\beq \label{gu-conds}
g_u(0) = (0), \ Dg_u(0) = 0,
\eeq 

\beq \label{gs-conds}
g_s(0) = (0), \ Dg_s(0) = 0,
\eeq 

and $W^u_{\delta}(0,T)$ and $W^s_{\delta}(0,T)$ are the graphs of $g_u$
and $g_s$, respectively.

Setting $R(x,y) = (x - g_s(y), y- g_u(x))$, it is readily
verified that, for $\delta > 0$ small enough, the map $R$ is $C^r$
diffeomorphism from $B_{\delta}$ onto its image, and, hence,  is a local flattening map for $T$. 

\begin{Rem} \label{Ws-par-dep}
In section \ref{par-dep}, we consider continuous dependence of our
linearizations on parameters. It will be necessary to have the maps
$g_s, g_u$ depend continuously on the parameters as well. The most elegant proof of this
result is given in \cite{Llave-Wayne} where the Irwin method
for proving the existence of invariant manifolds is generalized and
simplified. 
\end{Rem} 

Our discussion so far shows that, in moving toward the proof of
Theorem \ref{Main_Theorem_1}, we may assume the $T$ is locally flat
after a $C^{1,\alpha}$ coordinate change.  

It will be convenient to have a stronger condition which requires another definition.  

\begin{Def} \label{h-linear-def}
Let $U$ be an open neighborhood of $0$ in the Banach space $E$, and
let $T:U \rarrow V$ be a $C^{1,\alpha}$ diffeomorphism from $U$ onto
its image with a hyperbolic fixed point at $0$. 
Let $L = DT(0)$, and let $E = E^u \oplus E^s$ be the hyperbolic
splitting associated to $L$. 

We say that $T$ is {\bf hyperbolically linear } (or {\bf h-linear }) in $U$ if 

\beq \label{T-Eu-invar}
T(x) = Lx \mbox{ for } x \in \left[E^u \bigcup E^s\right] \bigcap U. 
\eeq

\end{Def}

Observe that if $T$ is $h-$linear on $U$ and $U$ is small enough, then
it is locally flat in $U$,  and $T^{-1}$ is h-linear on $T(U)$. 
  Also, for positive integers $m$, $T^m$ is h-linear on
$\bigcap_{j=0}^m T^{-j}U$, and $T^{-m}$ is h-linear on $\bigcap_{j=0}^m T^{j}U$, 

In general, a $C^r$ local diffeomorphism with a hyperbolic fixed point
at $0$ will not be conjugate to even a $C^1$ $h-$linear one.

However, the following lemma shows that a $C^{1,\alpha}$ local diffeomorphism with an $\alpha-$hyperbolic fixed
point at $0$ is, in fact,  strongly $C^{1,\alpha}$ conjugate to an h-linear one. 

As in section \ref{Contracting_Case}, given an open neighborhood $U$
of $0$, consider the space $C^{1,\alpha}_U$ of
$C^{1,\alpha}$ maps $f:U \rarrow E$ such that $f(0) = 0$ and 
$Df(0) = 0$ with finite  D-H{\"{o}}lder norm $\n{f}_{U} = Hol(Df,U) < \infty$.

\begin{Lem} \label{local-h-lin} 
 For $0 < \alpha < 1$, let $E$ be a real Banach space, and let $L \in Aut(E)$ be
an $\alpha-$hyperbolic linear automorphism with hyperbolic splitting $E
= E^u \oplus E^s$.  Let $U$ be an open
neighborhood of $0$ in $E$, let $f \in C^{1,\alpha}_U$, and let $T =
L + f$. 

Then, for any $\gre > 0$, there are neighborhoods $V$ and $W$ of $0$
in $E$ and a $C^{1,\alpha}$ diffeomorphism $R$ from
$V$ onto $W$ such that 

\beq \label{V-W-gre}
V \bigcup W \subset B_{\gre}(0) \subset U \bigcap T(U) \bigcap T^{-1}(U),
\eeq

\beq \label{RDR-1}
R(0) = 0, \  DR(0) = I,
\eeq

and the following conditions are satisfied. 

Letting $T_1 = RTR^{-1}$ and $f_1 = T_1 - L$, and defining 

\[ V_1 = R \left[T(V) \bigcap V \bigcap T^{-1}(V) \right], \]

we have 

\beq \label{T1-lin}
T_1 \mbox{ is defined and h-linear on } V_1
\eeq

\beq \label{Lip-T1-L}    
 Lip(f_1, V_1 ) < \gre, 
\eeq

and 

 \beq  \label{Holder-Df1}
   Hol(Df_1,V_1 )<  \infty.
 \eeq

Further, the map $T_1$ is a $C^{1,\alpha}$ diffeomorphism from 
$V_1$ onto its image $T(V_1) \eqdef W_1$, and, setting $g_1 = T_1^{-1}
- L^{-1}$, we have 

\beq \label{T1inv-lin}
T_1^{-1}  \mbox{ is defined and h-linear on } W_1
\eeq

\beq \label{Lip-g1}
 Lip(g_1, W_1) < \gre, 
\eeq

and 

 \beq  \label{Holder-Dg1}
   Hol(Dg_1, W_1) <  \infty.
 \eeq

\end{Lem}

\begin{proof}

Taking $\gre > 0$ sufficiently small,  we may assume that $T$ is locally flat in the
ball $B_{\gre} = B_{\gre}(0)$.   

We use the product representation for $T$ and write

\[ T(x,y) = (A x + X(x,y), By + Y(x,y)). \]

with $A$ $\alpha-$expanding, $B$ $\alpha-$contracting, and $X, Y$
$C^{1,\alpha}$ maps vanishing at $(0,0)$ together with their partial
derivatives. In these coordinates, of course,  

\[ L(x,y) = (Ax, By). \]

It  suffices to find a small neighborhoods $V, W$ of $(0,0)$ 
contained in $B_{\gre}(0,0)$ and a local $C^{1,\alpha}$ map $R$ from $V$
onto $W$  such that

\beq \label{RDR-2}
R(0,0) = (0,0), \  DR(0,0) = I,
\eeq

and,  setting $T_1 = RTR^{-1}$, we have

\beq \label{T1-h-lin} 
T_1(x,0) = (Ax, 0) \mbox{ and } T_1(0,y) = (0,By)
\eeq

for $(x,y) \in W$. 

Indeed, it is clear that, for $V$ and $W$ small enough, the conditions
(\ref{Lip-T1-L}), (\ref{Holder-Df1}),  (\ref{Lip-g1}), and
(\ref{Holder-Dg1}) will all be satisfied. 

Let us proceed to find the small neighborhoods $V, W$ and the appropriate map
$R$. 

Let $E^u_{\gre} = \pi^u(B_{\gre})$ and $E^s_{\gre} =
\pi^s(B_{\gre})$, and consider the maps 
$T^u:E^u_{\gre}  \rarrow E^u$ and $T^s:E^s_{\gre}  \rarrow E^s$
defined by 

\[ T^u(x) = Ax + X(x,0), \ T^s(y) = By + Y(0,y). \] 

Clearly,  $T^u$  has $0$ as an $\alpha-$expanding
fixed point, and $T^s$  has $0$  as an $\alpha-$contracting fixed point. 

By Theorem \ref{Lin_Contraction_2}, and the remark following it, there
are small neighborhoods $V^u$ of $0$ in $E^u$, $V^s$ of $0$ in $E^s$,
and $C^{1,\alpha}$ diffeomorphisms $R_1:V^u \rarrow R_1(V^u)$ and
$R_2:V^s \rarrow R_2(V^s)$ such that 

\[ R_1(0) = 0, \ DR_1(0) = I, R_2(0) = 0, \ DR_2(0) = I,  \]

\[ R_1 T^u R_1^{-1} = A \mbox{ on } R_1(V^u),  \]

and 
\[ R_2 T^s R_2^{-1} = B \mbox{ on  } R_2(V^s). \] 

Now, let $V = V^u \times V^s, \ W = R_1(V^u) \times R_2(V^s)$, and let
$R:V \rarrow R(V)$ be the product map

\[ R(x,y) = (R_1(x), R_2(y)). \] 

The map $R$ satisfies (\ref{RDR-2}), and the map $T_1 = RTR^{-1}$
satisfies

\[ T_1(x,0)  = R^{-1}TR(x,0) = (Ax, 0) \]

and 

\[ T_1(0,y)  = R^{-1}TR(0,y) = (0, By). \]

This completes the proof of Lemma \ref{local-h-lin}. 

\end{proof}

We will need the analog of Lemma \ref{local-h-lin}  for powers of $T$.

\begin{Lem} \label{local-h-lin-power}
Let $\gre > 0$, $n$ be a positive integer, and let $R, T_1, f_1, g_1,
V_1, W_1$ be as in the statement of Lemma \ref{local-h-lin}. 

Then, there are neighborhoods $V_n \subset V_1$ and $W_n \subset W_1$
such that, $T_1^n$ is an h-linear  diffeomorphism from $V_n$ onto $W_n$, and,
defining $f_n = T_1^n - L^n$ and $g_n = T_1^{-n} - L^{-n}$, we have 

\beq \label{Lip-fn-1}    
 Lip(f_n, V_n ) < \gre, 
\eeq

 \beq  \label{Holder-Dfn-1}
   Hol(Df_n,V_n )<  \infty,
 \eeq

\beq \label{Lip-gn-1}
 Lip(g_n, W_n) < \gre, 
\eeq

and 

 \beq  \label{Holder-Dgn-1}
   Hol(Dg_n, W_n) <  \infty.
 \eeq

\end{Lem}

\begin{proof}

As usual, writing $B_{\delta} = B_{\delta}(0)$ for any positive $\delta$,
first pick 
$\delta \in (0,\gre)$ so that 

\[ B_{\delta} \subset \bigcap_{k=-n}^n T_1^k(V_1 \bigcap W_1). \]

Next, choose $V_n = B_{\delta_n}(0)$ where $\delta_n$ is small enough so
that $T_1^j(V_n) \subset B_{\delta}$
for $\n{j} \leq n$. 

It is then clear that, setting $W_n =
T_1^n(V_n)$, we have  $V_n \subset V_1, \ W_n \subset W_1$, that $T_1^n$
is an h-linear diffeomorphism taking $V_n$ onto $W_n$. 

The only thing remaining is to get the Lipschitz estimates (\ref{Lip-fn-1}) and
(\ref{Lip-gn-1}).

But, an easy induction shows that $DT^n(0) = L^n$, which implies that $Df_n(0) = 0 =
Dg_n(0)$.

Then, it is clear that we can shrink $\delta_n$ enough so that (\ref{Lip-fn-1}) and
(\ref{Lip-gn-1}) hold.  In fact, using (\ref{Lip-f-1}), (\ref{Lip-fm}), and
(\ref{Lip-gm}), one can get explicit estimates of $\delta_n$ in terms of $Lip(f_1)$. 

\end{proof}

\subsection{Extension of local h-linear diffeomorphisms to global diffeomorphisms}

Let $T_1$ be the diffeomorphism obtained in Lemma \ref{local-h-lin}. 

We now wish to use Lemma \ref{Hol-Dg} to find an h-linear  $C^{1,\alpha}$
diffeomorphism $S$ defined on all of $E$  which agrees with $T_1$ on a
small neighborhood of $0$ and is linear off a slightly larger ball
$B_{\delta}$ where $\delta$ is small relative to $Hol(DS,E)$. 

Let $C^{1,\alpha}_0$ denote the space of $C^{1,\alpha}$ maps from $E$ into
$E$ such that $f(0) = 0$ and $Df(0) = 0$ with the norm 

\beq \label{f-global-C1-gra}
Hol(Df) = \sup_{x \neq y} \frac{\n{Df_x - Df_y}}{\n{x-y}^{\alpha}} < \infty.
\eeq

Recall that $Aut(E)$ is the set of linear automorphisms of $E$.

Let $D^{1,\alpha}_0$ denote the set of maps $T = L + f$ with $L \in
Aut(E)$ and $f \in C^{1,\alpha}_0$ such that 

\beq \label{Lip-f-E}
Lip(f,E) \leq \frac{1}{2\n{L^{-1}}}.
\eeq

From the Inverse Function Theorem, one sees that the maps 
$ T \in D^{1,\alpha}_0$ are bijective $C^{1,\alpha}$ maps with $C^{1,\alpha}$
inverses such that $T(0) = 0$ and $DT(0) = L$. 

By a $C^{1,\alpha}$ diffeomorphism of $E$, we will mean an element of
$D^{1,\alpha}_0$. More general diffeomorphisms can, of course, be
defined, but we don't need them in this paper. 

\begin{Def}
Given a Banach space $E$, let $T = L + f$ where $L$ is a bounded
linear map,  and let $f \in C^{1,\alpha}_0(E)$. 

We define the {\bf non-linear support of } $T$ to be the set of points
$x \in E$ such that $T(x) \neq L(x)$.  We denote this by 

\[ nsupp(T). \]

The {\it non-linear size} of $T$, denoted $N_s(T)$, is defined by

\beq \label{Ns-T}
N_s(T) = \inf\{\xi \in \R: nsupp(T) \subset B_{\xi}(0) \}.
\eeq

\end{Def}

It is clear that $N_s(T) = 0$ if and only if $T=L$ and that, in
general,  $N_s(T)$ can range from $0$ to $\infty$.

We will be interested in h-linear $C^{1,\alpha}$ maps $T$ for which $N_s(T)$ is
small compared to $Hol(DT,E)$.

Now, let $C_1(\grl) < C_2(\grl)$ be the constants given in  Lemma
\ref{Hol-Dg}, and 
consider the map $T_1 = L + f_1$ defined on the neighborhood $W_1$
obtained in Lemma \ref{local-h-lin}. 

For positive $\delta > 0$, let 

\beq \label{f2-T2-def}
f_2 = \grl_{\delta} f_1, \ \  T_2 = L + f_2
\eeq

where $\grl_{\delta}$ is the $\delta$-scaled version of the unit bump
function $\grl$ defined in (\ref{delta-scaled}). 

Let

\beq \label{gre-L-est}
0 < \gre < \frac{m(L)}{2} = \frac{1}{2\n{L^{-1}}}. 
\eeq

Since $Df_1(0) = 0$, Lemma (\ref{Hol-Dg}) and the Mean Value
Theorem say that, for any $\gre > 0$,  we may choose $\delta >0$ small enough so that
$B_{\delta}(0) \subset W_1$ and

\beq \label{Lip-f2-grd}
Lip(f_2,B_{\delta}(0)) \leq \delta^{\alpha} Hol(Df_2,B_{\delta}(0))
<\delta^{\alpha} C_2(\grl) Hol(Df_1,B_{\delta}(0)) < \gre.
\eeq

From (\ref{gre-L-est}) and (\ref{Lip-f2-grd}) it is easily seen that $T_2$ is a $C^{1,\alpha}$ global
diffeomorphism on $E$. 

Let $B_{\delta}^c = E \setminus B_{\delta}$ denote the set-theoretic
complement of $B_{\delta}$ in $E$. 

Since $f_1$ vanishes in 

\[ B_{\delta}(0,0) \bigcap \left[(E^u \times \{0\}) \bigcup (\{0\} \times E^s) \right],
\]

and $\grl_{\delta}$ vanishes in $B_{\delta}^c$, we have that 

\[ f_2 \mbox{ vanishes on } (E^u \times \{0\}) \bigcup (\{0\} \times
E^s) \bigcup B_{\delta}^c. 
\]

Thus, $T_2$ is a $C^{1,\alpha}$ h-linear
global diffeomorphism from $E$ onto $E$ with non-linear support in $B_{\delta}$.

Next, let $n$ be an integer greater than 1, and consider the powers
$T_2^n$ and $T_2^{-n}$. 

Let 

\[ B_n^+ = \bigcap_{k=0}^nT_2^{-k}B_{\delta}^c, \] 

and 

\[ B_n^- = \bigcap_{k=0}^nT_2^kB_{\delta}^c, \] 

In view of (\ref{recursive-fn}), we see that for each $1 \leq j \leq
n-1$, and $x \in B_n^+$, if $f_j(x) = 0$, then $f_{j+1}(x) = 0$ as
well.

In particular, 

\beq \label{fm-0}
f_n(x) = 0 \mbox{ for } x \in B_n^+.
\eeq

A similar argument shows that 

\beq \label{gm-0}
g_n(x) = 0 \mbox{ for } x \in B_n^-,
\eeq

which yields 

\beq \label{nsupp-Tm}
nsupp(T_2^n) \subset (B_n^+)^c = \bigcup_{k=0}^nT_2^{-k}B_{\delta},
\eeq

and 

\beq \label{nsupp-Tm-inv}
nsupp(T_2^{-n}) \subset (B_n^-)^c = \bigcup_{k=0}^nT_2^kB_{\delta}. 
\eeq

Now, defining $\delta_n$ by

\beq \label{grd-m-def}
\delta_n = \delta \ \max(Lip(T_2)^n, Lip(T_2^{-1})^n),
\eeq

we see that (\ref{nsupp-Tm}) and (\ref{nsupp-Tm-inv}) give 

\beq \label{nsupp-Tm-Tm-inv}
nsupp(T_2^n) \bigcup nsupp(T_2^{-n})   \subset  B_{\delta_n}.
\eeq

   Also, obviously, for a fixed positive integer $n$, we can choose
   the bump function scaling constant $\delta$ small enough so that
   the non-linear sizes  of $T_2^n$ and $T_2^{-n}$ are arbitrarily
   small.

With the above  definitions, let us summarize
the results obtained above in the following proposition. For notational convenience, we set $S = T_2$. 

\begin{Prop} \label{global-lin-prop}
Let $0 < \alpha < 1$, let $E$ be a $C^{1,\alpha}$ Banach space, and let $L$
 be an $\alpha-$hyperbolic linear automorphism of saddle type. 
Let $U$ be an open neighborhood of $0$,  let $f:U \rarrow E$ be a
 function in $C^{1,\alpha}_U$, and let $T = L + f$. 

Then, for any $\gre > 0$ and any non-zero integer $n$, there are a
$\delta_n = \delta_n(L,f,\gre) < \gre$ and a  global h-linear $C^{1,\alpha}$ diffeomorphism
$S \in D^{1,\alpha}_0$  which is strongly $C^{1,\alpha}$ 
conjugate to $T$ on $B_{\delta_n}$, and setting $f_n = S^n - L^n$, we have

\beq \label{Lip-fn-E}
Lip(f_n,E) < \gre
\eeq

and

\beq \label{nsupp-S}
nsupp(S^n) \subset B_{\delta_n}. 
\eeq

\end{Prop}

We will need the following simple consequence of Proposition \ref{global-lin-prop}. 

\begin{Lem} \label{Basic-X-Y-Lemma}
Let $S, n, \delta_n, f_n$ be as in Proposition \ref{global-lin-prop},
and,  using the product representation 
$E = E^u \times E^s$, write $f_n(x,y) = (X_n(x,y), Y_n(x,y))$. 

Then, for every $(x,y) \in E^u \times E^s$  we have

\beq \label{max_X_x_Y_x_gra_n}
 \max(\n{X_{n,x}(x,y)}, \n{Y_{n,x}(x,y)}) \leq Hol(Df_n) \min(\delta_n^{\alpha},\n{y}^{\alpha})
\eeq

and 

\beq \label{max_X_y_Y_y_gra_n}
 \max(\n{X_{n,y}(x,y)}, \n{Y_{n,y}(x,y)}) \leq Hol(Df_n)  \min(\delta_n^{\alpha},\n{x}^{\alpha}).
\eeq

\end{Lem} 

\begin{proof}

  We will only give the arguments for $X_n$ since those needed for
  $Y_n$ are similar. 

Because we are using the maximum norm on $E^u \times E^s$, we have 

\[ Hol(Df_n) = \max(Hol(DX_n), Hol(DY_n)). \]

Since $X_n(x,0) = 0$ for each $(x,0)$ we also have that $X_{n,x}(x,0)
= 0$ for each $(x,0)$.  Similarly, $X_n(0,y) = 0$ implies that $X_{n,y}(0,y) =
0$ for each $(y,0)$.  

Hence, 

\beq \label{X_x_0}
\n{X_{n,x}(x,y)} = \n{X_{n,x}(x,y) - X_{n,x}(0,0)} \leq Hol(DX_n)  \n{(x,y)}^{\alpha}, 
\eeq

\beq \label{X_y_0}
\n{X_{n,y}(x,y)} = \n{X_{n,y}(x,y) - X_{n,y}(0,0)} \leq Hol(DX_n)  \n{(x,y)}^{\alpha}, 
\eeq

\beq \label{X_x_est-1}
\n{X_{n,x}(x,y)} = \n{X_{n,x}(x,y) - X_{n,x}(x,0)} \leq Hol(DX_n)  \n{y}^{\alpha} 
\eeq

and 

\beq \label{X_y_est-1}
\n{X_{n,y}(x,y)} = \n{X_{n,y}(x,y) - X_{n,y}(0,y)} \leq Hol(DX_n)  \n{x}^{\alpha} 
\eeq

The non-linear support condition (\ref{nsupp-S}) for $X_n$ implies that $X_n(x,y)$ and its
partial derivatives vanish at points $(x,y)$ such that $\n{(x,y)} >
\delta_n$. So, (\ref{X_x_0}) and (\ref{X_y_0}) imply 

\beq
\max(\n{X_{n,x}(x,y)}, \n{X_{n,y}(x,y)}) \leq Hol(DX_n) \delta_n^{\alpha} 
\eeq

This, together with  (\ref{X_x_est-1}), (\ref{X_y_est-1}) and the fact that $Hol(DX_n) \leq Hol(Df_n)$,   gives
formulas (\ref{max_X_x_Y_x_gra_n}) and (\ref{max_X_y_Y_y_gra_n}). 

\end{proof}

\subsection{Norm Conditions induced by $\alpha$-hyperbolicity} \label{norm-gra-conds}

Assuming $L$ is $\gra-$hyperbolic, and using  the product representation $E = E^u \times E^s$, the map $L$ has
the form

\[ L(x,y) = (Ax, By). \]

with 

\beq \label{sp-1}
\max(\rho(A^{-1}), \ \rho(B)) < 1.
\eeq

By Lemma \ref{re_norm}, we may assume the maps $A$ and $B$ satisfy

\beq \label{norm-Ainv-B}
\n{A^{-1}} < 1 \mbox{  and  } \n{B} < 1.
\eeq

Next, using 
(\ref{gra-hyperbolic-def}) and formula (\ref{spectral-radius-formula}), we choose a positive integer $m$
such that 

\[ \left(\n{A^{-m}}\n{A^m}\n{B^m}^{\alpha} \right)^{\frac{1}{m}} < 1, \]

and 

\[ \left( \n{B^{-m}}\n{B^m}\n{A^{-m}}^{\alpha}\right)^{\frac{1}{m}} <
1. \]

Taking the $m$-th power, these imply  that 

\beq \label{m-s-cond}
 \n{A^{-m}}\n{A^m}\n{B^m}^{\alpha} < 1 
\eeq

and

\beq \label{m-u-cond}
\n{B^{-m}}\n{B^m}\n{A^{-m}}^{\alpha} < 1. 
\eeq

Following this, we choose $\gre > 0$ small enough so that

\beq \label{Am-gre-est}    
\min(m(A^m) - \gre, \ m(B^{-m}) - \gre) > 1,
\eeq

\beq \label{B+gre-est}
\max(\n{A^{-1}} + \gre, \n{B} + \gre) < 1, 
\eeq

\beq \label{norm-A-B-m-1}
\n{A^{-m}}(\n{A^m}+\gre)(\n{B^m} + \gre)^{\alpha} < 1, 
\eeq

and

\beq \label{norm-A-B-m-2}
\n{B^m}(\n{B^{-m}} + \gre)(\n{A^{-m}} + \gre)^{\alpha} < 1.
\eeq

The estimates (\ref{Am-gre-est}) and (\ref{B+gre-est}) of course imply that 

\beq \label{A+gre-est}    
m(A) - \gre > 1
\eeq

and 

\beq \label{Bm+gre-est}
\n{B^m} + \gre < 1. 
\eeq

Finally, we choose $0 < \eta < 1$ close enough to 1 so that

\beq \label{norm-A-B-m-1-eta}
\n{A^{-m}}(\n{A^m}+\gre)(\n{B^m} + \gre)^{\alpha
  \eta}(\n{A^m}+\gre)^{\alpha(1 - \eta)}   < 1, 
\eeq

and 

\beq \label{norm-A-B-m-2-eta}
\n{B^m}(\n{B^{-m}} + \gre)(\n{A^{-m}} + \gre)^{\alpha \eta}
(\n{B^{-m}} + \gre)^{\alpha(1 - \eta)}< 1.
\eeq

\section{Proof of Theorem \ref{Main_Theorem_1}} \label{Main_Proof} 

Let $E, T, L$ be as in the hypotheses of Theorem
\ref{Main_Theorem_1}, 
and let $E = E^u \oplus E^s$ be the hyperbolic
splitting associated to $L$.  

Applying Proposition \ref{global-lin-prop} for $n=1$ and arbitrary $\gre > 0$, 
there are a $\delta = \delta_1 \in (0,\gre)$ and  an h-linear map $S \in D^{1,\alpha}_0$ of the form $S = L + f_1$ with $f_1 \in
C^{1,\alpha}_0$ such that $S$ is strongly $C^{1,\alpha}$ conjugate to $T$
on $B_{\delta}$ with $nsupp(S) \subset B_{\delta}$. 

In subsection \ref{Global-Lin-Th} below, we will show that for $\gre$ satisfying conditions (\ref{Am-gre-est})-(\ref{norm-A-B-m-2})  and $\delta > 0$ sufficiently small,
 the map $S$  is globally strongly $C^1$ conjugate to $L$. That
is, there is a $C^1$ diffeomorphism $R$ mapping $E$ onto itself  such that $R(0) =
0, \ DR_0 = I$, and $RSR^{-1} = L$.  Then, in 
subsection \ref{D-Holder-Cont-of-R},  we will show that, adding conditions (\ref{norm-A-B-m-1-eta}) and (\ref{norm-A-B-m-2-eta}) and shrinking $\grd$ further,  the map $R$ is
$C^{1,\beta}$ on bounded subsets 
of the orbit of $B_{\delta}(0)$. 

Once these things are done, the restriction
to a small neighborhood of $0$ of the composition
of $R$ with the local $C^{1,\alpha}$ conjugacy from $S$ to $T$ provides
a local $C^{1,\beta}$ conjugacy from $T$ to $L$ and completes the
proof of Theorem \ref{Main_Theorem_1}.

\subsection{Global $C^1$ linearization of the map $S$.} \label{Global-Lin-Th}

To $C^1$ linearize $S = L+f_1$ on $E$, we will employ the method described in
section \ref{Associated-Ops}.   The map 
$R$ will have the form $R = I + \phi$, where 

\beq 
\phi = (I-H)^{-1}L^{-1}f_1
\eeq

where $H$ is an automorphism on a suitable Banach space of
functions $\sE = \sE^s \oplus \sE^u$ which we now define.

We identify $E$ with the space $E^u \times E^s$ as above.

Given  a Banach space $Z$ and a $C^1$ function
$\psi: E^u \times  E^s  \rarrow Z$, consider the following two real-valued functions.

\beq
\grg_1(\psi,x,y, \alpha) = \frac{\n{\psi_x(x,y)}}{\n{y}^{\alpha}} 
\eeq

\beq
\grg_2(\psi,x,y, \alpha) =  \frac{\n{\psi_y(x,y)}}{\n{x}^{\alpha}}
\eeq

Each of the preceding functions is defined to have value zero if its denominator vanishes. 
Otherwise, it is given by the indicated expression. 

Consider the corresponding suprema:

\[ \n{\grg_i(\psi,\alpha)} = \sup_{x \neq 0, y \neq 0} \grg_i(\psi,x,y,\alpha) \mbox{
for } i = 1,2, \]

and let $C^1_0(E,Z)$ be the set of $C^1$ functions from $\psi:E^u
\times E^s$ to $Z$ such that 

\beq \label{psi_cond_0}
\psi(x,0) = \psi(0,y) = 0 \ \forall \ (x,y) ,
\eeq

\beq \label{psi_xy_cond_0}
\psi_x(0,0) = \psi_y(0,0) = 0,
\eeq

and 

\beq \label{grg_i_finite}
\n{\grg_i(\psi)} = \n{\grg_i(\psi, \alpha)} < \infty \mbox{ for } i=1,2.
\eeq

For such $\psi$, set 

\beq \label{dnorm_psi}
\dnorm{\psi} = \dnorm{\psi}_{\alpha} = \max_{i=1,2} \n{\grg_i(\psi, \alpha)}.
\eeq

One can check that  the set $C^1_0(E,Z)$ is a linear space with the
usual pointwise operations of addition and scalar multiplication, and
that the function $\dnorm{\psi}$ provides a norm on it so that the
pair $(C^1_0(E,Z),\dnorm{\cdot})$ becomes a Banach space.  

Setting $Z$ to be $E^u$ and $E^s$, in turn, gives the two Banach spaces 
$\sE^s = C^1_0(E,E^u)$ and $\sE^u = C^1_0(E,E^s)$.  For $\grf =
(\grf_s, \grf_u) \in \sE = \sE^s \times \sE^u$, we define the norm 

\[ \dnorm{\grf} = \dnorm{(\grf_s, \grf_u)} = \max(\dnorm{\grf_s}, \dnorm{\grf_u}). \]

For $\grf = (\grf_s, \grf_u) \in \sE^s \times  \sE^u$, the map 
$H(\grf)(x,y) = L^{-1}\grf(S(x,y))$ is expressed as 

\[ H(\grf)(x,y) = (H_s(\grf_s)(x,y), H_u(\grf_u)(x,y)) \]

where 

\beq \label{H1-def}
H_s(\grf_s)(x,y) = A^{-1}\grf_s(Ax + X_1(x,y), By + Y_1(x,y))
\eeq

and 

\beq \label{H2-def}
H_u(\grf_u)(x,y) = B^{-1}\grf_u(Ax + X_1(x,y), By + Y_1(x,y)).
\eeq

Since $S$ is h-linear on $E$ and the composition of $C^{1,\alpha}$ maps is
again $C^{1,\alpha}$, we see that, at least as an operator on the function space 
$\sE^s \times \sE^u$ without consideration of norms, $H$ is represented as the direct sum of
operators $H = H_s \oplus H_u$.  


\begin{Prop} \label{Main_Prop}
Let $E, T, L, f$ be as in the statement of Theorem \ref{Main_Theorem_1},
and assume that $L(x,y) = (Ax, By)$ and $m > 0, \gre$ are such that
(\ref{Am-gre-est})-(\ref{Bm+gre-est}) 
are valid. There is a
$\delta \in (0,\gre)$  such that if 
$S$ as in Proposition \ref{global-lin-prop}, then, the associated maps $H_s$ and $H_u$
satisfy 

\beq \label{Hs-aut-norm}     
H_s \in Aut(\sE^s) \mbox{ and } \n{H_s^m} < 1
\eeq

and 

\beq \label{Hus-aut-norm}
H_u \in Aut(\sE^u) \mbox{ and } \n{H_u^{-m}} < 1.
\eeq

Hence, $\rho(H_s) < 1$ and $\rho(H_u^{-1}) < 1$, so $H$ is
hyperbolic. 

\end{Prop}

Assuming Proposition \ref{Main_Prop},  we use Lemma \ref{re_norm} to
find norms $\n{\cdot}_u$ and $\n{\cdot}_s$ on $\sE^u$ and $\sE^s$,
respectively so that

\beq \label{norm-Hs-less-one}
\n{H_s} < 1
\eeq

and

\beq \label{norm-Hu-inv-less-one}
\n{H_u^{-1}} < 1.
\eeq

Letting $I_s:\sE^s \rarrow \sE^s, I_u:\sE^u
\rarrow \sE^u$ denote the identity maps on $\sE^s, \ \sE^u$,
respectively, we then have 

\[ (I_s-H_s)^{-1} = \sum_{n \geq 0}   H_s^n \]

and 

\bqns
(I_u - H_u)^{-1} & = & H_u^{-1}H_u(I_u - H_u)^{-1} \\
& = & H_u^{-1}(H_u^{-1})^{-1}(I_u - H_u)^{-1} \\
& = & H_u^{-1}((I_u - H_u)H_u^{-1})^{-1} \\
& = & H_u^{-1}(H_u^{-1}-I_u)^{-1} \\
& = & -H_u^{-1}(I_u- H_u^{-1}) \\
 & = & -H_u^{-1}\sum_{n \geq 0}H_u^{-n} \\
 & = & -\sum_{n \geq 1}H_u^{-n}
\eqns

giving 

\beq \label{I-H-inv}
(I-H)^{-1} = (I_s-H_s)^{-1} \oplus (I_u - H_u)^{-1} = (\sum_{n \geq 0}   H_s^n, -\sum_{n \geq 1}H_u^{-n}).
\eeq

Then, formula (\ref{alg-soln}) in Section
\ref{Associated-Ops} shows that a (global) $C^1$
linearization of $S$ is given by $R = I + \psi$
where

\beq 
\psi = (I-H)^{-1}L^{-1}f = (\psi_s, \ \psi_u)
\eeq

with 

\beq \label{psi_s_def}
\psi_s = (I_s - H_s)^{-1}A^{-1}X_1 = \sum_{n \geq 0} H_s^n (A^{-1}X_1)
\eeq

and 
 
\beq \label{psi_u_def}
\psi_u = (I_u - H_u)^{-1}B^{-1}Y_1 = -\sum_{n \geq 1} H_u^{-n}(B^{-1}Y_1).
\eeq

Let us proceed to the proof of Proposition \ref{Main_Prop}. 

We first show that the statement (\ref{Hs-aut-norm}) is implied by
assumptions 
(\ref{B+gre-est}), (\ref{norm-A-B-m-1}), 
and
(\ref{A+gre-est}). 

Following this and using the expression 

\beq \label{Hu-inv-m}
 H_u^{-m}(\grf_u) = B^m\grf_u(A^{-m}x + X_{-m}(x,y), B^{-m}y +
Y_{-m}(x,y)), 
\eeq

it is easy to see that 
the same method works for (\ref{Hus-aut-norm}) by
replacing $S, A, B$ by $S^{-1}, B^{-1}, A^{-1}$, respectively.
using the assumptions 
(\ref{B+gre-est}),
(\ref{norm-A-B-m-2}), and 
(\ref{A+gre-est}).

Thus, it suffices to prove (\ref{Hs-aut-norm}). 

In the sequel, all real numbers $\delta_n$ will be assumed to be in
$(0,\gre)$ where $\gre$ satisfies
(\ref{norm-Ainv-B})-(\ref{norm-A-B-m-2}). 
 Also, once a condition is specified by a choice of $\delta_n$, it will
continue to hold for smaller $\delta_n$, so this will be assumed without
further mention. 

Moving to (\ref{Hs-aut-norm}), for a positive integer $n$, 
consider the following real numbers

\[ K_{n,1} = \n{A^{-n}}(\n{A^n}+\gre)(\n{B^n}+\gre)^{\alpha},  \]

\[ K_{n,2} = \n{A^{-n}}(\n{A^n}+\gre)^{\alpha} \delta_n^{\alpha} Hol(Df_n), \]

\[ K_{n,3} = \n{A^{-n}}(\n{B^n}+\gre)^{\alpha} \delta_n^{\alpha} Hol(Df_n), \]

\[ K_{n,4} = \n{A^{-n}}(\n{A^n} + \gre) (\n{B^n} + \gre),  \] 

and 

\[ \grt_n = \max(K_{n,1} + K_{n,2}, K_{n,3} + K_{n,4}). \] 

By (\ref{B+gre-est}), we see that $\n{B^n} + \gre < 1$, so $K_{n,4} < K_{n,1}$.

Letting  $n=m$ with $m$ as in 
(\ref{norm-A-B-m-1}) and  (\ref{norm-A-B-m-2}),
we have 

\[  K_{m,4} < K_{m,1} < 1. \]

So, we may choose $\delta_m$ small enough so that $\grt_m < 1$. 

Further, for $n=1$ or $n=m$, and $\delta_n$ small enough, we can guarantee that

\beq \label{Lip-X-Y-gre}
\max(Lip(X_1), Lip(Y_1)) <  \gre, 
\eeq

\beq \label{Lip-X-Ym-gre}
\max(Lip(X_m), Lip(Y_m)) <  \gre, 
\eeq

and 

\beq \label{Ns-X-Y-m}
\max(N_s(X_m), N_s(Y_m)) <  \delta_m^{\alpha} Hol(Df_m).
\eeq

In addition, since we are using the $\delta$-scaled version of the bump
function $\grl$ for $S$, we have that the functions $X_1, Y_1, X_m,
Y_m$ all vanish on 

\beq \label{X1-Y1-Xm-Ym}
\left(E^u \times \{0\} \bigcup \{0\} \times E^s \right) \bigcup B_{\delta}^c.
\eeq

Now, the main step in the proof of (\ref{Hs-aut-norm}) is

\begin{Prop} \label{pf-Hs-contraction}
Assume that conditions (\ref{Am-gre-est})-(\ref{norm-A-B-m-2})  on
$\gre$ are satisfied. 

Then, for any $\psi \in \sE^s$, and $n=1$ or $n=m$, 
we have 

\beq \label{grg1n-Hs-norm-est}
\grg_1(H_s^n(\psi)) \leq (K_{n,1} + K_{n,2}) \dnorm{\psi}
\eeq

and 

\beq \label{grg2n-Hu-norm-est}
\grg_2(H_s^n(\psi)) \leq (K_{n,3} + K_{n,4}) \dnorm{\psi}. 
\eeq

\end{Prop}

Since 

\[ \dnorm{\psi} =  \max(\grg_1(\psi), \grg_2(\psi)), \]

conditions (\ref{grg1n-Hs-norm-est}) and (\ref{grg2n-Hu-norm-est}) imply that, 
for $n=1$ or $n=m$, we have 

\[ \dnorm{H_s^n(\psi)} < \grt_n \dnorm{\psi}. \]

Thus, for $n=1$ and $n=m$, the maps $H_s^n$ are bounded linear maps and $\n{H_s^m} < 1$. 
They are also bijective with inverses 

\[ H_s^{-n}(\psi) = A^n \psi(S^{-n}). \]

By the open mapping theorem, the maps $H_s^n$ are automorphisms (i.e.,
have bounded inverses) and this proves (\ref{Hs-aut-norm}).

Let us proceed to the proof of Propostion \ref{pf-Hs-contraction};
i.e., to the proofs of (\ref{grg1n-Hs-norm-est}) and (\ref{grg2n-Hu-norm-est}).

Define the functions $u_n(x,y), v_n(x,y)$ by

\beq
u_n =  u_n(x,y) = A^nx + X_n(x,y),
\eeq

and 

\beq
v_n =  v_n(x,y) = B^ny + Y_n(x,y).
\eeq

\begin{Lem} \label{u-v-est}
For every $(x,y) \in E^u \times E^s$, we have

\beq \label{u-est}
\n{u_n(x,y)} \leq (\n{A^n} + \gre) \n{x}, \ \
\eeq

\beq \label{v-est}
 \n{v_n(x,y)} \leq (\n{B^n} + \gre) \n{y}, \ \
\eeq

\beq \label{ux-est}
\n{u_{n,x}} \leq \n{A^n} + \gre
\eeq

\beq \label{uy-est}
\n{u_{n,y}} \leq Hol(Df_n) \min(\delta_n^{\alpha}, \n{x}^{\alpha}), \ \
\eeq

\beq \label{vx-est}
 \n{v_{n,x}} \ \leq Hol(Df_n) \min(\delta_n^{\alpha}, \n{y}^{\alpha}), \ \ 
\eeq

and 

\beq \label{vy-est}
\n{v_{n,y}} \leq \n{B^n} + \gre. 
\eeq
\end{Lem}

\begin{proof}

The inequalities (\ref{ux-est}) and (\ref{vy-est}) follow from (\ref{Lip-X-Ym-gre}),
and  inequalities (\ref{uy-est}) and (\ref{vx-est}) are a restatement
of (\ref{max_X_x_Y_x_gra_n}) and (\ref{max_X_y_Y_y_gra_n}).

From (\ref{Lip-X-Ym-gre}), the Mean Value Theorem, and the vanishing
of $X_n, Y_n$ on

\[ E^u \times \{0\} \bigcup \{0\} \times E^s, \]

we get 

\bqns
 \n{u_n(x,y)} & = & A^nx + X_n(x,y) \\
 & \leq & \n{A^n}\n{x} + \n{X_n(x,y) - X_n(0,y)} \\
 & \leq & \n{A^n}\n{x} + \left( \sup_{0 < s < 1} \n{X_{n,x}(s x, y)}
\right) \n{x} \\
 & \leq & \n{A^n}\n{x} + \gre \n{x}, 
\eqns

and 

\bqns
 \n{v_n(x,y)} & = & B^ny + Y_n(x,y) \\
 & \leq & \n{B^n}\n{y} + \n{Y_n(x,y) - Y_n(x,0)} \\
 & \leq & \n{B^n}\n{y} + \left(\sup_{0 < s < 1} \n{Y_{n,y}(x, s y)} \right) \n{y} \\
 & \leq & \n{B^n}\n{y} + \gre \n{y}. 
\eqns

\end{proof}

\noindent {\bf Proof of (\ref{grg1n-Hs-norm-est}).}

Since $\n{v_{n,x}} = \n{Y_{n,x}(x,y)} = 0$ if $\n{(x,y)} > \delta_n$, we have that
$\n{v_{n,x}} \neq 0$ implies that $\n{x}^{\alpha} < \delta_n^{\alpha}$. 

With $K_{n,1}$ and $K_{n,2}$ defined above,  we have 

\bqns 
\n{A^{-n}} \n{v_n}^{\alpha}\n{u_{n,x}}  & \leq & \n{A^{-n}} (\n{B^n}+\gre)^{\alpha}\n{y}^{\alpha} \n{u_{n,x}} \\
   & \leq & \n{A^{-n}} (\n{B^n}+\gre)^{\alpha}\n{y}^{\alpha}(\n{A^n}+\gre) \\
 & \leq & K_{n,1} \n{y}^{\alpha}
\eqns

and 

\bqns
\n{A^{-n}} \n{u_n}^{\alpha} \n{v_{n,x}} & \leq  &  \n{A^{-n}} (\n{A^n} +
\gre)^{\alpha} \n{x}^{\alpha} \n{v_{n,x}} \\
 & \leq &  \n{A^{-n}} (\n{A^n} + \gre)^{\alpha}\delta_n^{\alpha} Hol(Df_n)
\n{y}^{\alpha} \\
 & \leq & K_{n,2} \n{y}^{\alpha}.
\eqns

For notational convenience, let us write $u = u_n, v = v_n$.

Then, we have 

\bqns 
\n{\partial_xH^n_s(\psi)(x,y)} & = & \n{ A^{-n} \partial_x \psi(u,v)}
\\
  & \leq  & \n{A^{-n}} \n{\psi_u(u,v)
u_x + \psi_v(u,v) v_x}  \\
    & \leq & \n{A^{-n}} \left[\frac{\n{\psi_u(u,v)}}{\n{v}^{\alpha}}
\n{v}^{\alpha}\n{u_x} + \frac{\n{\psi_v(u,v)}}{\n{u}^{\alpha}} \n{u}^{\alpha} \n{v_x} \right] \\
 & \leq & \n{A^{-n}} \grg_1(\psi) \n{v}^{\alpha}\n{u_x} + \n{A^{-n}}
\grg_2(\psi) \n{u}^{\alpha} \n{v_x} \\
    & \leq & K_{n,1} \grg_1(\psi) \n{y}^{\alpha} + K_{n,2} \grg_2(\psi)
\n{y}^{\alpha} \\
    & \leq & (K_{n,1} + K_{n,2}) \ \dnorm{\psi} \ \n{y}^{\alpha} 
\eqns

Dividing by $\n{y}^{\alpha}$ and taking the supremum over $x,y$ gives 

\[ 
\grg_1(H^n_s(\psi)) < (K_{n,1} + K_{n,2})  \dnorm{\psi}. 
\]

and  proves (\ref{grg1n-Hs-norm-est}).

\noindent {\bf Proof of (\ref{grg2n-Hu-norm-est}):}

As above, $\n{X_{n,y}} = 0$ for $\n{(x,y)} > \delta_n $. 

Thus, 

\bqns
\n{A^{-n}}\n{v}^{\alpha}\n{u_y} & \leq & \n{A^{-n}} (\n{B^n} + \gre)^{\alpha}\n{y}^{\alpha}  \n{X_{n,y}} \\
   & \leq & \n{A^{-n}} (\n{B^n} + \gre)^{\alpha} \delta_n^{\alpha}  \n{X_{n,y}} \\
   & \leq & \n{A^{-n}} (\n{B^n} + \gre)^{\alpha} \delta_n^{\alpha}  Hol(Df_n)  \n{x}^{\alpha} \\
   & =  & K_{n,3} \n{x}^{\alpha}
\eqns

and, using $(\n{A^n}+\gre)^{\alpha} \leq  \n{A^n} + \gre$, we have 

\bqns
\n{A^{-n}} \n{u}^{\alpha}\n{v_y} & \leq & \n{A^{-n}}(\n{A^n} + \gre)^{\alpha} \n{x}^{\alpha} \n{v_y} \\
   & \leq & \n{A^{-n}}(\n{A^n} + \gre) \n{x}^{\alpha} (\n{B^n} + \gre) \\
   & =  & K_{n,4} \n{x}^{\alpha}
\eqns

Hence, 

\bqns 
\n{\partial_yH^n_s(\psi)(x,y)} & = & \n{ A^{-n} \partial_y \psi(u,v)}
\\
  & \leq  & \n{A^{-n}} \n{\psi_u(u,v)
u_y + \psi_v(u,v) v_y}  \\
    & \leq & \n{A^{-n}} \left[\frac{\n{\psi_u(u,v)}}{\n{v}^{\alpha}}
\n{v}^{\alpha} \n{u_y} + 
\frac{\n{\psi_v(u,v)}}{\n{u}^{\alpha}} \n{u}^{\alpha}\n{v_y} \right] \\
 & \leq & K_{n,3} \grg_1(\psi) \n{x}^{\alpha} + \grg_2(\psi) K_{n,4}
\n{x}^{\alpha} \\
 & \leq & (K_{n,3} + K_{n,4}) \  \dnorm{\psi} \ \n{x}^{\alpha} 
\eqns 

Dividing by $\n{x}^{\alpha}$ and taking the supremum over $x,y$ gives 

\[  
\grg_2(H^n_s(\psi)) \leq (K_{n,3} + K_{n,4})  \ \dnorm{\psi},
\] 

and proves (\ref{grg2n-Hu-norm-est}). 

This completes the proof of Proposition \ref{Main_Prop} and gives the desired global $C^1$ linearization $R$ of the diffeomorphism $S$. 

\subsection{ D-H{\"{o}}lder continuity of $R$ near $0$} \label{D-Holder-Cont-of-R}

For a positive real number $\delta > 0$ small, let

\beq \label{pos-orbit}
O_+(\delta) = \bigcup_{n \geq 0} S^n(B_{\delta})
\eeq

and

\beq \label{neg-orbit}
O_-(\delta) = \bigcup_{n \leq 0} S^n(B_{\delta})
\eeq

be the {\it positive} and {\it negative} $S-$orbits of $B_{\delta}(0)$,
respectively,

and let 

\beq \label{full-orbit}
O(\delta) = \bigcup_{n \in \mathbb{Z}} S^n(B_{\delta})
\eeq

be the full $S-$orbit of $B_{\delta}(0)$.

The main result in this subsection is the following proposition. 

\begin{Prop} \label{D-Hol-psi}
Let $0 < \eta < 1$ be as in (\ref{norm-A-B-m-1-eta}) and (\ref{norm-A-B-m-2-eta}), and let $\beta = \alpha(1-\eta)$. Then, for $\delta$ sufficiently small, 

\beq \label{psi-s-Hold}
\mbox{ the map $\psi_s$ is $C^{1,\beta}$ on
$O_+(\delta)$, }
\eeq

and 

\beq \label{psi-u-Hold}
\mbox{ the map $\psi_u$ is $C^{1,\beta}$ on
$O_-(\delta)$.}
\eeq  

Hence, the map $\psi =  (\psi_s, \psi_u)$ is $C^{1,\beta}$ on
$B_{\delta}(0) = O_+(\delta) \bigcap O_-(\delta)$. 

\end{Prop}

Assuming  Proposition \ref{D-Hol-psi}, the proof of Theorem \ref{Main_Theorem_1} is
completed by simply restricting the map $R = (I_s + \psi_s, I_u + \psi_u)$ to the ball
$B_{\delta}(0)$.

Recalling that finite compositions of D-H{\"{o}}lder maps are again
D-H{\"{o}}lder, Proposition \ref{D-Hol-psi} implies 
that the map $R = (I_s + \psi_s, I_u + \psi_u)$  (which satisfies $R =
L^{-n}RS^n$) is $C^{1,\beta}$ on  $S^{-n}(B_{\delta}(0))$  for any finite integer $n$.

But, any bounded subset of $O(\delta)$ is contained in the set

\[ \bigcup_{\mid n \mid \leq N} S^n(B_{\delta})  \]

for some positive integer $N$, and we obtain

\begin{Prop} \label{DHolder-psi-s-u}
  For $\delta > 0$ small enough, the map $R$ is $C^{1,\beta}$ on bounded subsets of $O(\delta)$.
\end{Prop}

 To prove Proposition \ref{D-Hol-psi}, we first show that
 (\ref{psi-s-Hold}) holds.

 Then, as in the proof of (\ref{Hus-aut-norm}), we obtain
 (\ref{psi-u-Hold}) with the same method by 
 replacing $S, A, B$ with  $S^{-1}, B^{-1}, A^{-1}$, respectively.

 Let us use $\psi \mid Z$ to denote the restriction of the function
 $\psi$ to a set $Z$.

 Our proof of (\ref{psi_s_def}) implied  that the series

 \[ \psi_s = \sum_{n \geq 0}  H_s^n(A^{-1}X_1) \]

 converges uniformly in the $C^1$ topology, but, of course, it gives
 no information about D-H{\"{o}}lder convergence. 
 However, observe that the restriction $A^{-1}X_1 \mid O_+(\delta)$ 
is $C^{1,\alpha}$ on $O_+(\delta)$ since it vanishes off
 $B_{\delta}(0)$.  Moreover,
 
\[ \psi_s \mid O_+(\delta) = \sum_{n \geq 0}  H_s^n(A^{-1}X_1 \mid
 O_+(\delta)).  \]

Of course, we also have that $A^{-1}X_1 \mid O_+(\delta)$ is
$C^{1,\beta}$ on $O_+(\delta)$ since $0 < \delta < 1$ and
$0 < \beta <  \alpha$.  

We proceed to introduce a Banach space $\sE^+ = \sE^+(\delta)$ of
$C^{1, \beta}$ functions on  $O_+(\delta)$ with the property that if
$\delta > 0$ is sufficiently small, then function
$H_s(\psi) = A^{-1}\psi  \circ S$ is a well-defined bounded linear self-map of $\sE^+$ satisfying (\ref{Hs-m-less-1}) below. 

 It will follow from this that the series 

 \[ \sum_{n \geq 0}  H_s^n(A^{-1}X_1 \mid
 O_+(\delta)) \]

 converges uniformly in $\sE^+$.  
 
  This  will show that $\psi_s | O_{+}(\delta)$ is $C^{1,\beta}$, thus proving (\ref{psi-s-Hold}).

 \vs

Let us proceed to define the space $\sE^+$.

 For $(x,y), (h,k)  \in E^u \times E^s$, let

 \beq \label{M-star-def}
 M^{\star}_{x,h} = \max(\n{x},\n{h}),
\eeq

\beq \label{N-def}
N_{y,k} = \max(\n{y},\n{k}),
\eeq
 
\beq \label{M-def}
M_{x,h} = \min(1,M^{\star}_{x,y}), 
\eeq

and, let

\beq \label{U-def}
U(x,h,k) = M_{x,h}^{\alpha \eta} \n{(h,k)}^{\alpha(1-\eta)} 
\eeq

\beq \label{V-def}
V(y,h,k) = N_{y,k}^{\alpha \eta} \n{(h,k)}^{\alpha(1-\eta)}.
\eeq

Given a $C^1$ function
$\psi:O_+(\delta) \rarrow E^u$, and points $(x,y) \in O_+(\delta), (h,k) \in E^u \times E^s$, with
$(x+h,y+k)  \in O_+(\delta)$ and $\n{(h,k)} > 0$, consider the real-valued functions 

\beq
 \grg_3(\psi,x,y,h,k) = \frac{\n{\psi_x(x+h,y+k) - \psi_x(x,y)}}{V(y,h,k)}
\eeq

and 

\beq
 \grg_4(\psi,x,y,h,k) = \frac{\n{\psi_y(x+h,y+k) - \psi_y(x,y)}}{U(x,h,k)}
\eeq

and their suprema

\beq \label{norm-grg-3}
\n{\grg_3(\psi)} = \sup_
{
\stackrel{0 < \mid (h,k) \mid}
         {
          \stackrel{(x,y) \in O_+(\delta) }
                   {(x+h,y+k) \in O_+(\delta)}
         }
}
\grg_3(\psi,x,y,h,k).
\eeq

and 

\beq \label{norm-grg-4}
\n{\grg_4(\psi)} = \sup_
{
\stackrel{0 < \mid (h,k) \mid}
         {
          \stackrel{(x,y) \in O_+(\delta) }
                   {(x+h,y+k) \in O_+(\delta)}
         }
}
\grg_4(\psi,x,y,h,k).
\eeq

As in the case of the functions $\grg_1, \grg_2$, we define the functions in $\grg_3$ and $\grg_4$  to have the value 0 if the denominators vanish.

 We define $\sE^+ = \sE^+(\delta)$ to 
be the set of $C^1$ functions $\psi: O_+(\delta) \rarrow E^u$ such
that 

\beq \label{vanish-plus}
\psi(x,y) = 0 \mbox{ for } (x,y) \in O_+(\delta) \bigcap \left(E^u
\times \{0\} \bigcup \{0\} \times E^s \bigcup
B_{\delta}^c\right) 
\eeq

and 

\beq \label{grg-3-4-finite}
\sup(\grg_3(\psi), \grg_4(\psi)) < \infty.
\eeq

Defining the norm of $\psi \in \sE^+$ to be 

\beq \label{norm-psi-O_+}
\dnorm{\psi} = \max_{i = 3,4}(\n{\grg_i(\psi)}),
\eeq

we have that the pair $(\sE^+, \dnorm{\cdot})$ becomes a Banach
space.

For $\psi \in \sE^+$, we use the same formula as in (\ref{H1-def}): 

\[ H_s(\psi) =  A^{-1}\psi \circ S. \]

The proof of (\ref{psi-s-Hold}) is obtained from

\begin{Lem}  \label{Hs-E+-Lemma}
Assume that conditions (\ref{Am-gre-est})-(\ref{norm-A-B-m-2-eta}) are
satisfied. 

Then, for $\delta > 0$ small enough, the map $H_s$ is a well defined
bounded linear  map from $\sE^+$ into itself such that

\beq \label{Hs-m-less-1}
\n{H_s^m} < 1.
\eeq

\end{Lem}

Note that, since $S=L$ off $B_{\delta}(0)$ and $S$ is Lipschitz close
to $L$, we have

\[ S(O_+(\delta) \setminus B_{\delta}(0)) \subset O_+(\delta)
\setminus B_{\delta}(0) \]
for $\delta$ small enough.

This, and the linearity of $S$ on $E^u
\times \{0\} \bigcup \{0\} \times E^s$, imply that, if $\psi$
satisfies (\ref{vanish-plus}), then $\psi \circ S$ also satisifies
(\ref{vanish-plus}).  Hence, the $C^1$ function $H_s(\psi)$ is
well-defined as a map in the function space $\sE^+$.  

We need to show that $\n{H_s} < \infty$ and that (\ref{Hs-m-less-1})
holds.

These will follow from inequality (\ref{dnorm-Hn-est}) below applied in the cases
$n=1$ and $n=m$, repectively. 

Proceeding toward this goal, let 
$n \neq 0$ be a non-zero  integer, let $\delta_n, f_n$
be as in Proposition \ref{global-lin-prop}, and write $f_n = (X_n, Y_n)$ as in
Lemma \ref{Basic-X-Y-Lemma}.  

We will need  the following functions. 

Define

\beq \label{DelX-def}
\Del X_n = \Del X_n(x,y,h,k) = X_n(x + h, y+k) - X_n(x,y),
\eeq

\beq \label{DelY-def}
\Del Y_n = \Del Y_n(x,y,h,k) = Y_n(x + h, y+k) - Y_n(x,y),
\eeq

\beq \label{DelX_nx-def}
\Del X_{n,x} = \Del X_{n,x}(x,y,h,k) = X_{n,x}(x+h, y+k) - X_{n,x}(x,y),
\eeq

\beq \label{DelX_ny-def}
\Del X_{n,y} = \Del X_{n,y}(x,y,h,k) = X_{n,y}(x+h, y+k) - X_{n,y}(x,y),
\eeq

\beq \label{DelY_nx-def}
\Del Y_{n,x} = \Del Y_{n,x}(x,y,h,k) = Y_{n,x}(x+h, y+k) - Y_{n,x}(x,y),
\eeq

\beq \label{DelY_ny-def}
\Del Y_{n,y} = \Del Y_{n,y}(x,y,h,k) = Y_{n,y}(x+h, y+k) - Y_{n,y}(x,y),
\eeq

\begin{Lem} \label{MN-Xn-Yn-Lemma}
Let $0 < \alpha < 1, \ 0 < \eta < 1$, $n \in \{1,m\}$, and let $\gre$ satisfy conditions 
(\ref{Am-gre-est})-(\ref{norm-A-B-m-2-eta}).   Let $X_n, Y_n$ be as in Proposition \ref{Basic-X-Y-Lemma}, and let
$\delta = \delta_n \in (0,\min(\gre,1/3))$ be small enough that

\beq \label{delta-cond-1}
Hol(Df_n)\delta^{\alpha} < \frac{\gre}{3}.
\eeq

Then, for  $(x,y) \in O_+(\delta)$, $(x+h, y+k) \in O_+(\delta)$, and $\n{(x,y)} < \delta$, we have 

\beq \label{max-MN-DelX-DelY}
\max(\n{\Del X_n}, \n{\Del Y_n}) \leq  \gre  \min(\n{(h,k)},M_{x,h},
N_{y,k}), 
\eeq

\beq \label{max-V-Xnx-Ynx}
\max(\n{\Del X_{n,x}}, \n{\Del Y_{n,x}})  \leq 3 Hol(Df_n) \min  \left(\n{(h,k)}^{\alpha} , 
V(y,h,k) \right)
\eeq

\beq \label{max-U-Xny-Yny}
\max( \n{\Del X_{n,y}}, \n{\Del Y_{n,y}})  \leq 3 Hol(Df_n) \min \left( \n{(h,k)}^{\alpha}, U(x,h,k)\right)
\eeq

\end{Lem} 

{\bf Proof.}
We only give the arguments for $\Del X_n, \ \Del X_{n,x}$ and $\Del X_{n,y}$  leaving the similar arguments for
$\Del Y_n, \ \Del Y_{n,x}, \ \Del Y_{n,y}$ to the  reader.

 For notational convenience, let $\gre_1 = \frac{\gre}{3}$.

Since $X_n(x,0) = 0$ for each $(x,0)$ we also have that $X_{n,x}(x,0) =0 $ for
each $(x,0)$.  Similarly,  $X_n(0,y) = 0$ implies that $X_{n,y}(0,y) =
0$ for each $(0,y)$.  

Letting $n=1$ or $n=m$, (\ref{delta-cond-1}) and the vanishing of $X_n$ off $B_{\delta}(0)$ imply that 

\beq \label{Lip-fn-gre-div-3}
Lip(X_n) \leq Hol(Df_n)\delta^{\alpha} <
\gre_1. 
\eeq

From this, the  Mean Value Theorem gives 


\bqns
 \n{\Del X_n} & = & \n{X_n(x+h,y+k) - X_n(x,y)} \\
  & \leq & \n{X_n(x+h,y+k) -X_n(x,y+k)} + \n{X_n(x,y+k) - X_n(x,y)} \\
  & \leq & \sup_{0 \leq t \leq 1} \n{X_{n,x}(x + th, y+k)}\n{h} + \sup_{0 \leq t
\leq 1} \n{X_{n,y}(x,y+tk)}\n{k} \\
 &  \leq & 2 \gre_1 \n{(h,k)}.
\eqns

In addition, using the Mean Value Theorem again and the fact that
$X_n(0,y) = X_n(x,0) = 0$ for all $x,y$, we have

\bqns
\n{ \Del X_n} & = & \n{X_n(x+h, y+k) - X_n(x,y)}  \\
        & = & \n{X_n(x+h, y+k) -X_n(x+h,y) + X_n(x+h,y) - X_n(x,y)} \\ 
 & \leq &  (\sup_{0 \leq t \leq 1} \n{X_{n,y}(x+h,y+tk)})\n{k}
 + \n{X_n(x+h,y) - X_n(x,y)} \\
        & \leq & \gre_1 \n{k} + \n{X_n(x+h,y) - X_n(x+h,0)} + \n{X_n(x,y)-X_n(x,0)} \\
        & \leq & \gre_1 \n{k} + \gre_1 \n{y} + \gre_1 \n{y} \\
        & \leq & 3\gre_1 \max(\n{y}, \n{k}) \\
        & \leq & 3\gre_1 N_{y,k} \\
        & = & \gre N_{y,k}. 
\eqns

Similarly, 

\bqns
\n{ \Del X_n}& = & \n{X_n(x+h, y+k) - X_n(x,y)}  \\ 
        & = & \n{X_n(x+h, y+k) -X_n(x,y+k) + X_n(x,y+k) - X_n(x,y)}  \\
 & \leq & (\sup_{0 \leq t \leq 1} \n{X_{n,x}(x+th,y+k})\n{h} +  \n{X_n(x,y+k) - X_n(0,y+k)} + \n{X_n(x,y)-X_n(0,y)}  \\
        & \leq & \gre_1 \n{h} + \gre_1 \n{x} + \gre_1 \n{x} \\
         & \leq & 3 \gre_1 \max(\n{x}, \n{h}) \\
         & \leq &  \gre M^{\star}_{x,h}. 
 \eqns

 Since $M_{x,h} = \min(1,M^{\star}_{x,h})$, to get

\beq \label{Del-X-n-est}
\n{ \Del X_n}  \leq \gre M_{x,h}, 
 \eeq

 it suffices to show that if $M^{\star}_{x,h} > 1$, then 

 \beq \label{Del-X-n-est-2}
\n{ \Del X_n}  \leq  \gre. 
 \eeq

 But, if $M^{\star}_{x,h} > 1$, then, since $\n{(x,y)} < \delta < \half$,
 we have

 \[ \n{(x+h, y+k)}  \geq \n{x+h} > \n{h} - \delta > 1 - \delta >
 \delta, \] 
 giving  $X_n(x+h,y+k) = 0$, and

\[ \n{ \Del X_n} = \n{X_n(x+h, y+k) - X_n(x,y)} = \n{X_n(x,y)} \leq \gre_1 
 \n{(x,y)} \leq \gre_1 \delta \leq  \gre. \] 
 
Moving to (\ref{max-V-Xnx-Ynx}), we proceed as follows. 

First, since $Hol(DX_n) \leq Hol(Df_n)$, we  have 

\beq \label{Hol-DXn}
\n{\Del X_{n,x}} \leq Hol(Df_n) \n{(h,k)}^{\alpha}.
\eeq

Next, using (\ref{max_X_x_Y_x_gra_n})
at the points $(x+h,y)$ and $(x,y)$, we have 

\beq
\n{X_{n,x}(x+h,y)  -X_{n,x}(x,y)}  \leq 2 Hol(Df_n) \n{y}^{\alpha} 
\eeq

for each $x, y, h$.  

Hence, 

\bqn
\n{\Del X_{n,x}} & = & \n{X_{n,x}(x+h,y+k) -X_{n,x}(x,y)} \nonumber \\
  & \leq & \n{X_{n,x}(x+h,y+k) -X_{n,x}(x+h,y)} \nonumber \\
    & & + \n{X_{n,x}(x+h,y)  -X_{n,x}(x,y)} \nonumber \\ 
    & \leq & Hol(Df_n) \n{k}^ {\alpha} + 2 Hol(Df_n) \n{y}^{\alpha} \nonumber \\
     & \leq & 3 Hol(Df_n) N_{y,k}^{\alpha} \label{Del-X-1} 
\eqn

Putting this together with (\ref{Hol-DXn}) gives

\bqn
\n{\Del X_{n,x}} & =  & \n{\Del X_{n,x}} ^{\eta + 1 - \eta} \nonumber \\
   & \leq & \left[3 Hol(Df_n)N_{y,k}^{\alpha}\right]^{\eta} \left[Hol(Df_n)
\n{(h,k)}^{\alpha}\right]^{(1-\eta)} \nonumber \\ 
   & \leq & Hol(Df_n) \left(3 N_{y,k}^{\alpha}\right)^{\eta} \left(\n{(h,k)}^{\alpha}\right)^{(1-\eta)}  \nonumber \\ 
   & = & 3^{\eta} Hol(Df_n) N_{y,k}^{\alpha \eta} \n{(h,k)}^{\alpha(1-\eta)} \nonumber 
\\
   & \leq & 3 Hol(Df_n)) V(y,h,k),  \label{Del-Xn-x}
\eqn

and  this proves 
(\ref{max-V-Xnx-Ynx}).

The proof of (\ref{max-U-Xny-Yny}) is similar to that for (\ref{max-V-Xnx-Ynx}), except we have to use $M_{x,h} = \min(1, M^{\star}_{x,h})$ instead of $N_{y,k}$.

Working with $X_{n,y}$ we have

\beq \label{Del-Xn-y}
\n{\Del X_{n,y}} \leq Hol(Df_n) \n{(h,k)}^{\alpha},
\eeq

\[\n{X_{n,y}(x,y+k)  -X_{n,y}(x,y)}  \leq 2 Hol(Df_n) \n{x}^{\alpha},  \]

and 

\bqns
\n{\Del X_{n,y}} & = & \n{X_{n,y}(x+h,y+k) -X_{n,y}(x,y)} \\
  & \leq & \n{X_{n,y}(x+h,y+k) -X_{n,y}(x,y+k)} \\ 
    & & + \n{X_{n,y}(x,y+k)  -X_{n,y}(x,y)} \\ 
    & \leq & Hol(Df_n) \n{h}^ {\alpha} + 2 Hol(Df_n) \n{x}^{\alpha} \\
     & \leq & 3 Hol(Df_n) (M^{\star}_{x,h})^{\alpha}.
\eqns

As in the proof of (\ref{Del-X-n-est-2}), if $M^{\star}_{x,h} > 1$, then $\n{X_{n,y}(x+h,y+k)} = 0, \ M_{x,h} = 1$, and 

\bqns
\n{\Del X_{n,y}} & = & \n{X_{n,y}(x,y)}  \leq  Hol(Df_n) \delta^{\alpha} \\
& \leq  & 3 Hol(Df_n)  \\ 
& \leq & 3 Hol(Df_n) M_{x,h}^{\alpha} \label{Del-X-3}.
\eqns

Now, replacing $N_{y,k}$ by $M_{x,h}$ in the technique used to get 
(\ref{Del-Xn-x}), we obtain the statement involving $\n{\Del X_{n,y}}$
in (\ref{max-U-Xny-Yny}), and this completes the proof of Lemma
\ref{MN-Xn-Yn-Lemma}. 

\vs

Now, set 

\beq \label{tKn1}
\tilde{K}_{n,1} = \n{A^{-n}}(\n{A^n}+\gre)(\n{B^n}+ \gre)^{\alpha \eta}(\n{A^n} + \gre)^{\alpha(1-\eta)},
\eeq

\beq \label{tkn2}
\tilde{K}_{n,2} = \n{A^{-n}}(\n{A^n}+ \gre) \delta^{\alpha(1-\eta)} Hol(Df_n),
\eeq

\beq \label{tkn3}
\tilde{K}_{n,3} = \n{A^{-n}}(\n{B^n} + \gre)^{\alpha \eta}(\n{A^n} + \gre)^{\alpha(1-\eta)} Hol(Df_n) \delta^{\alpha(1-\eta)}, 
\eeq

\beq \label{tkn4}
\tilde{K}_{n,4} = \n{A^{-n}}(\n{A^n} + \gre) (\n{B^n} + \gre),
\eeq

and

\beq \label{tgrtn}
\tilde{\grt}_n = \max(\tilde{K}_{n,1} + \tilde{K}_{n,2}, \tilde{K}_{n,3} + \tilde{K}_{n,4}).
\eeq

Observe that we can choose $\eta$ close enough to 1 so that 

\[ \tilde{K}_{m,4} < \tilde{K}_{m,1} < 1. \]

Then, we can choose $\delta$ small enough, depending on $1-\eta$, such that 

\beq \label{tkn1+tkn2}
\tilde{K}_{m,1} + \tilde{K}_{m,2} < 1, 
\eeq

and

\beq \label{tkn3+tkn4}
\tilde{K}_{m,3} + \tilde{K}_{m,4} < 1.
\eeq

This gives

\beq \label{tkn-est}
\tilde{\grt}_m < 1.
\eeq

Next, we define some new functions and make some more  estimates.

Write  

\[ u = u_n = A^n x + X_n(x,y), \ v = v_n = B^ny + Y_n(x,y), \]

\[  u_{1,n} = A^n(x+h) + X_n(x+h,y+k), \mbox{ and } v_{1,n} = B^n(y+k) + Y_n(x+h,y+k) \]

so that

\[ 
T^n(x,y) = (u_n, v_n) \mbox{ and } T^n(x+h, y+k) = (u_{1,n}, v_{1,n}).
\]

Then, setting 
\[ \Del u = u_{1,n} - u_n  = A^nh + \Del X_n  \]

and

\[ \Del v = v_{1,n}  -v_n  = B^nk + \Del Y_n, \]

we have, for $\psi \in \sE^+$,

\beq \label{Hn-psi-formula} 
H^n_s(\psi)(x+h,y+k) - H^n_s(\psi)(x,y)  =  A^{-n}\psi(u+\Del u, v+\Del v) - A^{-n}\psi(u,v).
\eeq

The proof of Lemma \ref{Hs-E+-Lemma}  involves estimating the $x$ and
$y$ partial derivatives of the left side of (\ref{Hn-psi-formula}) in
terms of those of  $\psi$ and will be concluded in the inequalities (\ref{grg-3-est-1}) and (\ref{grg-4-est-1}) below.

First,  we have the estimates

\bqns
\n{\Del u} & = & \n{A^nh + \Del X_n} \\
 & \leq & (\n{A^n}\n{h} + \gre \min(M_{x,h}, N_{y,k}, \n{(h,k)}) \\
 & \leq & (\n{A^n}\n{h} + \gre \min(M^{\star}_{x,h}, N_{y,k}, \n{(h,k)}) \\
 & \leq & (\n{A^n} + \gre) \min(M^{\star}_{x,h},\n{(h,k)}),
\eqns

and

\bqns
\n{\Del v} & = & \n{B^nk + \Del Y_n} \\
 & \leq & \n{B^n}\n{k} + \gre \min(M_{x,h}, N_{y,k},\n{(h,k)}) \\
& \leq & (\n{B^n} + \gre) \min(N_{y,k},\n{(h,k)}). 
\eqns

Putting these together with 

\[ \n{u} = \n{A^nx + X_n(x,y)} \leq (\n{A^n} + \gre_1)\n{x} \leq (\n{A^n} + \gre_1)M^{\star}_{x,h} \]

and 

\[ \n{v} = \n{B^ny + Y_n(x,y)} \leq (\n{B^n} + \gre_1)\n{y} \leq (\n{B^n} + \gre_1)N_{y,k}, \] 

we have 

\beq \label{u-Delu}
M^{\star}_{u,\Del u} \leq (\n{A^n} + \gre)M^{\star}_{x,h},
\eeq

\beq \label{v-Delv}
N_{v,\Del v} \leq (\n{B^n} + \gre)N_{y,k},
\eeq

and 

\beq \label{Delu-Delv}
\n{(\Del u, \Del v)} \leq (\n{A^n} + \gre) \n{(h,k)}.
\eeq

From (\ref{u-Delu}), we have

\beq \label{min-u-Del-u-a}
\min(1,M^{\star}_{u,\Del u}) \leq \min(1, (\n{A^n} + \gre)M^{\star}_{x,h}).
\eeq

Observe that, if $a, b, c$ are non-negative real numbers and $c \geq 1$, then

\beq \label{min-1-abc}
min(a, bc) \leq c \cdot min(a,b).
\eeq

Applying this with $a = 1, \ b = M^{\star}_{x,h}$ and $c = \n{A^n} +
\gre$, we get 

\[ \min(1,M^{\star}_{u,\Del u}) \leq (\n{A^n} + \gre) \min(1,M^{\star}_{x,h}), \] 

which, by definition, is 

\beq \label{M-u-Del-u}
M_{u,\Del u} \leq (\n{A^n} + \gre) M_{x,h}.
\eeq

Consider the compositions of $U(x,h,k)$ and $V(y,h,k)$ with $u, v,
\Del u, \Del v$:

\beq \label{U1-def}
U_1 = U_1(x,h,k) = U(u,\Del u, \Del v) = M_{u,\Del u}^{\alpha \eta}\n{(\Del u, \Del v)}^{\alpha(1-\eta)}, 
\eeq

\beq \label{V1-def}
V_1 = V_1(y,h,k) = V(v, \Del u, \Del v) = N_{v,\Del v}^{\alpha \eta}\n{(\Del u, \Del v)}^{\alpha(1-\eta)}.
\eeq

Using $(\n{A^n} + \gre) > 1$, we have 

\beq \label{U1-est}  
U_1  \leq (\n{A^n} + \gre)^{\alpha \eta} (\n{A^n} +
\gre)^{\alpha(1-\eta)} U = (\n{A^n} + \gre)^{\alpha})U \leq (\n{A^n} + \gre)U
\eeq

and 
\beq \label{V1-est}  
V_1  \leq  (\n{B^n} + \gre)^{\alpha \eta} (\n{A^n} +
\gre)^{\alpha(1-\eta)} V .
\eeq

Let us observe that, with $0 < \delta \leq  \half$, we have 

\beq \label{Mxh-Nyk-leq-1}
 \max(M_{x,h}, N_{y,k}) \leq 1. 
\eeq

Indeed, the definition of $M_{x,h}$ makes it no larger than 1.  If $N_{y,k} > 1$, then, since both $(x,y)$ and $(x+h, y+k)$ are in $O_+(\delta)$, we would get 
that $\n{y} < \delta \leq  \half$, which, in turn, would imply that $N_{y,k} = \n{k} > 1$.

This would give 

\[ \delta > \n{y+k} \geq  \n{k} - \n{y} >  1 - \delta, \]

contradicting the condition that $ \delta \leq  \half$. 

Now, we have $v_x = Y_{n,x}$ and $u_y = X_{n,y}$. 

Since $\max(\n{X_{n,y}}, \n{Y_{n,x}}) = 0$ when $\n{(x,y)} > \delta$,  we have that 

\[ U \n{v_x} >  0 \mbox{ or } V \n{u_y} > 0  \Rightarrow \n{(x,y)} \leq \delta.  \]

This, together with (\ref{Mxh-Nyk-leq-1}), gives

\bqns 
\n{A^{-n}} U_1 \n{v_x} & = & \n{A^{-n}}U_1 \n{Y_{n,x}} \\
    & \leq  & \n{A^{-n}}(\n{A^n} +  \gre) U \n{Y_{n,x}} \\
    & \leq & \n{A^{-n}}(\n{A^n} +  \gre) 
M_{x,h}^{\alpha \eta}\n{(h,k)}^{\alpha(1 - \eta)} Hol(Df_n)  \n{y}^{\alpha} \\
    & \leq & \n{A^{-n}}(\n{A^n} +  \gre)  Hol(Df_n) \n{y}^{\alpha} \n{(h,k)}^{\alpha(1 - \eta)} \\
    & = & \n{A^{-n}}(\n{A^n} +  \gre)  Hol(Df_n) \n{y}^{\alpha \eta} \n{y}^{\alpha(1- \eta)} \n{(h,k)}^{\alpha(1 - \eta)} \\
    & \leq & \n{A^{-n}}(\n{A^n} +  \gre)  Hol(Df_n) \n{y}^{\alpha \eta} \delta ^{\alpha(1- \eta)} \n{(h,k)}^{\alpha(1 - \eta)} \\
    & = & \tilde{K}_{n,2}  \n{y}^{\alpha \eta}  \n{(h,k)}^{\alpha(1 - \eta)} \\
    & \leq & \tilde{K}_{n,2}  N_{y,k}^{\alpha \eta} \n{(h,k)}^{\alpha(1 - \eta)} \\
    & \leq  & \tilde{K}_{n,2}  V.
\eqns

In addition,

\bqns 
 V \n{u_y} & = & V \n{X_{n,y}} \\
 & \leq & N_{y,k}^{\alpha \eta} \n{(h,k)}^{\alpha(1-\eta)} Hol(Df_n)
\n{x}^{\alpha} \\
    & \leq & Hol(Df_n) \n{x}^{\alpha} \n{(h,k)}^{\alpha(1 - \eta)} \\
    & = & Hol(Df_n) \n{x}^{\alpha \eta} \n{x}^{\alpha(1 - \eta)} \n{(h,k)}^{\alpha(1 - \eta)} \\
    & \leq & Hol(Df_n) \delta^{\alpha(1-\eta)} M_{x,h}^{\alpha \eta} \n{(h,k)}^{\alpha(1 - \eta)} \\
    & = & Hol(Df_n) \delta^{\alpha(1-\eta)} U
\eqns

which implies 

\[ \n{A^{-n}} V_1 \n{u_y} \leq \n{A^{-n}} (\n{B^n} + \gre)^{\alpha \eta}(\n{A^n} + \gre)^{\alpha(1-\eta)} V \n{u_y}  \leq \tilde{K}_{n,3} U.  \]

We now are in a position to prove that

\beq \label{grg-3-est-1}
\n{\grg_3(H^n_s(\psi))} \leq \tilde{\grt}_n  \dnorm{\psi}
\eeq

and

\beq \label{grg-4-est-1}
\n{\grg_4(H^n_s(\psi))} \leq \tilde{\grt}_n  \dnorm{\psi}.
\eeq

These inequalities will imply that

\beq \label{dnorm-Hn-est}
\dnorm{H^n_s(\psi)} = \max\left(\grg_3(H^n_s(\psi)), \grg_4(H^n_s(\psi))\right) \leq \tilde{\grt}^n \dnorm{\psi}.
\eeq

As we mentioned above, applying this, respectively, for $n=1$ and
$n=m$, will prove Lemma \ref{Hs-E+-Lemma}, and, hence, complete the
proof of Theorem \ref{Main_Theorem_1}.

Let us proceed to the proofs of (\ref{grg-3-est-1}) and (\ref{grg-4-est-1}). 

For $(x,y) \in O_+(\delta), \ (h,k) \in E^u \times E^s$ such that $(x+h, y+k) \in O_+(\delta)$ and $\n{(h,k)} > 0$,  let

\bqns \Del H^n_s(\psi)_x & = & \partial_x( H^n_s(\psi))(x+h, y+k) - \partial_x(H^n_s(\psi)) (x,y)  \\
                         & = &  \partial_x( H^n_s(\psi)(x+h, y+k) - H^n_s(\psi)(x,y))
\eqns

and

\bqns \Del H^n_s(\psi)_y & = & \partial_y( H^n_s(\psi))(x+h, y+k) - \partial_y(H^n_s(\psi)) (x,y)  \\
                         & = &  \partial_y( H^n_s(\psi)(x+h, y+k) - H^n_s(\psi)(x,y)).
\eqns

  Then, from (\ref{Hn-psi-formula}), we have 
  
\bqns
 \n{\Del H^n_s(\psi)_x} & \leq  & \n{A^{-n}} \n{\partial_x(\psi(u+\Del
   u, v + \Del v) - \psi(u,v))} \\
 & \leq & \n{A^{-n}} \n{\partial_u (\psi(u + \Del u, v + \Del v)
  -  \psi(u,v))} \n{u_x} \\
 & & \ \ \ +  \n{A^{-n}} \n{\partial_v (\psi(u + \Del u, v + \Del v) -
  \psi(u,v))} \n{v_x} \\
 & \leq & \n{A^{-n}} \grg_3(\psi) V_1 \n{u_x} + \n{A^{-n}} \grg_4(\psi) U_1 \n{v_x} \\
 & \leq & \n{A^{-n}} \grg_3(\psi) V_1 (\n{A^n} + \gre) + \grg_4(\psi) \tilde{K}_{n,2} V \\
 & \leq & \n{A^{-n}} (\n{A^n}
+ \gre)\grg_3 (\psi)(\n{B^n} + \gre)^{\alpha \eta}(\n{A^n} + \gre)^{\alpha(1-\eta)}  V +
\grg_4(\psi) \tilde{K}_{n,2} V \\
 & \leq & \tilde{K}_{n,1} \grg_3(\psi) V + \tilde{K}_{n,2}
\grg_4(\psi) V \\
& \leq  & \tilde{\grt}_n  \dnorm{\psi} V.
\eqns

Dividing by $V$ and taking the supremum as $(x,y)$ and $(x+h, y+k)$ vary gives (\ref{grg-3-est-1}). 

Similarly, 

\bqns
 \n{\Del H^n(\psi)_y} & \leq  & \n{A^{-n}} \n{\partial_y(\psi(u+\Del
   u, v + \Del v) - \psi(u,v))} \\
 & \leq & \n{A^{-n}} \n{\partial_u (\psi(u + \Del u, v + \Del v)
  -  \psi(u,v))} \n{u_y} \\
 & & \ \ \ +  \n{A^{-n}} \n{\partial_v (\psi(u + \Del u, v + \Del v) -
  \psi(u,v))} \n{v_y} \\
 & \leq & \n{A^{-n}} \grg_3(\psi) V_1 \n{u_y} + \n{A^{-n}} \grg_4(\psi) U_1 \n{v_y} \\
 & \leq & \tilde{K}_{n,3} \grg_3(\psi) U  + \n{A^{-n}} \grg_4(\psi) U_1 \n{v_y} \\
 & \leq & \tilde{K}_{n,3} \grg_3(\psi) U  + \n{A^{-n}} \grg_4(\psi) (\n{A^n} +
\gre)U (\n{B^n} + \gre) \\
 & \leq &  \tilde{K}_{n,3} \grg_3(\psi)  U  + \tilde{K}_{n,4} \grg_4(\psi) U \\
 & \leq & \tilde{\grt}_n \dnorm{\psi} U
\eqns

Dividing by $U$ and, again, taking the supremum as $(x,y)$ and $(x+h, y+k)$ vary gives (\ref{grg-4-est-1}). 

Observe that, in the above estimates, it may be assumed that $\n{(\Del u, \Del v)} > 0$. Otherwise, the given inequalities would trivially be satisfied, since $\Del H^n(\psi)_x$ and $\Del H^n(\psi)_y$ would vanish.

\section{Continuous Dependence on Parameters} \label{par-dep}

In this section, we formulate a version of Theorem
\ref{Main_Theorem_1} for continuous families of maps.  The proofs use
more or less standard methods and will only be sketched.

Let $\Grl$ be a topological space, and let $X$ be a complete metric
space.  Given a map $\Grf:\Grl \times X \rarrow X$ and a point $\grl
\in \Grl$, define the $\grl$-section map $\Grf_{\grl}:
X \rarrow X$ to be the map given by 

\[ \Grf_{\lambda}(x) = \Grf(\grl,x) \mbox{ for } x \in X. \]

A map $\Grf:\Grl \times X \rarrow X$ is called a {\it uniform contraction
map} on $\Grl \times X$ if it is continuous and there is a constant $0
< \mu < 1$ such that each $\grl$-section $\Grf_{\grl}$ is a contraction
map with Lipschitz constant less than $\mu$. 

The following theorem is well-known and easy to prove.

\begin{Theorem} Let $\Grf: \Grl \times X \rarrow X$ be a uniform
contraction map on $\Grl \times X$. For each $\grl \in \Grl$, let
$p_{\grl}$ be the unique fixed point of the $\grl$-section map
$\Grf_{\grl}$. Then, the map $\grl \rarrow p_{\grl}$ is continuous. 
\end{Theorem}

For $0 < \alpha < 1$, and a $C^1$ function $f:U \rarrow E$, define

\[ 
\n{f}_{1,\alpha} = \sup_{x \neq y, x, y \in U} \max\left(\n{f(x)}, \n{Df(x)}, \frac{\n{Df(x) -
Df(y)}}{\n{x-y}^{\alpha}}\right),
\]

and let $C^{1,\alpha}(U,E)$ be the set of $C^1$ functions from
$U$ into $E$ such that

\[ \n{f}_{1,\alpha} < \infty. \]

Clearly, the quantity $\n{f}_{1,\alpha}$ defines a norm in $C^{1,\alpha}(U,E)$ making it into a
real Banach space. 

We now state the parametrized version of Theorem
\ref{Main_Theorem_1} 

\begin{Theorem} \label{par-Main-Theorem}
Let $0 < \alpha < 1$, let $E$ be a $C^{1,\alpha}$ Banach
space, $p \in E$, and let $U$ be an open neighborhood of $p$
 in $E$. Let $\Grl$ be a topological space, 
let $\grl_0$ be a point in $\Grl$, and let $F: \Grl \times
C^{1,\alpha}(U,E)$ be a continuous map such that $p = p_{\grl_0}$ is an
$\alpha-$hyperbolic fixed point of $F_{\grl_0}$. 

Then, there are neighborhoods $\sL$ of $\grl_0$ in $\Grl$ and $V \subset
U$ of $p$ in $E$ and a real number $\beta \in (0, 1)$ 
such that for each $\grl \in \sL$, the map $F_{\grl}$ has a unique
$\alpha-$hyperbolic fixed point $p_{\grl}$ in $V$ and there is a $C^{1,\beta}$ 
diffeomorphism $R_{\grl}$ from $V$ onto a neighborhood of $0$ in $E$ such that

\beq
DF_{\grl}(p_{\grl})(R_{\grl}(x)) = R_{\grl}(F_{\grl}(x)) \mbox{ for } x \in V \bigcap
F_{\grl}^{-1} V. 
\eeq

Moreover, the map $\grl \rarrow (p_{\grl}, R_{\grl})$ is a continuous
map from $\sL$ into the product space $V \times C^{1,\beta}(V,E)$. 

\end{Theorem}

{\bf Sketch of Proof.}

Replacing $F_{\grl_0}$ by $F_{\grl_0}(x+p) - p$, we may assume that
$p=0$, and $U$ is an open connected neighborhood of $0$. 

Letting $L = DF_{\grl_0}(0)$, we write

\[ F_{\grl_0}(x) = Lx + f_{\grl_0}(x) \]

with $f_{\grl_0}(0) = 0$ and $Df_{\grl_0}(0) = 0$.

The definitions and statements below require that $\grl$ be close to $\grl_0$ in $\Grl$, and we assume this without further mention.

Let $f_{\grl}$ be the $C^{1,\alpha}$ map defined
by 

\beq \label{f-grl-def}
f_{\grl} = F_{\grl} -L.
\eeq

Since the operator $L$ is hyperbolic, the operator $I - L$ has a
bounded inverse.

The equation for a fixed point $x$ of
$F_{\grl}$ has the following forms

\[ x = Lx + f_{\grl}(x), \]

\[ (I - L)x = f_{\grl}(x), \]

or 

\[ x = (I - L)^{-1}f_{\grl}(x). \]

Let $\bar{B}_{\delta}$ denote the closed ball of radius $\delta$ about
$0$ in $E$.  

Let $G_{\grl} = (I - L)^{-1}f_{\grl}$. 

Given $\gre_1 \in (0,1)$, choose $\delta > 0$ such that,
$\bB_{\delta} \subset U$, and 

\[ \n{DG_{\grl_0}} = \n{(I - L)^{-1} Df_{\grl_0}(x)}  < \gre_1, \]

for $x \in \bB_{\delta}$. 

Next, let $\delta_1 = (1 - \gre_1)\delta$ and  let 
$\sL$ be a neighborhood of $\grl_0$ in $\Grl$ so that, for $\grl
\in \sL$ and $x \in \bB_{\delta}$, we have 

\beq \label{eqn-11}
\n{G_{\grl}(0)} = \n{(I - L)^{-1}f_{\grl}(0)} < \delta_1,
\eeq

and 

\beq \label{eqn-10}
\n{DG_{\grl}(x)} = \n{(I - L)^{-1} Df_{\grl}(x)}  < \gre_1.
\eeq

From (\ref{eqn-10}), we see that, for $\grl \in \sL$, 

\beq \label{eqn-12}
Lip(G_{\grl}, \bB_{\delta}) \leq \gre_1.
\eeq

Now, for each $\grl \in \sL$ and each $x \in \bB_{\delta}$, we have 

\bqns
 \n{G_{\grl}(x)} & = & \n{G_{\grl}(x) - G_{\grl}(0) + G_{\grl}(0)} \\
 & \leq & \gre_1 \n{x} + \delta_1 \\
 & \leq  & \gre_1 \delta + (1 - \gre_1)\delta \\
 & \leq   & \delta,
\eqns

This and (\ref{eqn-12}) show that the map $(\grl,x) \rarrow G_{\grl}(x)$ is a uniform contraction map
on $\sL \times \bar{B}_{\delta}$.

Hence, there is a unique fixed point $p_{\grl}$ of $G_{\grl}$ in
$\bB_{\delta}(0)$ which depends continuously on $\grl \in \sL$. 

Since the set $Hyp_{\alpha}(E,E)$ of $\alpha-$hyperbolic automorphisms of
$E$ is an open subset of the set $L(E,E)$ of bounded linear maps on
$E$, we may assume that the point $p_{\grl}$ is an $\alpha-$hyperbolic fixed
point of $F_{\grl}$.

Conjugating $F_{\grl}$ with the translations $x \rarrow x + p_{\grl}$
we may assume that $F_{\grl}(0) = 0$ for each $\grl$ near $\grl_0$. 

Letting $L_{\grl}$ be the derivative of $F_{\grl}$ at $0$, it can be
shown that the subspace $E^u_{\grl}$ is the graph of a bounded linear
map $P^u_{\grl}:E^u_{\grl_0} \rarrow E^s_{\grl_0}$ of small norm. The map
$\grl \rarrow P^u_{\grl}$ is continuous with the obvious topologies. 
A similar statement holds for $E^s_{\grl}$.  

Next, choose coordinates so that the hyperbolic splitting 

\[ E^u_{\grl_0} \times E^s_{\grl_0} \]

coincides with the coordinate subspaces 

\beq \label{constant-splitting}
 E^u \times \{0\} \times \{0\} \times E^s.
\eeq

  For $\grl$ near $0$, there is a $\grl-$dependent family of linear coordinate changes $\{A_{\lambda} \}$ near the identity, such that 
  $A_{\lambda}L_{\grl}A_{\lambda}^{-1}$ preserves the splitting (\ref{constant-splitting})  for each $\grl$.  Replacing $L_{\lambda}$ with
  $A_{\lambda}L_{\grl}A_{\lambda}^{-1}$, we may assume that $L_{\lambda}$ preserves this splitting. 

Now, all of the constructions in the proof of Theorem \ref{Main_Theorem_1}
vary continuously with $\grl$. Here we use the
Irwin method as in the paper of de la Llave and Wayne
\cite{Llave-Wayne} to get the local stable and unstable manifolds via
the Implicit Function Theorem.  Thus, these constructions vary
continuously with $\grl$. 

This leads to  continuously varying linearization operators  $\grl
\rarrow H_{\grl}$ which are contractions on appropriate function spaces as in the proofs of 
(\ref{Hs-aut-norm}), (\ref{Hus-aut-norm}), (\ref{psi-s-Hold}), and (\ref{psi-u-Hold}).

This gives a family $R_{\grl}$ of $C^{1,\beta}$ linearizations
at $p(\grl)$ of $T_{\grl}$ depending continuously on $\grl$ as required for
Theorem \ref{par-Main-Theorem}. 

\section{Vector Fields and Motivation} \label{Motivation}

Our main motivation for the results presented here is concerned with
the study of vector fields. 

For simplicity, we consider the finite dimensional case. 

In this section $r$ will denote a real number larger than 1.

Let $M$ be a $C^{r+1}$ finite dimensional real manifold, and let $X$ be a $C^r$ vector field defined in
an open subset $U \subset M$. Let $p$ be a critical point of $X$;
i.e., $X(p) = 0$. 
  Assume that the local flow $\grf(t,x)$ associated to $X$ is
defined on the product $(-2,2) \times U$, and let $T(x) = \grf(1,x)$
be the time one map of $X$. We choose an open subset $U_1$ of $U$ so
that $T$ and $T^{-1}$ are defined and $C^r$ on $U_1$. 

For $\alpha \in (0,1)$, and each of the properties $\alpha-$contracting,
$\alpha-$expanding, $\alpha-$hyperbolic, bi-circular, etc., we define $p$
to have the corresponding property if it has that property as a fixed
point of $T$. 

For instance, if $M$ is the Euclidean space $R^N$, with $N$ a positive
integer, then, using that
$exp(DX_p) = DT_p$, it can be seen that  $p$ is bi-circular if the set consisting of the real parts of the eigenvalues
of $DX_p$ is a two element set $\{a, b\}$ with $a < 0 < b$. 

A well-known technique of Sternberg \cite{Sternberg-Loc-Con-Theorem-Poincare}, page 817,  implies that a $C^r$ vector
field has a local $C^k$ (here $r, k \geq 1$) linearization near a critical point $p$ if and
only if its time-one map $T$ has such a linearization at $p$. 

As a simple application of Theorem \ref{Main_Theorem_1}, let us consider the Shilnikov
three-dimensional saddle focus theorem as in \cite{Shilnikov-1965}. 

One has  a vector field $X$ in $R^3$ with a saddle type critical
point at $0$ such that $DX_0$ has a pair of complex conjugate eigenvalues $a \pm 
ib$ with $a < 0, b \neq 0$ and a real positive eigenvalue $c > 0$ such
that $a+c > 0$. It is also assumed that there is a
homoclinic orbit (an orbit which is forward and backward asymptotic to
$0$). Shilnikov then showed that there are infinitely many saddle
periodic orbits near the homoclinic orbit. In fact, the geometry shows that 
horseshoe type dynamics (including symbolic systems) appear. 

 There is a geometric description of this in 
section 6.5 in \cite{Guckenheimer-Holmes}. The treatment assumes that the
flow is linear in a neighborhood of the critical point. It is
not difficult to see that a $C^1$ linearization would suffice, so our
Theorem \ref{Main_Theorem_1} in the bi-circular case can be applied to
show that Shilnikov's Theorem holds even if the vector field is only
$C^{1,\alpha}$.

More recently, higher dimensional systems with homoclinic orbits involving saddle foci have
been studied by Ovsyannikov and Shilnikov in
\cite{Ovsyannikov-Shilnikov-1987} and
\cite{Ovsyannikov-Shilnikov-1992}. 
For related material, the reader is encouraged to look at the recent books
\cite{Shil-Shil-Tur-Chua-I} and \cite{Shil-Shil-Tur-Chua-II} which
contain a beautiful and fairly complete
treatment of much of the work carried out by the so-called Shilnikov
School in city of Nizhny-Novogorod (formerly Gorky). See
\cite{Shilnikov-Editorial} for a detailed description of the work of
Shilnikov and his many students and collaborators. 
Other treatments which include closely related work are in the
books of Ilyashenko and Li \cite{Ilyashenko-Li} and Guckenheimer and
Holmes \cite{Guckenheimer-Holmes}. 

Many of the proofs in works dealing with homoclinic orbits in
high dimension (e.g., greater than four) are complicated. We expect that
simplifications can be obtained using $C^{1,\beta}$ (or even $C^1$) linearizations on Lyapunov center manifolds 
as follows. 

Consider a $C^r$ vector field $X$ on the
Euclidean space $R^N$ with a hyperbolic critical point of saddle type
at $0$, with $r > 1$. Thus, the real parts of the eigenvalues of $DX_0$ are
non-zero and intersect both the positive and negative sets of reals.

Let $L  = DX_0$, and let 

\[ a_m > a_{m-1} >  \ldots > a_1 > 0 > b_1 > b_2 > \ldots > b_n \]

denote the distinct real parts of the eigenvalues of $L$. 

We express $R^N$ in the direct sum
decomposition 

\[ R^N = E_1 \oplus E_2 \oplus E_3 \oplus E_4 \]

such that $L(E_i) = E_i$ for $i=1,2,3,4$, and, writing $L_i$ for the
restriction $L \mid E_i$, we have 

\be
 \item the eigenvalues of $L_1$ have real parts equal to $a_1$, 
 \item the eigenvalues of $L_2$ have real parts equal to $b_1$, 
 \item the eigenvalues of $L_3$ have real parts greater than $a_1$,
and 
 \item the eigenvalues of $L_4$ have real parts less than $b_1$.
\ee

Following terminology in \cite{Shil-Shil-Tur-Chua-I} we
call the eigenvalues with real parts $a_1$ or $b_1$ the {\it leading
eigenvalues} and the corresponding subspaces $E_1, E_2$ the {\it
  leading eigenspaces}.

The numbers $a_1$ and $b_1$ are the Lyapunov exponents of $L_1$ and $L_2$, respectively.

That is, for  each $v \in E_1 \setminus \{0\}$ and $w \in E_2 \setminus \{0\}$, we have 

\[ \chi(v) = \lim_{t \rightarrow \pm \infty} \frac{1}{t} \log \n{e^{tL} v} = a_1 \]

and

\[ \chi(w) = \lim_{t \rightarrow \pm \infty} \frac{1}{t} \log \n{e^{tL} w} = b_1. \]

Accordingly, we will call the direct sum $E_1 \oplus E_2$ the {\it Lyapunov Center Subspace}, and we denote it by $E^{cc}$.  

The theory of invariant manifolds, e.g. as in \cite{Llave-Wayne},
\cite{Hirsch-Pugh-Shub},   shows that there are submanifolds 
$W^{cc}$, invariant by the local flow of $X$, tangent at $0$ to
the subspace $E^{cc}$.  We will call these {\it Lyapunov Center Manifolds}.  They are not unique, but each will be 
$C^{1,\alpha}$ near $0$ (i.e., the transition maps on local coordinate charts are
$C^{1,\alpha}$) provided that

\beq \label{C1-alpha-cond}
(1+\alpha)a_1 < a_2 \mbox{ and } (1+\alpha)b_1 > b_2.
\eeq

The derivative $DX_0$, restricted to $E_1 \times E_2$
is bi-circular, so by Theorem
\ref{Main_Theorem_1} and the Sternberg technique mentioned above, the flow of $X$ restricted to each such $W^{cc}$ is $C^{1,\beta}$ linearizable at
$0$.  Moreover, by Theorem \ref{par-Main-Theorem}, the linearizations can be chosen to  depend continuously 
on external parameters.

Note that even for linear differential equations,  most of the manifolds $W^{cc}$ may not be smoother than $C^{1,\alpha}$ for some $0 < \alpha < 1$. 
For instance, consider the three dimensional linear system $\dot{x} =
x, \ \dot{y} = (1+\alpha) y, \ \dot{z} = -z$ for $0 < \alpha < 1$. Each surface 

\[ y = C \n{x}^{1 + \alpha} \]

  for some non-zero constant $C$ can be taken as one of the manifolds $W^{cc}$.


\edoc